\documentclass{amsart}
\usepackage[ margin=1.2in]{geometry}

\usepackage{lineno,hyperref}
\usepackage{amssymb,amsmath,amsfonts,amsthm}
\usepackage{caption}
\modulolinenumbers[5]
\usepackage{subcaption}

\newtheorem{theorem}{Theorem}[section]
\newtheorem{corollary}[theorem]{Corollary}
\newtheorem{lemma}[theorem]{Lemma}
\newtheorem{definition}[theorem]{Definition}

\theoremstyle{remark}
\newtheorem{remark}[theorem]{Remark}
\newtheorem{example}[theorem]{Example}

\usepackage{color}
\usepackage{hyperref}

\usepackage{bigints}

\usepackage{overpic}
\usepackage{graphicx}
\graphicspath{{Figures/}}
\usepackage{array}

\usepackage{amssymb}
\usepackage{tikz}

\usepackage{comment}

\usetikzlibrary{matrix,arrows.meta}
\newlength{\ml}
\renewcommand{\epsilon}{\varepsilon}
\renewcommand{\phi}{\varphi}
\def\R{\mathbb R}

\newcommand{\sign}{\operatorname{sign}}

\newcommand{\ZZ}{\mathbb Z}

\newcommand{\Span}{\mathop{\rm span}}

\newcommand{\orangle}{\overrightarrow{\angle}}
\newcommand{\pr}{\mathrm{pr}}
\newcommand{\const}{\mathrm{const}}

\numberwithin{equation}{section}

\begin{document}


\title{Isometric Deformations of Discrete and Smooth T-surfaces}


\author{Ivan Izmestiev }
\author{Arvin Rasoulzadeh }
\author{Jonas Tervooren }

\address{Institute of Discrete Mathematics and Geometry\\
Vienna University of Technology \\
Wiedner Hauptstrasse 8--10\\
1040 Vienna, AUSTRIA}
\email{izmestiev@dmg.tuwien.ac.at.}
\email{rasoulzadeh@geometrie.tuwien.ac.at.}
\email{jtervooren@geometrie.tuwien.ac.at.}





\begin{abstract}
Quad-surfaces are polyhedral surfaces with quadrilateral faces and the combinatorics of the square grid.
A generic quad-surface is rigid.
T-hedra is a class of flexible quad-surfaces introduced by Graf and Sauer in 1931.
Particular examples of T-hedra are the celebrated Miura fold, discrete surfaces of revolution, and discrete molding surfaces.
We provide an explicit parametrization of the isometric deformation of a T-hedron.
T-hedra have a smooth analog, T-surfaces.
For these we provide a synthetic and an analytic description, both similar to the corresponding descriptions of T-hedra.
We also parametrize the isometric deformations of T-surfaces and discuss their deformability range.

\vspace{.5cm}
\end{abstract}

\maketitle

\section{Introduction}
\begin{figure}[ht]
\begin{center}
\includegraphics[width=.28\textwidth]{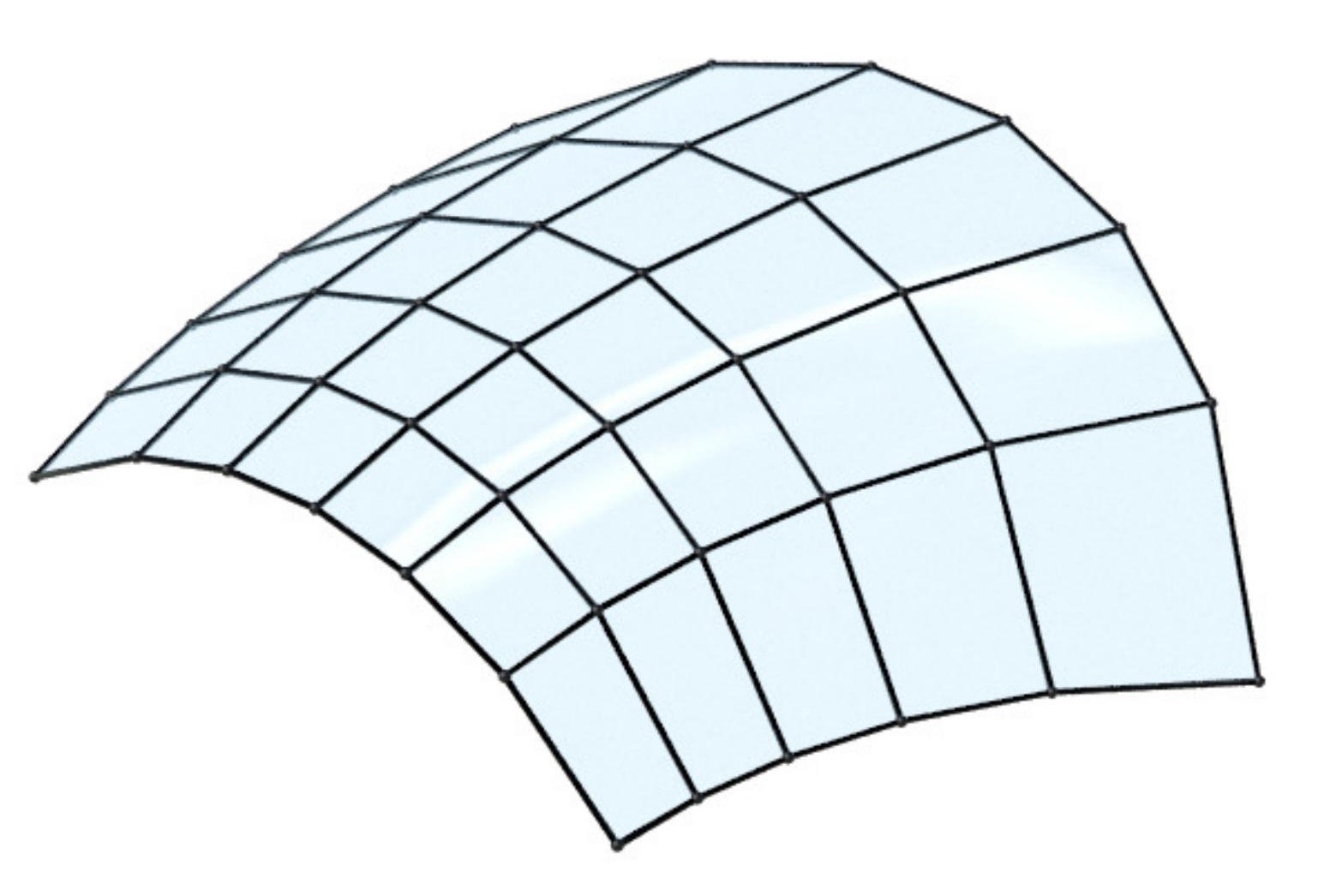}
\hspace{.04\textwidth}
\includegraphics[width=.28\textwidth]{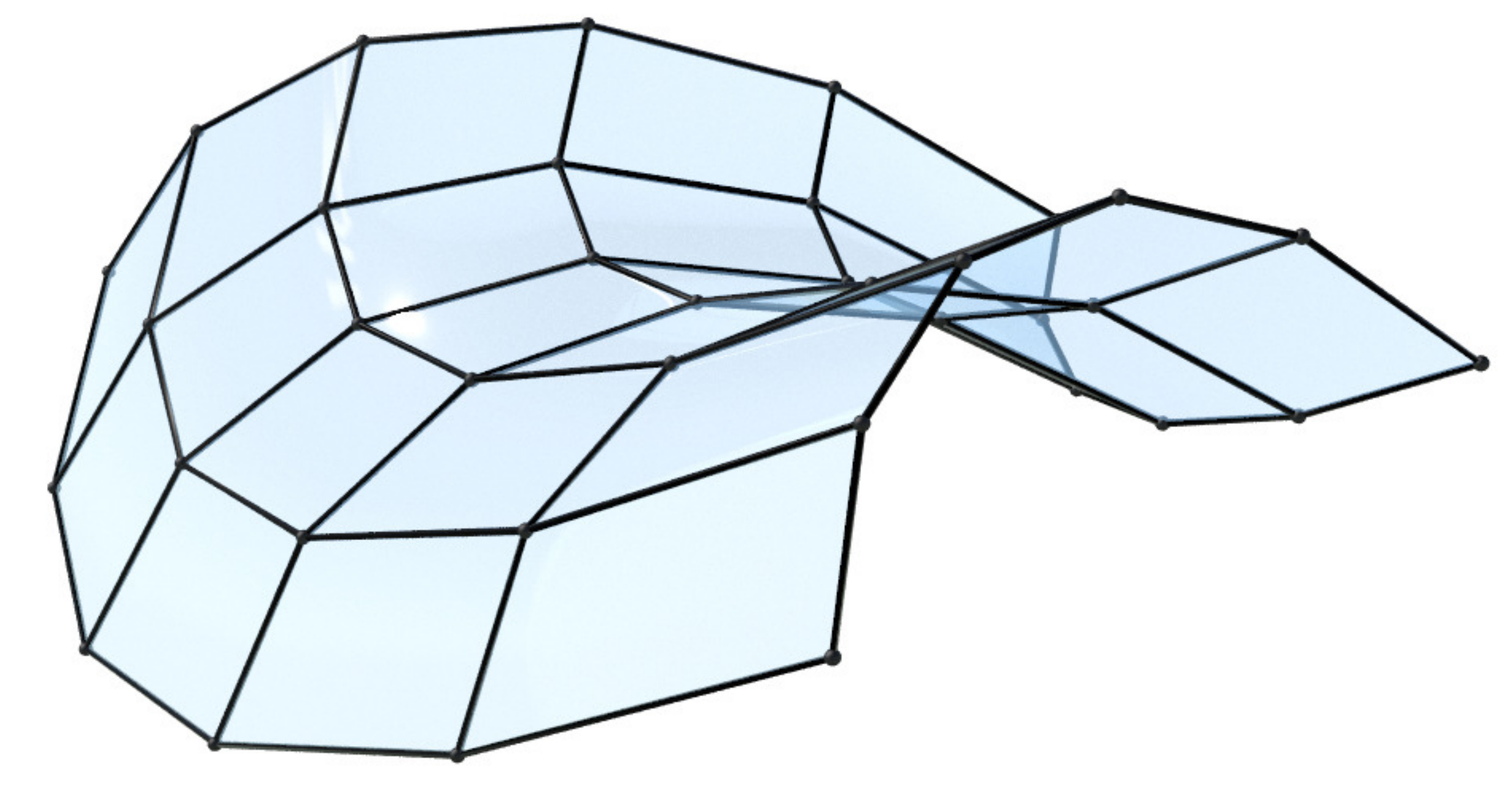}
\hspace{.04\textwidth}
\includegraphics[width=.28\textwidth]{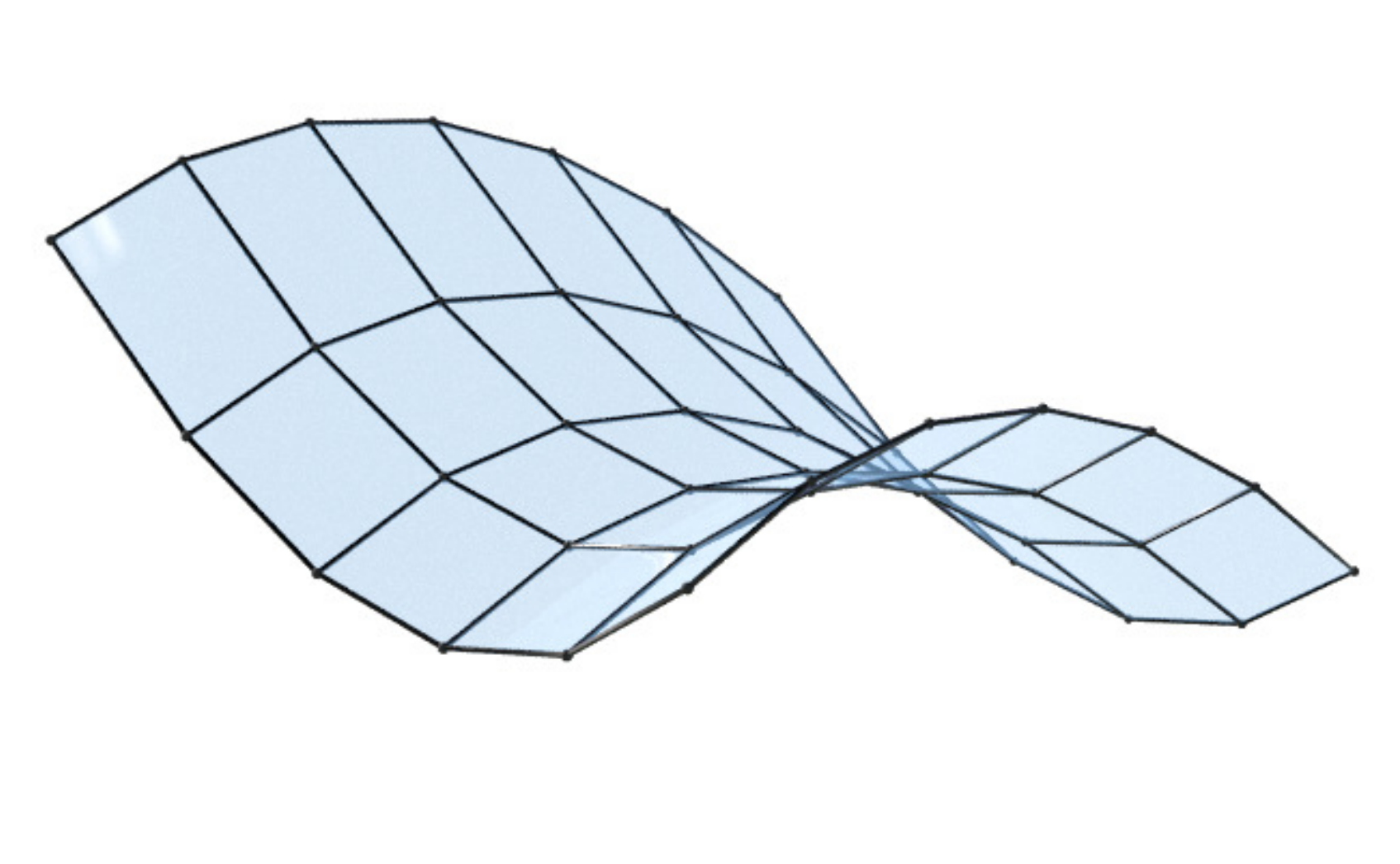}
\end{center}
\caption{Examples of quad-surfaces.}
\label{fig:QuadSurface}
\end{figure}

A \emph{quad-surface} is a polyhedral surface with quadrilateral faces connected in the combinatorics of the square grid, see Figure \ref{fig:QuadSurface} for some examples.
There are less restrictive versions of this definition which allow curved quadrilaterals or vertices where a different number of quadrilaterals meet, but for the purposes of this article we adopt the above definition.
Quad-surfaces are popular in the freeform architecture for several reasons.
For example, the fact that the average vertex degree of a quad surface is four, and not six as for a triangulated surface, simplifies the manufacturing of beam connections.
Design of quad-surfaces of desired shape and local properties leads to interesting and deep mathematical questions.
We refer the reader to the PhD thesis of Davide Pellis \cite{pellis2019phd} for a survey of main aspects and literature.

A generic quad-surface is rigid, which means that it cannot be deformed in such a way that every face moves as a rigid body while the dihedral angles between the faces change.
By contrary, a triangulated surface is usually flexible with the number of degrees of freedom approximately equal to the number of the boundary vertices.
Flexible quad-surfaces are rare, deform usually with one degree of freedom, and are an attractive object for the transformable design in architecture and engineering.
A well-known example is the Miura-ori or Miura fold shown in Figure \ref{fig:MiuraOri} which was used in the design of solar panels in the Japanese space program \cite{Miura}.
In the last decades several modifications of the Miura fold were proposed, \cite{Tachi, LangHowell, Feng}.

The study of flexible quad-surfaces was initiated almost a century ago by Graf and Sauer \cite{sauer1931flachenverbiegung} who described two classes of flexible quad-surfaces.
Miura fold belongs to both of these classes.
At about the same time Kokotsakis \cite{kokotsakis1933bewegliche} studied the flexibility and the infinitesimal flexibility of a special class of polyhedra which includes $3 \times 3$ quad-surfaces and found for the latter a series of flexible configurations different from those of Graf and Sauer.
The importance of study of flexible $3 \times 3$ quad-surfaces lies in the fact that an $m \times n$ quad-surface is flexible if and only if each of its $3 \times 3$ subsurfaces is flexible.
The description of all flexible $3 \times 3$ quad-surfaces remained for a long time an open problem.
A novel approach to this problem was proposed by Stachel in \cite{stachel2010kinematic} and further developed by Nawratil \cite{Naw11,Naw12}.
A complete solution was given by the first-named author in \cite{Izm_Koko}.
It includes many special cases and opens a vast perspective for the use of flexible quad-surfaces, but a lot of work remains to be done in order to make the abstract results of \cite{Izm_Koko} accessible to design. 

The present article makes a first step in this direction by parametrizing deformations of T-hedra, one of the flexible classes of quad-surfaces discovered by Graf and Sauer \cite{sauer1931flachenverbiegung,sauer1970differenzengeometrie}.
One of the properties of a T-hedron is that all of its faces are trapezoids, hence the letter T in the name.
Graf and Sauer define T-hedra constructively, through a certain series of parallel projections of an open planar polygon.
We start with a descriptive definition and show it to be equivalent to the original one.
The flexibility of T-hedra was proved in \cite{sauer1931flachenverbiegung} with an elegant geometric argument.
Building upon this we derive explicit formulas for the vertex coordinates of a T-hedron in dependence on the deformation parameter.
A first version of these formulas was used by Kiumars Sharifmoghaddam in a Rhino/Grasshopper plugin for practical design of T-hedra and visualization of their deformation, see \cite{sharifmoghaddam2021using}.
Formulas in the present article differ from those in \cite{sharifmoghaddam2021using} because we use a somewhat different set of data to compute the vertex coordinates.

Polyhedral surfaces can be viewed as discrete analoga of smooth surfaces, and quad-surfaces represent a very special discretization.
The grid lines of a quad-surface play a role similar to that of coordinate curves on a parametrized smooth surface.
The quads being flat is then a discrete analog of the coordinate net being conjugate.
This analogy goes back at least to Peterson \cite[\S 27]{Peterson1868} and was further developed by Graf and Sauer \cite{sauer1931flachenverbiegung, sauer1970differenzengeometrie}.
The smooth analoga of T-hedra, which we call T-surfaces were also mentioned by Graf and Sauer.
In a complete analogy to the discrete situation, a T-surface can be isometrically deformed in such a way that a designated conjugate net remains conjugate during the deformation.
In this article we give a synthetic and an analytic description of T-surfaces which are very similar to the corresponding descriptions of T-hedra.
We also derive formulas for isometric deformations of T-hedra.
The problem of isometric deformations of smooth surfaces is an old problem in differential geometry, going back to the mid-19th century, and most of our deformation results are not new.
However we find it useful to have a general treatment of all T-surfaces and instructive to draw a direct comparison between the polyhedral and the smooth situation.

\section{T-hedra}
\subsection{T-hedra and planar T-nets}
\label{sec:ThedraTnets}

\begin{definition}
\label{dfn:THedron}
A \emph{T-hedron} is a quad-surface with vertices given by a map
\[
\sigma \colon I \times J \to \R^3, \quad I = \{0, 1, \ldots, m\}, \quad J = \{0, 1, \ldots, n\}, \quad \sigma(i,j) =: \sigma_{ij}
\]
that satisfies the following conditions.
\begin{enumerate}
\item
For every $i \in I \setminus \{0\}$ and $j \in J \setminus \{0\}$ the quadrilateral $\sigma_{i-1,j-1} \sigma_{i-1,j} \sigma_{ij} \sigma_{i,j-1}$ is planar and non-self-intersecting.
\item
Each of the coordinate polygons
\[
\sigma_{\bullet j} := (\sigma_{0j}, \sigma_{1j}, \ldots, \sigma_{mj}), \quad
\sigma_{i \bullet} := (\sigma_{i0}, \sigma_{i1}, \ldots, \sigma_{in})
\]
is planar but not contained in a line.
\item
Denote by $P_j$ the plane spanned by the coordinate polygon $\sigma_{\bullet j}$ and by $Q_i$ the plane spanned by the coordinate polygon $\sigma_{i \bullet}$.
Then for all $i \in I$ and $j \in J$ the planes $P_j$ and $Q_i$ are orthogonal.
\item
Consecutive planes in the same family are different: $P_j \ne P_{j+1}$, $Q_i \ne Q_{i+1}$.
\end{enumerate}
\end{definition}

The above definition does not exclude self-intersections of the polyhedral surface arising from $\sigma$.
It would be natural to impose this restriction, but self-intersections are hard to detect and hard to avoid during an isometric deformation.

\begin{lemma}
\label{lem:OrthPlanesDiscr}
Let $\{P_j\}_{j\in J}$ and $\{Q_i\}_{i\in I}$ be two families of planes such that $P_j \perp Q_i$ for all $i$ and $j$.
Then in at least one of the families all planes are parallel to each other.
\end{lemma}
\begin{proof}
Assume that there is a pair of non-parallel planes $P_j, P_k$.
Then each of the planes $Q_i$ is orthogonal to the intersection line of $P_j$ and $P_k$, thus all planes in the family $\{Q_i\}_{i \in I}$ are parallel to each other.
\end{proof}

This lemma implies that the faces of a T-surface are trapezoids, which explains the letter T in the name of this class of quad-surfaces.

\begin{figure}[ht]
\begin{center}
\includegraphics[width=.7\textwidth]{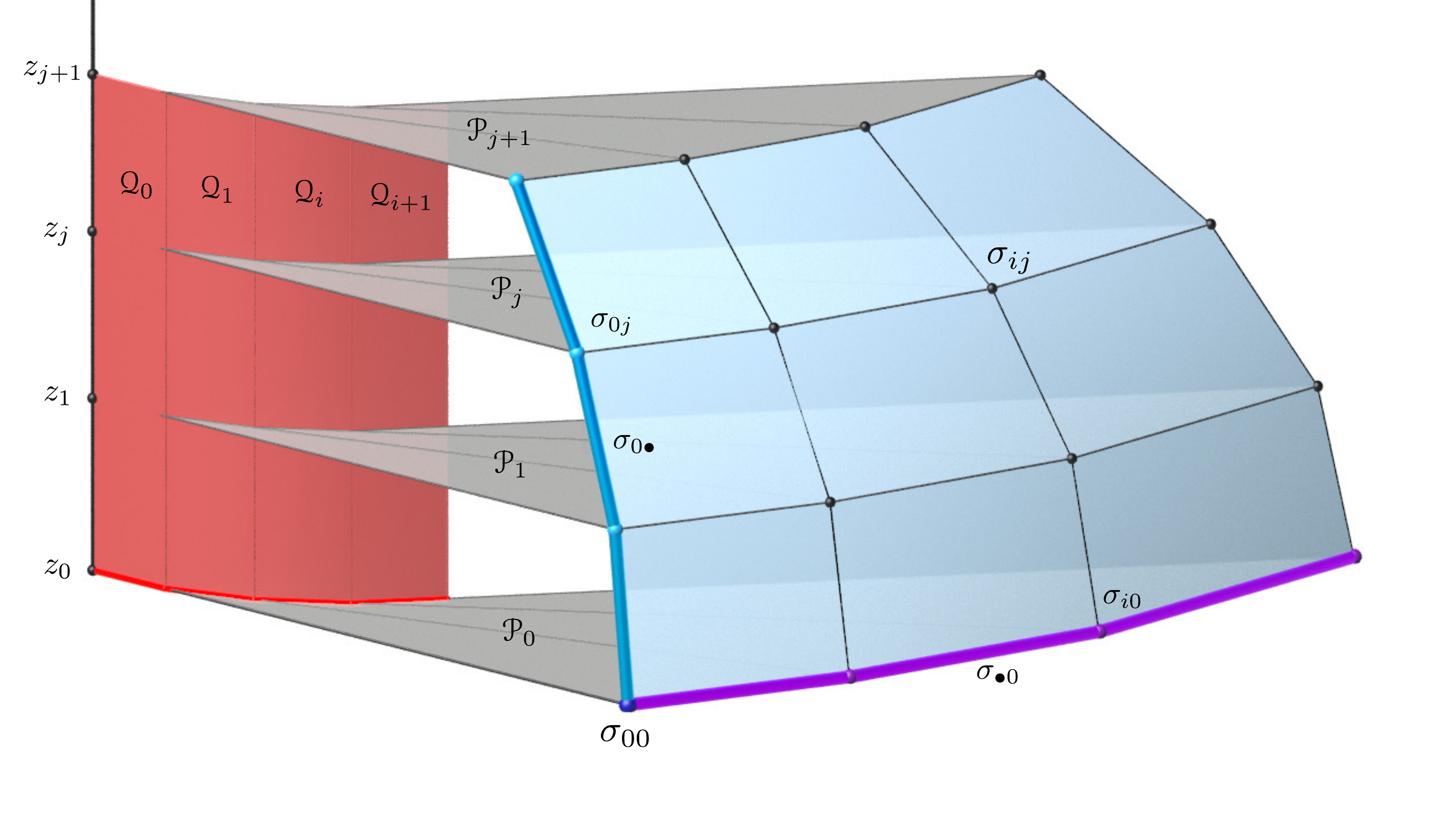}
\end{center}
\caption{Trajectory planes $P_j$, profile planes $Q_i$, a trajectory polygon $\sigma_{\bullet 0}$ and a profile polygon $\sigma_{0\bullet}$.}
\label{fig:TrajProPlanes}
\end{figure}

Without loss of generality assume that the planes $P_j$ are parallel to each other.
In this case $P_j$ are called \emph{trajectory planes} and $Q_i$ are called \emph{profile planes} of the T-surface.
The polygons $\sigma_{i\bullet} \subset Q_i$ are called the \emph{profile polygons}, and the polygons $\sigma_{\bullet j} \subset P_j$ the \emph{trajectory polygons} of $\sigma$.
See Figure \ref{fig:TrajProPlanes}.

Denote by $\pr \colon \R^3 \to P_0$ the orthogonal projection to the plane $P_0$.
The composition $\tau = \pr \circ \sigma \colon I \times J \to \R^2$ will be called the \emph{ground view} of $\sigma$.
Without loss of generality one can assume that the plane $P_0$ is the $xy$-coordinate plane, so that if $\sigma_{ij} = (x_{ij}, y_{ij}, z_{ij})$, then $\tau_{ij} = (x_{ij}, y_{ij})$.
Since for a fixed $j$ all $\sigma_{ij}$ belong to a horizontal plane $P_j$, the value of $z_{ij}$ depends on $j$ only, and we change the notation to $z_j := z_{ij}$.
The T-hedron $\sigma$ is thus uniquely determined by its ground view and a collection of numbers
\[
z_0 = 0, z_1, \ldots, z_n, \quad z_{j-1} \ne z_j \text{ for all }j \in \{1, \ldots, n\}.
\]
The condition $z_{j-1} \ne z_j$ is equivalent to the condition $P_{j-1} \ne P_j$ in Definition \ref{dfn:THedron}.

\begin{figure}[ht]
\begin{center}
\includegraphics[width=.6\textwidth]{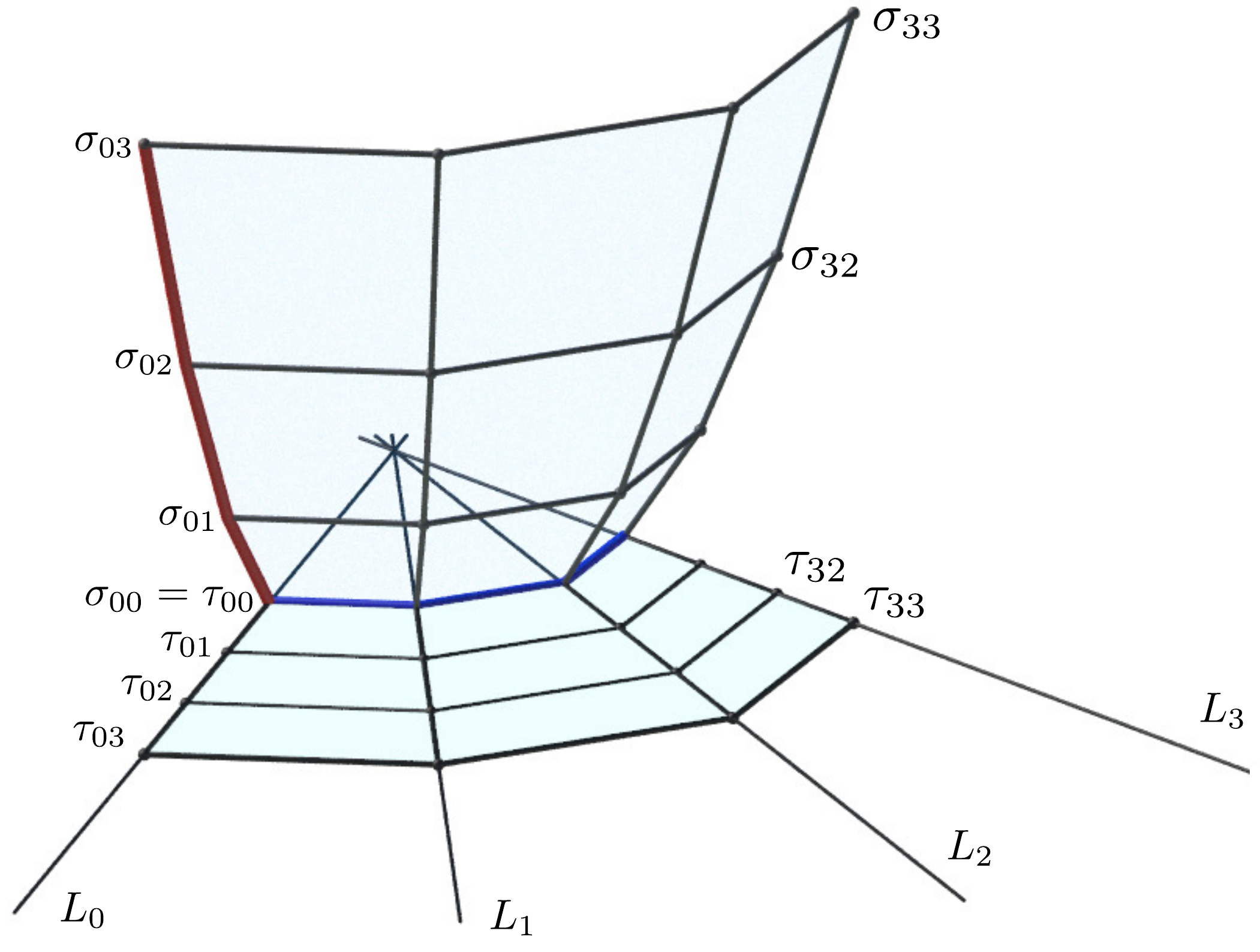}
\end{center}
\caption{A T-hedron and its ground view. }
\label{fig:THedronTNet}
\end{figure}

Since every profile plane $Q_i$ is parallel to the $z$-axis, it projects to a line $L_i \in \R^2$.
Every face of $\sigma$ projects to a trapezoid with legs on two consecutive lines $L_{i-1}$ and $L_i$.
A degeneration is possible: if a face is parallel to the $z$-axis, then its projection is a segment with endpoints on $L_{i-1}$ and $L_i$.
Figure \ref{fig:THedronTNet} shows a T-hedron together with its ground view.

\begin{definition}
A \emph{(planar) T-net} is a map
\[
\tau \colon I \times J \to \R^2, \quad I = \{0, 1, \ldots, m\}, \quad J = \{0, 1, \ldots, n\}, \quad \tau(i,j) =: \tau_{ij}
\]
that satisfies the following conditions.
\begin{enumerate}
\item
For every $i \in I$ the points $\tau_{i0}, \tau_{i1}, \ldots, \tau_{in}$ span a line, which we call a \emph{profile line} and denote by $L_i$.
\item
For every $j \in J$ the points $\tau_{0j}, \tau_{1j}, \ldots, \tau_{mj}$ are not collinear.
\item
For every $i \in I \setminus\{0\}$ the lines $L_{i-1}$ and $L_i$ are different.
\item
Every quadrilateral $\tau_{i-1,j-1}\tau_{i-1,j}\tau_{ij}\tau_{i,j-1}$ is a (non-self-intersecting) trapezoid with legs on the lines $L_{i-1}$ and $L_i$.
This trapezoid may degenerate to a segment with endpoints on $L_{i-1}$ and $L_i$.
\end{enumerate}
\end{definition}

The following lemma is straightforward, see the discussion preceding the above definition.
\begin{lemma}
\label{lem:TNetTHedron}
The ground view of every T-hedron is a T-net.
Every T-net $\tau$ and a generic collection of real numbers $0 = z_0 \ne z_1 \ne \cdots \ne z_n$ determine a T-hedron $\sigma_{ij} = (\tau_{ij}, z_j)$ whose ground view is $\tau$.
\end{lemma}
The genericity condition on $z_j$ is that the profile polygons $(\tau_{ij}, z_j)_{j=0}^n$ are not contained in a line.

\subsection{Analytic description of T-hedra and T-nets}
\label{sec:AnDescr}
\begin{figure}[ht]
\begin{picture}(0,0)%
\includegraphics{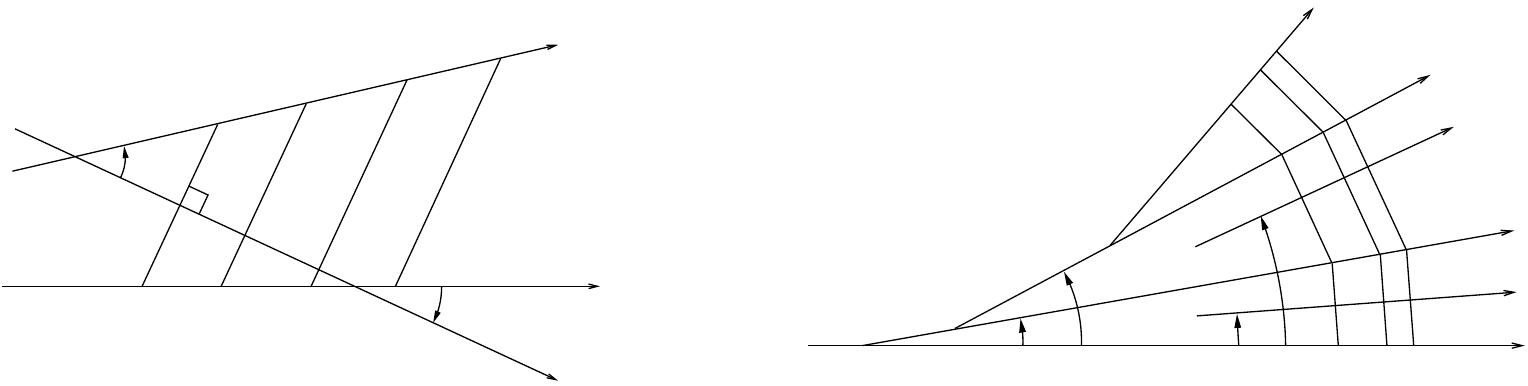}%
\end{picture}%
\setlength{\unitlength}{3315sp}%
\begingroup\makeatletter\ifx\SetFigFont\undefined%
\gdef\SetFigFont#1#2#3#4#5{%
	\reset@font\fontsize{#1}{#2pt}%
	\fontfamily{#3}\fontseries{#4}\fontshape{#5}%
	\selectfont}%
\fi\endgroup%
\begin{picture}(8729,2199)(781,-774)
\put(9454,-648){\makebox(0,0)[lb]{\smash{{\SetFigFont{8}{9.6}{\rmdefault}{\mddefault}{\updefault}{\color[rgb]{0,0,0}$L_0$}%
}}}}
\put(9383,136){\makebox(0,0)[lb]{\smash{{\SetFigFont{8}{9.6}{\rmdefault}{\mddefault}{\updefault}{\color[rgb]{0,0,0}$L_1$}%
}}}}
\put(8853,1030){\makebox(0,0)[lb]{\smash{{\SetFigFont{8}{9.6}{\rmdefault}{\mddefault}{\updefault}{\color[rgb]{0,0,0}$L_2$}%
}}}}
\put(8282,1290){\makebox(0,0)[lb]{\smash{{\SetFigFont{8}{9.6}{\rmdefault}{\mddefault}{\updefault}{\color[rgb]{0,0,0}$L_3$}%
}}}}
\put(9460,-257){\makebox(0,0)[lb]{\smash{{\SetFigFont{8}{9.6}{\rmdefault}{\mddefault}{\updefault}{\color[rgb]{0,0,0}$M_1$}%
}}}}
\put(9087,654){\makebox(0,0)[lb]{\smash{{\SetFigFont{8}{9.6}{\rmdefault}{\mddefault}{\updefault}{\color[rgb]{0,0,0}$M_2$}%
}}}}
\put(6670,-506){\makebox(0,0)[lb]{\smash{{\SetFigFont{8}{9.6}{\rmdefault}{\mddefault}{\updefault}{\color[rgb]{0,0,0}$\phi_1$}%
}}}}
\put(7895,-512){\makebox(0,0)[lb]{\smash{{\SetFigFont{8}{9.6}{\rmdefault}{\mddefault}{\updefault}{\color[rgb]{0,0,0}$\psi_1$}%
}}}}
\put(6953,-249){\makebox(0,0)[lb]{\smash{{\SetFigFont{8}{9.6}{\rmdefault}{\mddefault}{\updefault}{\color[rgb]{0,0,0}$\phi_2$}%
}}}}
\put(8075,-11){\makebox(0,0)[lb]{\smash{{\SetFigFont{8}{9.6}{\rmdefault}{\mddefault}{\updefault}{\color[rgb]{0,0,0}$\psi_2$}%
}}}}
\put(1540,424){\makebox(0,0)[lb]{\smash{{\SetFigFont{8}{9.6}{\rmdefault}{\mddefault}{\updefault}{\color[rgb]{0,0,0}$\theta_i$}%
}}}}
\put(3955,1050){\makebox(0,0)[lb]{\smash{{\SetFigFont{8}{9.6}{\rmdefault}{\mddefault}{\updefault}{\color[rgb]{0,0,0}$L_i$}%
}}}}
\put(3954,-714){\makebox(0,0)[lb]{\smash{{\SetFigFont{8}{9.6}{\rmdefault}{\mddefault}{\updefault}{\color[rgb]{0,0,0}$M_i$}%
}}}}
\put(4276,-511){\makebox(0,0)[lb]{\smash{{\SetFigFont{8}{9.6}{\rmdefault}{\mddefault}{\updefault}{\color[rgb]{0,0,0}$L_{i-1}$}%
}}}}
\put(3322,-401){\makebox(0,0)[lb]{\smash{{\SetFigFont{8}{9.6}{\rmdefault}{\mddefault}{\updefault}{\color[rgb]{0,0,0}$\eta_i$}%
}}}}
\end{picture}%
\caption{Signed angles in a T-net.}
\label{fig:SignedAngles}
\end{figure}
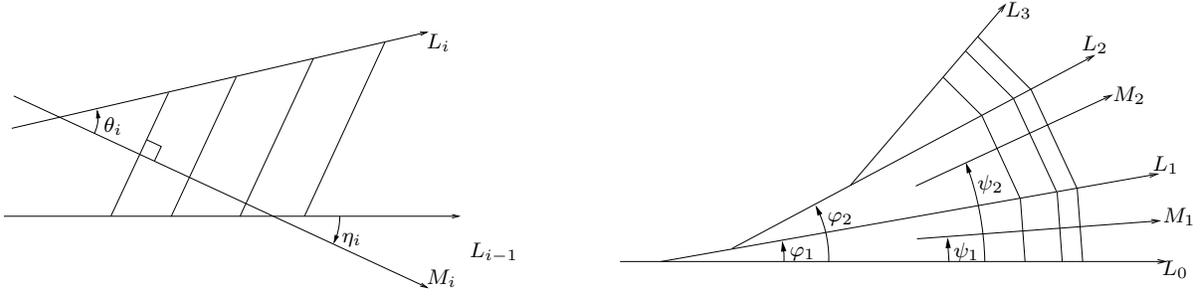

Let $\tau$ be the ground view of a T-hedron $\sigma$.
Choose an orientation of the profile line $L_0$.
It induces orientations of all other profile lines according to the rule that the legs of all trapezoids must be oriented in the same way.
For $i = 1, \ldots, m$ choose a line $M_i$ orthogonal to the bases of trapezoids with legs on the lines $L_{i-1}$ and $L_i$.
Orient $M_i$ in the way compatible with the orientations of $L_{i-1}$ and $L_i$.
Introduce for $i=1, \ldots, m$ the signed angles
\begin{equation}
\label{eqn:EtaTheta}
\begin{aligned}
\eta_i = \orangle(L_{i-1}, M_i) \in \left(-\frac{\pi}2, \frac{\pi}2\right), \quad
\theta_i = \orangle(M_i, L_i) \in \left(-\frac{\pi}2, \frac{\pi}2\right),\\
\phi_i = \orangle(L_0, L_i) = \sum_{\alpha=1}^i \eta_\alpha + \sum_{\alpha=1}^i \theta_i, \quad
\psi_i = \orangle(L_0, M_i) = \sum_{\alpha=1}^i \eta_\alpha + \sum_{\alpha=1}^{i-1} \theta_\alpha.
\end{aligned}
\end{equation}
See Figure \ref{fig:SignedAngles} for an illustration.

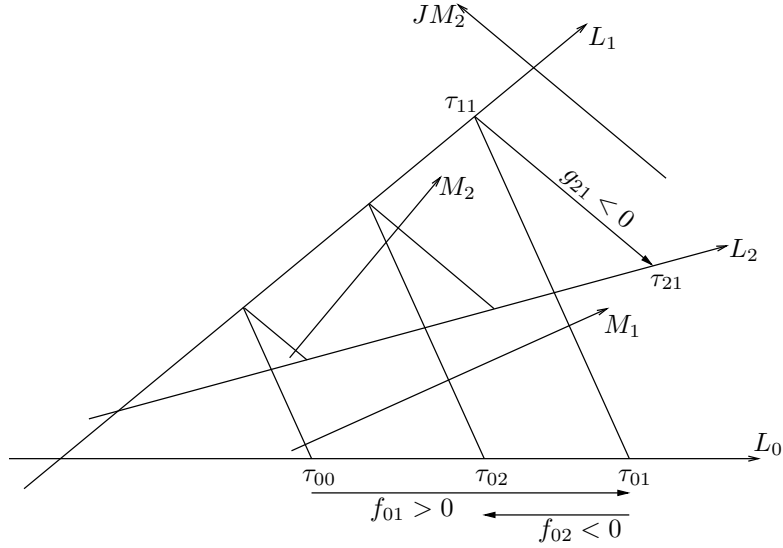
\begin{figure}[ht]
\begin{center}
\begin{picture}(0,0)%
\includegraphics{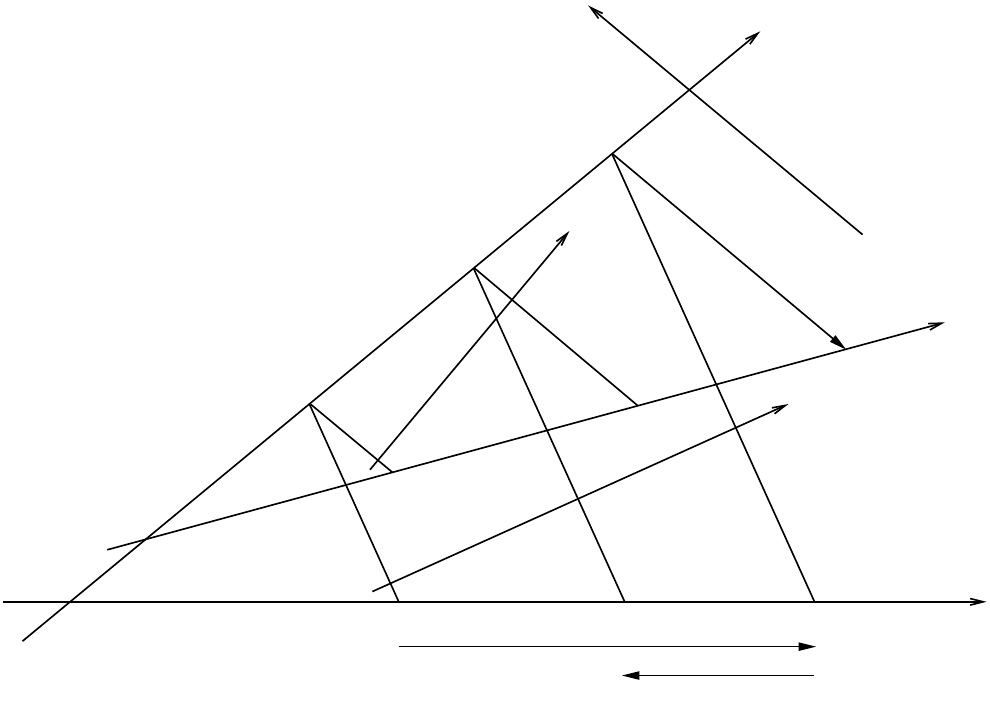}%
\end{picture}%
\setlength{\unitlength}{4144sp}%
\begingroup\makeatletter\ifx\SetFigFont\undefined%
\gdef\SetFigFont#1#2#3#4#5{%
	\reset@font\fontsize{#1}{#2pt}%
	\fontfamily{#3}\fontseries{#4}\fontshape{#5}%
	\selectfont}%
\fi\endgroup%
\begin{picture}(4524,3275)(-11,-1945)
\put(4456,-1366){\makebox(0,0)[lb]{\smash{{\SetFigFont{10}{12.0}{\rmdefault}{\mddefault}{\updefault}{\color[rgb]{0,0,0}$L_0$}%
}}}}
\put(3466,1064){\makebox(0,0)[lb]{\smash{{\SetFigFont{10}{12.0}{\rmdefault}{\mddefault}{\updefault}{\color[rgb]{0,0,0}$L_1$}%
}}}}
\put(3837,-374){\makebox(0,0)[lb]{\smash{{\SetFigFont{10}{12.0}{\rmdefault}{\mddefault}{\updefault}{\color[rgb]{0,0,0}$\tau_{21}$}%
}}}}
\put(2791,-1546){\makebox(0,0)[lb]{\smash{{\SetFigFont{10}{12.0}{\rmdefault}{\mddefault}{\updefault}{\color[rgb]{0,0,0}$\tau_{02}$}%
}}}}
\put(1756,-1546){\makebox(0,0)[lb]{\smash{{\SetFigFont{10}{12.0}{\rmdefault}{\mddefault}{\updefault}{\color[rgb]{0,0,0}$\tau_{00}$}%
}}}}
\put(3646,-1546){\makebox(0,0)[lb]{\smash{{\SetFigFont{10}{12.0}{\rmdefault}{\mddefault}{\updefault}{\color[rgb]{0,0,0}$\tau_{01}$}%
}}}}
\put(3166,-1881){\makebox(0,0)[lb]{\smash{{\SetFigFont{10}{12.0}{\rmdefault}{\mddefault}{\updefault}{\color[rgb]{0,0,0}$f_{02}<0$}%
}}}}
\put(2160,-1756){\makebox(0,0)[lb]{\smash{{\SetFigFont{10}{12.0}{\rmdefault}{\mddefault}{\updefault}{\color[rgb]{0,0,0}$f_{01}>0$}%
}}}}
\put(2571,168){\makebox(0,0)[lb]{\smash{{\SetFigFont{10}{12.0}{\rmdefault}{\mddefault}{\updefault}{\color[rgb]{0,0,0}$M_2$}%
}}}}
\put(3309,283){\rotatebox{320.0}{\makebox(0,0)[lb]{\smash{{\SetFigFont{10}{12.0}{\rmdefault}{\mddefault}{\updefault}{\color[rgb]{0,0,0}$g_{21}<0$}%
}}}}}
\put(4310,-204){\makebox(0,0)[lb]{\smash{{\SetFigFont{10}{12.0}{\rmdefault}{\mddefault}{\updefault}{\color[rgb]{0,0,0}$L_2$}%
}}}}
\put(3560,-660){\makebox(0,0)[lb]{\smash{{\SetFigFont{10}{12.0}{\rmdefault}{\mddefault}{\updefault}{\color[rgb]{0,0,0}$M_1$}%
}}}}
\put(2607,687){\makebox(0,0)[lb]{\smash{{\SetFigFont{10}{12.0}{\rmdefault}{\mddefault}{\updefault}{\color[rgb]{0,0,0}$\tau_{11}$}%
}}}}
\put(2411,1195){\makebox(0,0)[lb]{\smash{{\SetFigFont{10}{12.0}{\rmdefault}{\mddefault}{\updefault}{\color[rgb]{0,0,0}$JM_2$}%
}}}}
\end{picture}%
\end{center}
\caption{An example of a T-net with some negative edge lengths.}
\label{fig:SignedLentghs}
\end{figure}

Let $JM_i$ be the line $M_i$ rotated by $\pi/2$.
Now every edge of the T-net is equipped with an orientation coming from $L_i$ or from $JM_i$.
Let us call this \emph{geometric orientation} of edges.
Every edge also has a \emph{combinatorial orientation}, namely from $\tau_{i-1,j}$ to $\tau_{i,j}$ and from $\tau_{i,j-1}$ to $\tau_{i,j}$.
Equip the length of each edge with a sign according to whether its geometric and combinatorial orientations coincide or differ.
This gives rise to two collections of real numbers
\begin{equation}
\label{eqn:fg}
\begin{aligned}
f_{ij} &=
\begin{cases}
|\tau_{ij} - \tau_{i,j-1}|, &\text{ if } \overrightarrow{\tau_{i,j-1}\tau_{ij}} \text{ has the direction of }L_i,\\
-|\tau_{ij} - \tau_{i,j-1}|, &\text{ if } \overrightarrow{\tau_{i,j-1}\tau_{ij}} \text{ has the direction opposite to }L_i.
\end{cases}\\
g_{ij} &=
\begin{cases}
|\tau_{ij} - \tau_{i-1,j}|, &\text{ if } \overrightarrow{\tau_{i-1,j}\tau_{ij}} \text{ has the direction of }JM_i,\\
-|\tau_{ij} - \tau_{i-1,j}|, &\text{ if } \overrightarrow{\tau_{i-1,j}\tau_{ij}} \text{ has the direction opposite to }JM_i.
\end{cases}
\end{aligned}
\end{equation}
Observe that for every $j$ all numbers $f_{0j}, \ldots, f_{mj}$ are of the same sign: this is due to the consistent choice of orientations of the lines $L_i$.
Besides, the condition that the net quads are non-self-intersecting implies that for every $i$ all numbers $g_{i0}, \ldots, g_{in}$ are of the same sign.
In Figure \ref{fig:SignedLentghs} an example of a $2\times 2$ T-net with $f_{i2} < 0$ and $g_{2j} < 0$ is shown.
Observe that although this net is self-overlapping, its lifts are non-self-intersecting.

\begin{figure}[ht]
\begin{center}
\includegraphics[width=.45\textwidth]{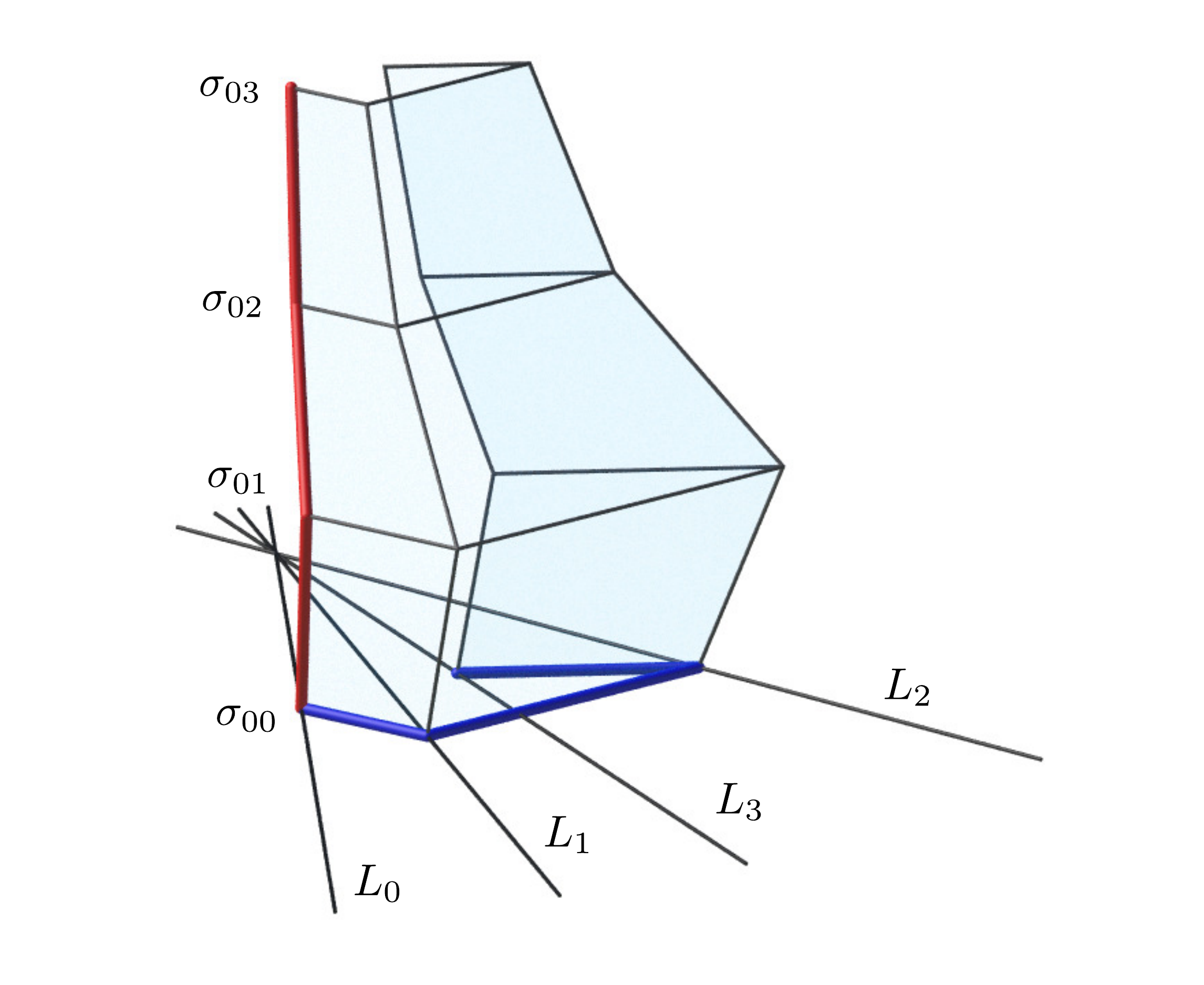}
\includegraphics[width=.45\textwidth]{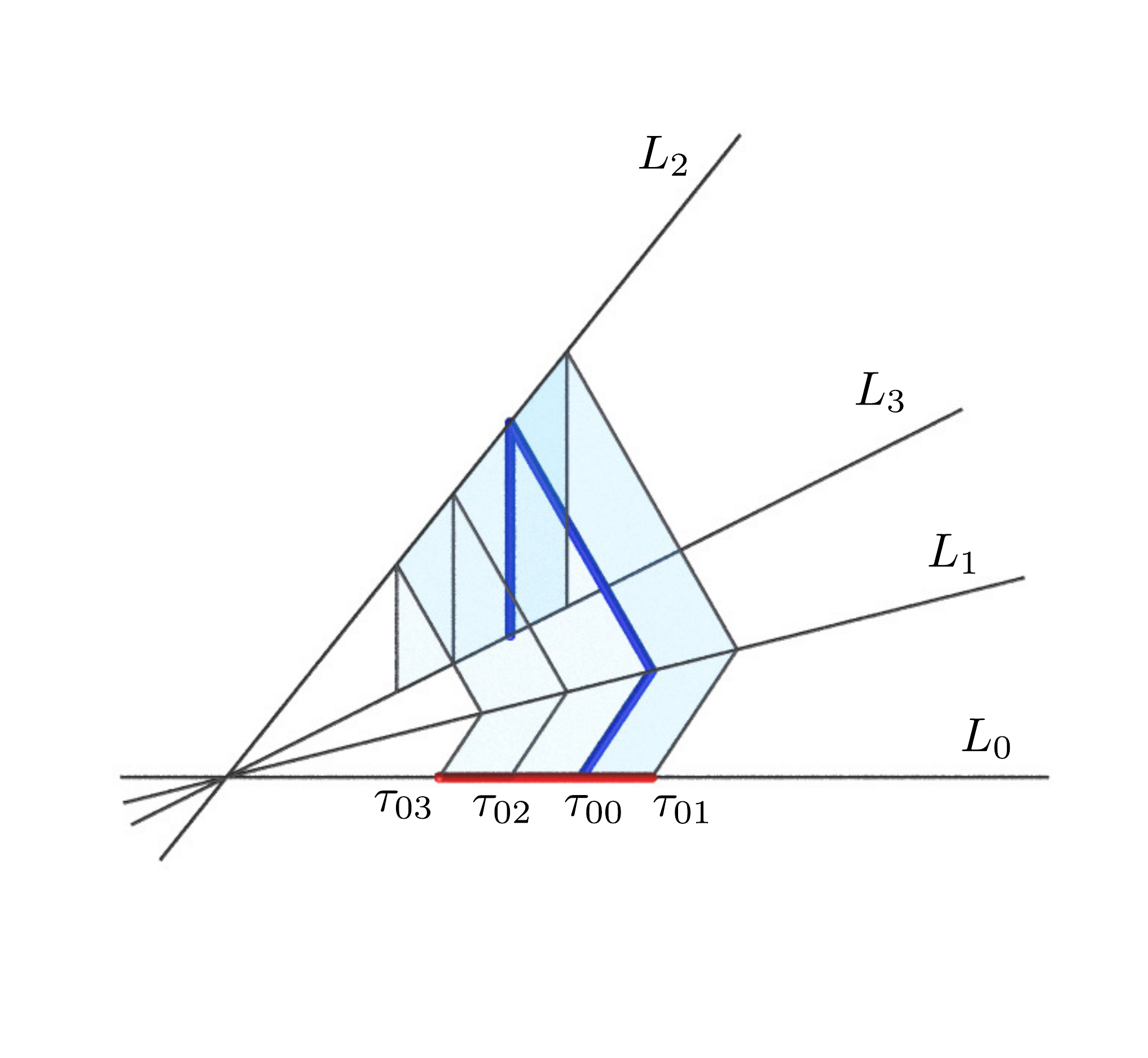}
\end{center}
\caption{The ground view of a non-self-intersecting T-hedron may be self-overlapping.}
\label{fig:NonSelfOverlap}
\end{figure}

Figure \ref{fig:NonSelfOverlap} shows another example of a non-self-intersecting T-hedron with a highly self-overlapping T-net.

\begin{theorem}
\label{thm:THedronCoord}
Let $\sigma$ be a T-hedron.
Choose a coordinate system with the origin at $\sigma_{00}$, with the $xy$-plane in the trajectory plane $P_0$, the $xz$-plane in the profile plane $Q_0$.
Then the coordinates of the vertices of $\sigma$ are given by
\begin{equation}
\label{eqn:TSurfCoord}
\sigma_{ij} = \begin{pmatrix} \tau_{ij}\\ z_j \end{pmatrix}, \quad
\tau_{ij} = \sum_{\alpha=1}^i g_{\alpha 0}
\begin{pmatrix} -\sin\psi_\alpha\\ \cos\psi_\alpha \end{pmatrix} + C_iF_j
\begin{pmatrix} \cos\phi_i\\ \sin\phi_i \end{pmatrix},
\end{equation}
where
\[
C_i = \prod_{\alpha=1}^i c_\alpha, \quad c_\alpha = \frac{\cos\eta_\alpha}{\cos\theta_\alpha}, \quad
F_j = \sum_{\beta=1}^j f_{0\beta},
\]
and $f_{0\beta}, g_{\alpha 0}, \phi_\alpha, \psi_\alpha, \eta_\alpha, \theta_\alpha$ are defined in \eqref{eqn:EtaTheta} and \eqref{eqn:fg}.
It is assumed that the line $L_0 = P_0 \cap Q_0$ is equipped with the orientation of the $x$-axis.

In particular, a T-hedron is uniquely determined by a profile polygon $\sigma_{0\bullet}$, a trajectory polygon $\sigma_{\bullet 0}$, and the directions of the profile planes $Q_1, \ldots, Q_m$.

Conversely, any set of data
\begin{gather*}
\phi_1, \ldots, \phi_m \in \R, \quad \psi_1, \ldots, \psi_m \in \R,\\
f_{01}, \ldots, f_{0n} \in \R, \quad g_{10}, \ldots, g_{m0} \in \R \setminus \{0\},\\
0 = z_0 \ne z_1 \ne \cdots \ne z_n \in \R
\end{gather*}
determines by the above formulas a T-surface, provided that the following conditions are satisfied:
\begin{enumerate}
\item
$\eta_i := \psi_i - \phi_{i-1} \in \left(-\frac{\pi}2, \frac{\pi}2 \right), 
\theta_i := \phi_i - \psi_i \in \left(-\frac{\pi}2, \frac{\pi}2 \right)$;
\item
for every $i$ the numbers $g_{i0}, \ldots, g_{im}$ computed from $\tau_{ij}$ by the formulas \eqref{eqn:fg} are of the same sign;
\item
the profile and the trajectory polygons are not contained in a line.
\end{enumerate}
\end{theorem}
\begin{proof}
Due to Lemma \ref{lem:TNetTHedron} for the first part of the theorem it suffices to show that the coordinates of the vertices of any T-net are given by the second formula of \eqref{eqn:TSurfCoord}.
One has
\[
\tau_{i0} = \sum_{\alpha=1}^i g_{\alpha 0}
\begin{pmatrix} -\sin\psi_\alpha\\ \cos\psi_\alpha \end{pmatrix}
\]
for every $i\in I$ (including $i=0$ when the sum is zero) by definition of $\psi_\alpha$ and $g_{\alpha 0}$, see \eqref{eqn:EtaTheta} and \eqref{eqn:fg}.

For every $(i,j) \in I \times J$ the point $\tau_{ij}$ lies on the line $L_i$ which goes through the point $\tau_{i0}$ and has the direction $\phi_i$.
The signed distance from $\tau_{i0}$ to $\tau_{ij}$ is $\sum_{\beta=1}^j f_{i\beta}$.
It follows that
\[
\tau_{ij} = \tau_{i0} + \sum_{\beta=1}^j f_{i\beta}
\begin{pmatrix} \cos\phi_i\\ \sin\phi_i \end{pmatrix}.
\]
Next, observe that for every $i \in I \setminus\{0\}$ and $j \in J \setminus\{0\}$ one has
\[
f_{i-1,j} \cos\eta_i = f_{ij} \cos\theta_i,
\]
since both products are equal to the (signed) height of a net trapezoid.
This implies
\[
f_{ij} = c_i f_{i-1,j} \Rightarrow f_{ij} = C_i f_{0j} \Rightarrow \sum_{\beta=1}^j f_{i\beta} = C_i F_j,
\]
and the formula \eqref{eqn:TSurfCoord} is proved.

The profile polygon $\sigma_{0\bullet}$ contains the information on $z_0, \ldots, z_n$ as well as on $F_0, \ldots, F_n$.
The trajectory polygon $\sigma_{\bullet 0}$ contains the information on $g_{10}, \ldots, g_{m0}$ and on $\psi_1, \ldots, \psi_m$.
Finally, the directions of the planes $Q_1, \ldots, Q_m$ contain the information on the angles $\phi_1, \ldots, \phi_m$.

The last statement of the theorem is straightforward from the construction.
\end{proof}

\subsection{Special classes of T-hedra}
\subsubsection{Translational T-hedra}
\label{sec:TranslTHedra}
A T-hedron is called a \emph{translational T-hedron} if all of its profile planes are parallel to each other.
That is, in the two families of planes spanned by the coordinate polygons of a translational T-hedron the planes from the same family are parallel, and the planes from different families are orthogonal.

Since all faces of translational T-hedra are parallelograms, these polyhedra are discrete translation surfaces.
Let us briefly recall the basic properties of the latter.

\begin{definition}
A quad-surface is called a \emph{discrete translation surface} if all of its faces are parallelograms.
\end{definition}
This local characterization can be globalized.
\begin{lemma}
A quad-surface $\sigma$ is a discrete translation surface if and only if one has
\[
\sigma_{ij} + \sigma_{kl} = \sigma_{il} + \sigma_{kj}
\]
for all $i, j, k, l$.
\end{lemma}
\begin{proof}
For $k=i+1$ and $l=j+1$ this equation is equivalent to all faces of $\sigma$ being parallelograms.
Thus if the equation holds, then $\sigma$ is a discrete translation surface.

For the opposite direction, assume
\[
\sigma_{ij} + \sigma_{i+1,j+1} = \sigma_{i,j+1} + \sigma_{i+1,j}
\]
for all $i, j$.
By induction on the index $k$ starting from $k=i+1$ one proves
\[
\sigma_{ij} + \sigma_{k,j+1} = \sigma_{i,j+1} + \sigma_{kj}.
\]
Then by induction on the index $l$ one obtains the statement of the lemma.
\end{proof}

A special case is the formula $\sigma_{ij} = \sigma_{i0} + \sigma_{0j} - \sigma_{00}$,
which can be interpreted in the following way.
\begin{corollary}
A discrete translation surface can be obtained from any two of its coordinate polygons $\sigma_{i\bullet}$ and $\sigma_{\bullet j}$ by translating one of them along the other.

Conversely, any two polygons with a common vertex and such that no edge of the first polygon is parallel to an edge of the second polygon generate a discrete translation surface.
\end{corollary}

The following theorem is straightforward.
\begin{theorem}
Translational T-hedra are exactly discrete translation surfaces with planar generatrices (coordinate polygons) lying in orthogonal planes.
\end{theorem}

In order to compute the vertices of a translational T-hedron choose a coordinate system with the origin at $\sigma_{00}$, with the $xy$-plane in the trajectory plane $P_0$ and the $xz$-plane in the profile plane $Q_0$.
Then due to $\sigma_{ij} = \sigma_{i0} + \sigma_{0j}$ one has
\begin{equation}
\label{eqn:DiscrTransl}
\sigma_{ij} =
\begin{pmatrix}
x_{i0} + x_{0j}\\ y_i\\ z_j
\end{pmatrix}, \quad
x_{00} = y_0 = z_0 = 0.
\end{equation}

In terms of Section \ref{sec:AnDescr} the translational T-hedra are characterized by the condition $\phi_i = 0$ for all $i$.
This implies $\theta_i = -\eta_i$ and $C_i = 1$ for all~$i$.
Equation \eqref{eqn:TSurfCoord} rewrites as
\[
\sigma_{ij} =
\begin{pmatrix}
F_j - \sum_{\alpha = 1}^i g_{\alpha 0} \sin\psi_\alpha\\
\sum_{\alpha = 1}^i g_{\alpha 0} \cos\psi_\alpha\\
z_j
\end{pmatrix},
\]
that is
\[
x_{0j} = F_j, \quad x_{i0} = - \sum_{\alpha = 1}^i g_{\alpha 0} \sin\psi_\alpha, \quad y_i = \sum_{\alpha = 1}^i g_{\alpha 0} \cos\psi_\alpha.
\]

\subsubsection{Discrete molding surfaces}
Discrete molding surfaces are T-hedra made of isosceles trapezoids.
This is equivalent to all trapezoids in the ground view being isosceles which, in turn, is equivalent to $\theta_i = \eta_i$.
This implies $\psi_i = \frac{\phi_{i-1} + \phi_i}2$ and $C_i = 1$ for all $i$ and simplifies equation \eqref{eqn:TSurfCoord} to
\begin{equation}
\label{eqn:AnDescrMoldHedra}
\sigma_{ij} = \begin{pmatrix} \tau_{ij}\\ z_j \end{pmatrix}, \quad
\tau_{ij} = \sum_{\alpha=1}^i g_{\alpha 0}
\begin{pmatrix} -\sin\frac{\phi_{\alpha-1} + \phi_\alpha}2\\ \cos\frac{\phi_{\alpha-1} + \phi_\alpha}2 \end{pmatrix} + F_j
\begin{pmatrix} \cos\phi_i\\ \sin\phi_i \end{pmatrix}.
\end{equation}

A common property of translational T-hedra and discrete molding surfaces is that their profile polygons are congruent to each other.
Similarly to translational T-hedra, a discrete molding surface is uniquely determined by one of its profile polygons and one of its trajectory polygons.

\subsubsection{Axial T-hedra}
If all profile planes intersect in a line, then we call this line the axis of the T-hedron.
In the corresponding T-net all lines $L_i$ meet at a point.

\begin{theorem}
A T-surface $\sigma$ is axial if and only if
\begin{equation}
\label{eqn:AxialCrit}
\frac{g_{i0}}{g_{10}} = \frac{\sin(\eta_i+\theta_i) \cos\theta_1}{\sin(\eta_1+\theta_1)\cos\eta_i} C_i \quad \text{for all }i.
\end{equation}
If $\sigma$ is axial, then the signed distance along the line $L_0$ from the intersection point of the axis with the ground plane $P_0$ to the point $\tau_{00}$ is equal to $\frac{g_{10}\cos\eta_1}{\sin(\eta_1+\theta_1)}$.
\end{theorem}
\begin{proof}
Denote the intersection point of the lines $L_{i-1}$ and $L_i$ by $\rho_i$ (if these lines are parallel, then $\sin(\eta_i+\theta_i) = 0$ and the equation \eqref{eqn:AxialCrit} cannot hold because $g_{i0} \ne 0$).
The lines $L_{i-1}, L_i$, and $L_{i+1}$ meet at a point if and only if $\rho_i = \rho_{i+1}$.
One computes the signed lengths of $\rho_i \tau_{i0}$ and $\rho_{i+1}\tau_{i0}$ by applying the sine law to the triangles $\rho_i \tau_{i-1,0} \tau_{i0}$ and $\rho_{i+1} \tau_{i0} \tau_{i+1,0}$ and obtains
\[
\rho_i = \rho_{i+1} \Leftrightarrow
\frac{g_{i0} \cos\eta_i}{\sin(\eta_i+\theta_i)} =
\frac{g_{i+1,0} \cos\theta_{i+1}}{\sin(\eta_{i+1}+\theta_{i+1})} \Leftrightarrow
\frac{g_{i+1,0}}{g_{i0}} =
\frac{\sin(\eta_{i+1}+\theta_{i+1}) \cos\eta_i}{\sin(\eta_i+\theta_i) \cos\theta_{i+1}}.
\]
The latter equation holds if \eqref{eqn:AxialCrit} is true.
In the opposite direction, if the latter equation holds for all $i$, then \eqref{eqn:AxialCrit} follows by induction on $i$.
\end{proof}

Choose the common point of the lines $L_i$ as the origin and extend the range of the T-net to $(\{-1\} \cup I) \times J$ by setting $\tau_{-1,j} = 0$ for all $j$.
The trapezoids in the additional trajectory strip degenerate to triangles.
The coordinates of the vertices can be obtained by appending to the sequence $f_{01}, \ldots, f_{0n}$ the distance $f_{00} = \frac{g_{10}\cos\eta_1}{\sin(\eta_1+\theta_1)}$ from the origin to the vertex $\tau_{00}$, by setting $g_{i,-1} = 0$ for all $i$ and computing the coordinates of all vertices as described in Theorem \ref{thm:THedronCoord} from the profile polygon $\tau_{0\bullet}$ and the degenerate trajectory polygon $\tau_{\bullet,-1}$:
\begin{equation}
\label{eqn:AnDescrAxisThedra}
\sigma_{ij} =
\begin{pmatrix}
C_iF_j \cos\phi_i\\ C_iF_j \sin\phi_i\\ z_j
\end{pmatrix}, \quad
C_i = \prod_{\alpha=1}^i c_\alpha, \quad c_\alpha = \frac{\cos\eta_\alpha}{\cos\theta_\alpha}, \quad
F_j = \sum_{\beta=0}^j f_{0\beta}.
\end{equation}

Axial discrete molding surfaces are \emph{discrete surfaces of revolution}.
Their profile polygons are obtained by rotating the initial profile polygon about the axis.
The vertices of a discrete surface of revolution are given by
\begin{equation}
\label{eqn:DiscrSurfRev}
\sigma_{ij} =
\begin{pmatrix} F_j \cos\phi_i\\ F_j \sin\phi_i\\ z_j \end{pmatrix}.
\end{equation}

\begin{figure}[ht]
\begin{center}
\includegraphics[width=.9\textwidth]{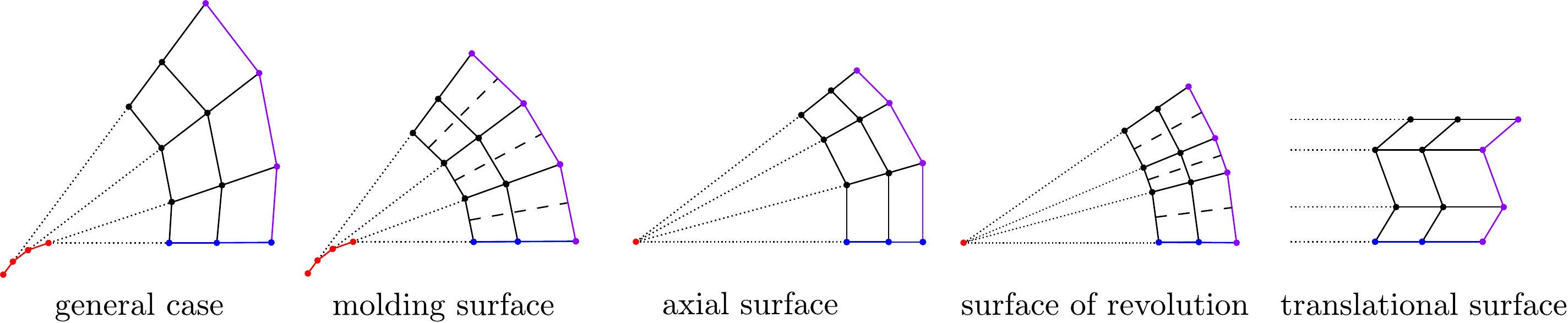}
\end{center}
\caption{The ground views of different classes of T-hedra.}
\label{fig:Classification}
\end{figure}

Figure \ref{fig:Classification} shows different classes of T-hedra in their ground views.

\section{Isometric deformations of T-hedra}
\subsection{Deformations of general T-hedra}
\begin{theorem}
\label{thm:IsomDefTHedron}
Every T-hedron $\sigma \colon I \times J \to \R^3$ allows an isometric deformation within the class of T-hedra.
That is, there is a one-parameter family of T-hedra $\sigma(t) \colon I \times J$ for $t \in (-\epsilon, \epsilon)$ such that
\begin{enumerate}
\item
$\sigma(0) = \sigma$;
\item
for every $t \in [a,b]$ and for every $i \in I \setminus \{0\}$, $j \in J \setminus \{0\}$ the quadrilateral $\sigma_{i-1,j-1}(t)\sigma_{i,j-1}(t)\sigma_{ij}(t)\sigma_{i-1,j}(t)$ is congruent to the quadrilateral $\sigma_{i-1,j-1}\sigma_{i,j-1}\sigma_{ij}\sigma_{i-1,j}$;
\item
the quad-surfaces $\sigma(t)$, $t \ne 0$ are not congruent to the quad-surface $\sigma$, that is there is no family $\Phi(t)$ of isometries of $\R^3$ such that $\sigma(t) = \Phi(t) \circ \sigma$.
\end{enumerate}

If the vertices of $\sigma$ are given by \eqref{eqn:TSurfCoord}, then the vertices of $\sigma(t)$ are as follows:
\begin{equation}
\label{eqn:THedronDeform}
\sigma_{ij}(t) =
\begin{pmatrix}
\tau_{ij}(t)\\ z_{j}(t)
\end{pmatrix}, \quad
\tau_{ij}(t) = \sum_{\alpha=1}^i g_{\alpha 0}
\begin{pmatrix}
-\sin \psi_{\alpha}(t)\\ \cos \psi_\alpha (t)
\end{pmatrix}
+ C_i(t) F_j
\begin{pmatrix}
\cos\phi_i(t)\\ \sin\phi_i(t)
\end{pmatrix},
\end{equation}
where
\begin{gather*}
\psi_i(t) = \sum_{\alpha=1}^i \eta_\alpha(t) + \sum_{\alpha=1}^{i-1} \theta_\alpha(t), \quad \phi_i(t) = \sum_{\alpha=1}^i \eta_\alpha(t) + \sum_{\alpha=1}^i \theta_\alpha(t),\\
\sin\eta_i(t) = \frac{C_{i-1}}{C_{i-1}(t)}\sin\eta_i, \quad \sin\theta_i(t) = \frac{C_i}{C_i(t)}\sin\theta_i, \quad
C_i(t) = \sqrt{C_i^2 + t},\\
z_{j}(t) = \sum_{\beta=1}^j \sign(\Delta_\beta z) \sqrt{(\Delta_\beta z)^2 - t (\Delta_\beta F)^2}.
\end{gather*}
\end{theorem}

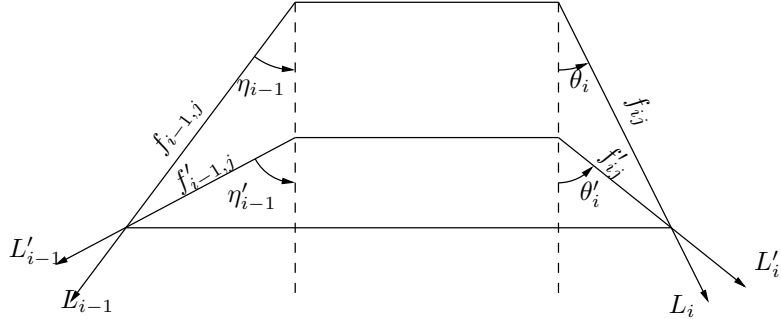
\begin{figure}[ht]
\begin{center}
\begin{picture}(0,0)%
\includegraphics{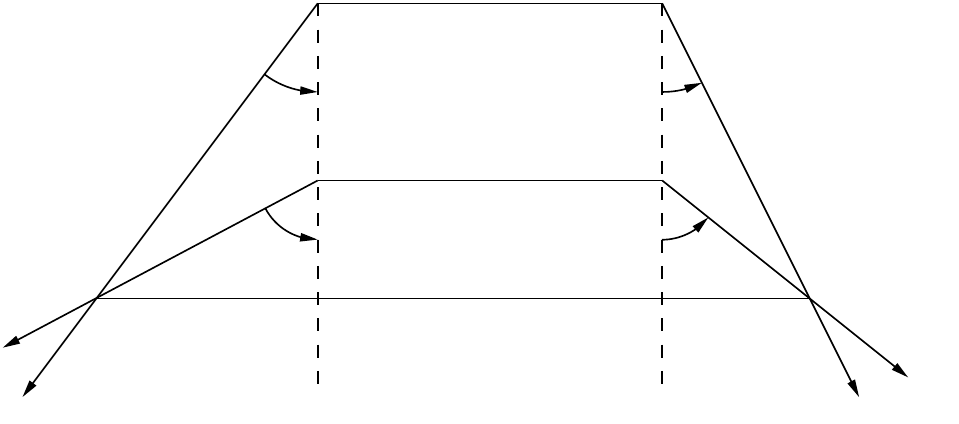}%
\end{picture}%
\setlength{\unitlength}{4144sp}%
\begingroup\makeatletter\ifx\SetFigFont\undefined%
\gdef\SetFigFont#1#2#3#4#5{%
	\reset@font\fontsize{#1}{#2pt}%
	\fontfamily{#3}\fontseries{#4}\fontshape{#5}%
	\selectfont}%
\fi\endgroup%
\begin{picture}(4374,1933)(-5015,268)
\put(-1168,328){\makebox(0,0)[rb]{\smash{{\SetFigFont{10}{12.0}{\rmdefault}{\mddefault}{\updefault}{\color[rgb]{0,0,0}$L_i$}%
}}}}
\put(-656,585){\makebox(0,0)[rb]{\smash{{\SetFigFont{10}{12.0}{\rmdefault}{\mddefault}{\updefault}{\color[rgb]{0,0,0}$L'_i$}%
}}}}
\put(-4644,359){\makebox(0,0)[rb]{\smash{{\SetFigFont{10}{12.0}{\rmdefault}{\mddefault}{\updefault}{\color[rgb]{0,0,0}$L_{i-1}$}%
}}}}
\put(-4954,664){\makebox(0,0)[rb]{\smash{{\SetFigFont{10}{12.0}{\rmdefault}{\mddefault}{\updefault}{\color[rgb]{0,0,0}$L'_{i-1}$}%
}}}}
\put(-1795,1674){\makebox(0,0)[rb]{\smash{{\SetFigFont{10}{12.0}{\rmdefault}{\mddefault}{\updefault}{\color[rgb]{0,0,0}$\theta_i$}%
}}}}
\put(-3620,1656){\makebox(0,0)[rb]{\smash{{\SetFigFont{10}{12.0}{\rmdefault}{\mddefault}{\updefault}{\color[rgb]{0,0,0}$\eta_{i-1}$}%
}}}}
\put(-1733,1003){\makebox(0,0)[rb]{\smash{{\SetFigFont{10}{12.0}{\rmdefault}{\mddefault}{\updefault}{\color[rgb]{0,0,0}$\theta'_i$}%
}}}}
\put(-3674,991){\makebox(0,0)[rb]{\smash{{\SetFigFont{10}{12.0}{\rmdefault}{\mddefault}{\updefault}{\color[rgb]{0,0,0}$\eta'_{i-1}$}%
}}}}
\put(-1601,1149){\rotatebox{320.0}{\makebox(0,0)[rb]{\smash{{\SetFigFont{10}{12.0}{\rmdefault}{\mddefault}{\updefault}{\color[rgb]{0,0,0}$f'_{ij}$}%
}}}}}
\put(-4130,1587){\rotatebox{54.0}{\makebox(0,0)[rb]{\smash{{\SetFigFont{10}{12.0}{\rmdefault}{\mddefault}{\updefault}{\color[rgb]{0,0,0}$f_{i-1,j}$}%
}}}}}
\put(-3928,1262){\rotatebox{32.0}{\makebox(0,0)[rb]{\smash{{\SetFigFont{10}{12.0}{\rmdefault}{\mddefault}{\updefault}{\color[rgb]{0,0,0}$f'_{i-1,j}$}%
}}}}}
\put(-1488,1444){\rotatebox{298.0}{\makebox(0,0)[rb]{\smash{{\SetFigFont{10}{12.0}{\rmdefault}{\mddefault}{\updefault}{\color[rgb]{0,0,0}$f_{ij}$}%
}}}}}
\end{picture}%
\end{center}
\caption{Deformation of angles and lengths in the ground view.}
\label{fig:GroundViewDeformLengths}
\end{figure}

\begin{proof}
The proof consists of two parts.
In the first part we show that any isometric deformation of a T-hedron within the class of T-hedra is subject to a very restrictive condition.
Using this, in the second part we derive the formulas for $\sigma_{ij}(t)$.

Assume that $\sigma'$ is a T-hedron isometric to $\sigma$.
Applying, if needed, a rigid motion to $\sigma'$, we can assume that the planes $P_0$ and $P'_0$ coincide, so that the ground views $\tau$ and $\tau'$ lie in the same plane.
Every trapezoid of $\sigma$ or $\sigma'$ has its bases parallel to the ground plane.
Therefore in the ground view it gets shrinked in the direction orthogonal to the bases.
In particular, two corresponding trapezoids of $\tau$ and $\tau'$ differ by a scaling in the direction orthogonal to their bases, see Figure \ref{fig:GroundViewDeform}, left.

\begin{figure}[ht]
\begin{center}
\includegraphics{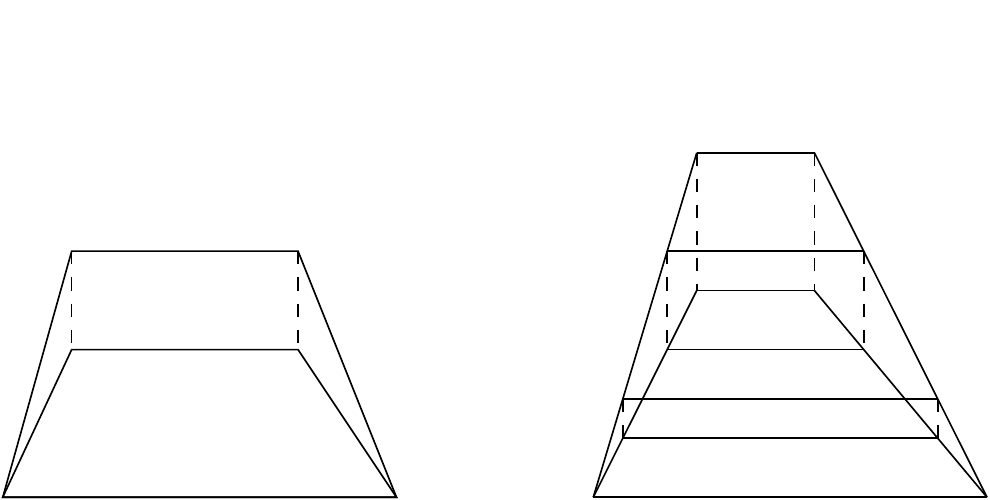}
\end{center}
\caption{The ground views of corresponding faces (left) and of corresponding profile strips (right) of two isometric T-hedra.}
\label{fig:GroundViewDeform}
\end{figure}

The $i$-th profile strips of $\tau$ and $\tau'$ consist of piles of trapezoids with their legs on a pair of lines.
It follows that all trapezoids in a profile strip undergo a scaling by the same factor, see Figure \ref{fig:GroundViewDeform}, right.
This leads to a recursive construction of $\tau'$.
Scale the first profile strip of $\tau$ in the direction of the line $M_1$ by some factor.
This scales all the edge lengths on the line $L_1$ by some factor $k_1$.
Scale the second profile strip so that the induced scaling of the line $L_1$ is also by factor $k_1$.
Then scale the third profile strip so that $L_2$ is scaled in a compatible way, and so on.
The scaled profile strips can be glued along the scaled lines $L_i$, which yields a T-net $\tau'$.

The above argument gives only a necessary condition on the ground view of an isometrically deformed T-hedron.
One still has to prove that the deformed T-net can be lifted to a T-hedron isometric to $\sigma$.
The possibility of such a lift can be proved geometrically: the scaling of profile strips can be realized by changing distances between the trajectory planes $P_j$ as long as the scaling factors do not differ too much from $1$.
Since we will need an analytic description anyway, let us turn to formulas.
Recall the notations from Sections \ref{sec:ThedraTnets} and \ref{sec:AnDescr}.
The vertical component of the $j$-th trajectory polygon was denoted by $z_j$, so that $\sigma_{ij} = (\tau_{ij}, z_j)$.
Denote
\[
\Delta_j z = z_j - z_{j-1}, \quad j = 1, \ldots, n.
\]
Then according to \eqref{eqn:fg} one has
\[
|\sigma_{ij} - \sigma_{i,j-1}| = \sqrt{f_{ij}^2 + (\Delta_j z)^2}, \quad
|\sigma_{ij} - \sigma_{i-1,j}| = |g_{ij}|.
\]
Our construction of $\tau'$ preserves the signed lengths $g_{ij}$.
Therefore it remains to find a sequence $\Delta_1 z', \ldots, \Delta_n z'$ such that
\begin{equation}
\label{eqn:Z'}
(f'_{ij})^2 + (\Delta_j z')^2 = f_{ij}^2 +(\Delta_j z)^2 \quad \text{for all } i \in I, j \in J \setminus \{0\}.
\end{equation}
Let $k_i$ be the scaling factor of $L_i$, $i \in I$.
Then one has
\[
f'_{ij} = k_i f_{ij}, \quad i \in I, j \in J \setminus\{0\}.
\]
Solving equation \eqref{eqn:Z'} for $\Delta_jz'$ one gets
\begin{equation}
\label{eqn:ZDeform}
(\Delta_j z')^2 = (1-k_i^2) f_{ij}^2 + (\Delta_j z)^2.
\end{equation}

This makes sense if and only if the product $(1-k_i^2) f_{ij}^2$ is independent of $i$ and if the right hand side is non-negative.
The former is indeed the case as follows from Figure \ref{fig:GroundViewDeformLengths}.
By Pythagorean theorem one has
\[
f_{i-1,j}^2 - (f'_{i-1,j})^2 = f_{ij}^2 - (f'_{ij})^2\\
\Rightarrow (1 - k_{i-1}^2)f_{i-1,j}^2 = (1 - k_i^2)f_{ij}^2.
\]
The nonnegativity of $(1-k_0^2)f_{0j}^2 + h_j^2$ for all $j$ can be ensured by choosing $k_0$ close enough to $1$.

One has
\[
(1 - k_{i-1}^2)f_{i-1,j}^2 = (1 - k_i^2)f_{ij}^2\ \forall i \Leftrightarrow 1-k_i^2 = \frac{f_{0j}^2}{f_{ij}^2}(1-k_0^2) = \frac{1-k_0^2}{C_i^2}.
\]
If one choses $k_0(t) = \sqrt{1+t}$, then $k_i(t) = \frac{\sqrt{C_i^2 + t}}{C_i}$.
Stretching the line $L_i$ by this factor results in replacing the factor $C_iF_j$ in the formula \eqref{eqn:TSurfCoord} by the factor $\sqrt{C_i^2+t} F_j$ as in formula \eqref{eqn:THedronDeform}.

It remains to compute how do the angles $\eta_i$ and $\theta_i$ change under the above deformation.
This can be done by applying the sine law in the Figure \ref{fig:GroundViewDeformLengths}:
\[
f_{i-1,j} \sin\eta_i = f'_{i-1,j} \sin\eta'_i, \quad
f_{ij} \sin\theta_i = f'_{ij} \sin\theta'_i.
\]
It follows that
\[
\frac{\sin\eta_i(t)}{\sin\eta_i} = \frac{f_{i-1,j}}{f_{i-1,j}(t)} = \frac{1}{k_{i-1}(t)} = \frac{C_{i-1}}{C_{i-1}(t)}, \quad
\frac{\sin\theta_i(t)}{\sin\theta_i} = \frac{f_{ij}}{f_{ij}(t)} = \frac{1}{k_i(t)} = \frac{C_i}{C_i(t)},
\]
as stated in the theorem.

Finally, substituting the obtained value of $k_i(t)$ in \eqref{eqn:ZDeform} one gets
\[
(\Delta_jz(t))^2 = \left(1 - \frac{C_i^2+t}{C_i^2}\right) f_{ij}^2 + (\Delta_jz)^2 = -\frac{t}{C_i^2} C_i^2 f_{0j}^2 + (\Delta_jz)^2 = (\Delta_jz)^2 - t f_{0j}^2.
\]
The square root must be extracted in such a way that $\Delta_jz(t)$ depends continuously on $t$.
This finishes the proof of the theorem.
\end{proof}

\subsection{Deformation of translational T-hedra}
In a translational T-hedron one has
\[
C_i=1 \Rightarrow C_i(t) = \sqrt{1+t}, \quad \psi_i(t) = \eta_i(t), \quad \sin\eta_i(t) = \frac{\sin\eta_i}{\sqrt{1+t}}.
\]
Together with the formulas from Section \ref{sec:TranslTHedra} this can be used to write the isometric deformation in terms of $x_{i0}, x_{0j}, y_i, z_j$.
The result looks more symmetric after reparametrization $\sqrt{1+t} \to e^t$.
Although it is a consequence of the more general Theorem \ref{thm:IsomDefTHedron}, we provide below a simple direct proof.

\begin{theorem}
\label{thm:DeformTranslTHedron}
A translational T-hedron with vertices
\[
\sigma_{ij} =
\begin{pmatrix} x_{i0} + x_{0j}\\ y_i\\ z_j \end{pmatrix}, \quad
x_{00} = y_0 = z_0 = 0
\]
allows an isometric deformation
\begin{gather*}
x_{i0}(t) = e^{-t}x_{i0}, \quad x_{0j}(t) = e^tx_{0j},\\
y_i(t) = \sum_{\alpha=1}^i \sign(\Delta_\alpha y) \sqrt{(\Delta_\alpha y)^2 + (1-e^{-2t})(\Delta_{\alpha 0}x)^2},\\
z_j(t) = \sum_{\beta=1}^j \sign(\Delta_\beta z) \sqrt{(\Delta_\beta z)^2 + (1-e^{2t})(\Delta_{0\beta}x)^2},
\end{gather*}
where
\[
\Delta_\alpha y = y_{\alpha}-y_{\alpha-1}, \quad \Delta_\beta z = z_{\beta}-z_{\beta-1}, \quad
\Delta_{\alpha 0}x = x_{\alpha 0} - x_{\alpha-1,0}, \quad
\Delta_{0\beta}x = x_{0\beta} - x_{0,\beta-1}.
\]
\end{theorem}
\begin{proof}
Consider the parallelogram in the translational T-net $\sigma(t)$ spanned by the vectors
\begin{align*}
\overrightarrow{\sigma_{i-1,j}(t)\sigma_{ij}(t)} &=
\begin{pmatrix}
e^{-t}\Delta_{i0}x\\
\pm\sqrt{(\Delta_i y)^2 + (1-e^{-2t})(\Delta_{i0}x)^2}\\
0
\end{pmatrix},\\
\overrightarrow{\sigma_{i,j-1}(t)\sigma_{ij}(t)} &=
\begin{pmatrix}
e^t\Delta_{0j}x\\
0\\
\pm\sqrt{(\Delta_j z)^2 + (1-e^{2t})(\Delta_{0j}x)^2}
\end{pmatrix}.
\end{align*}
The norms of these vectors and their inner product do not depend on $t$.
Therefore all faces of $\sigma$ move as rigid bodies.
\end{proof}
While the proof of the above Theorem is straightforward, the formulas are not so obvious.
They can be obtained by transforming back and forth between the $x,y,z$- and the $\eta, f, g, h$-data, see end of Section \ref{sec:TranslTHedra}.

Formulas in Theorem \ref{thm:DeformTranslTHedron} make apparent the boundaries of deformability of a translational T-hedron.
For $t<0$ the deformation stops as soon as one of the radicands $(\Delta_i y)^2 + (1-e^{-2t})(\Delta_{i0}x)^2$ vanishes.
This corresponds to $\Delta_i y(t) = 0$, that is the $(i-1)$-st and the $i$-th profile planes become coincident, and the $i$-th profile strip parallel to the profile planes.
Similarly, for $t>0$ the deformation stops as soon as two consecutive trajectory planes become coincident, that is a trajectory strip becomes parallel to trajectory planes.

\begin{figure}[ht]
\begin{center}
\includegraphics[width=.8\textwidth]{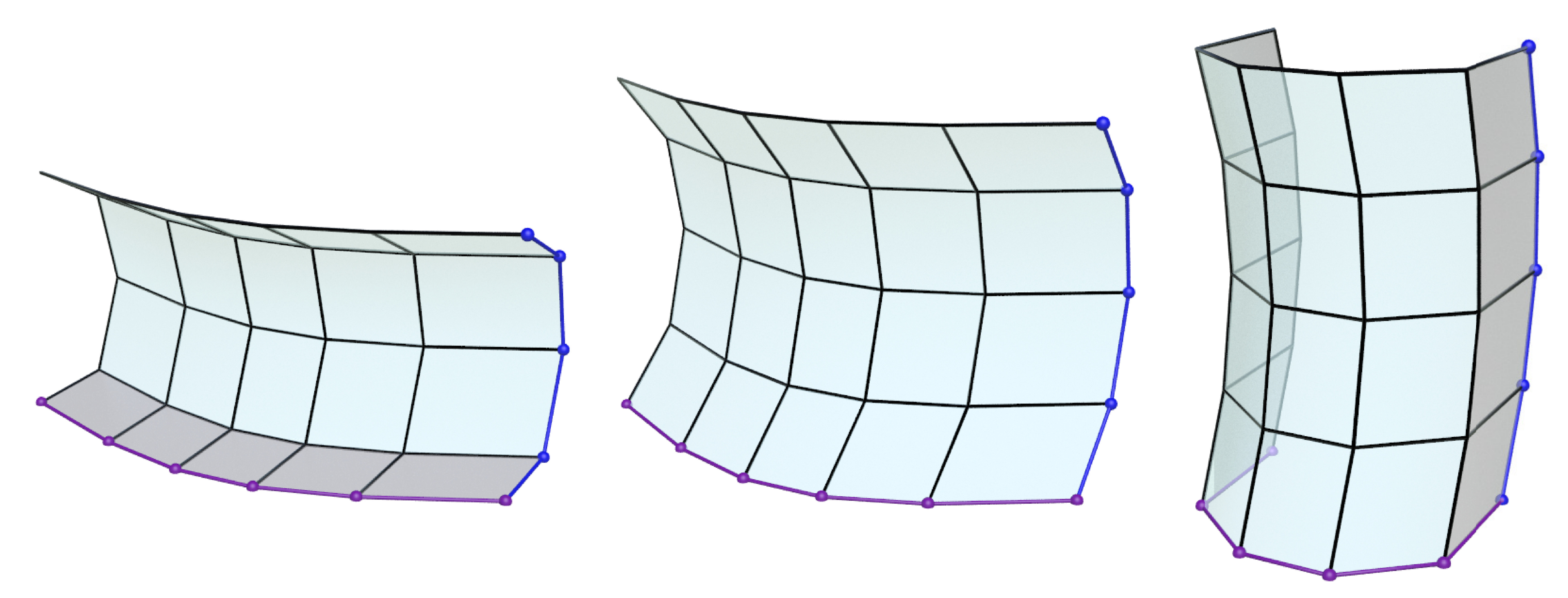}
\end{center}
\caption{Isometric deformation of a translational T-hedron. Shaded faces on the left hand side are parallel to the trajectory planes, shaded faces on the right hand side are parallel to the profile planes; they prevent further deformation.}
\label{fig:BoundFlex}
\end{figure}

\begin{figure}[ht]
\begin{center}
\includegraphics[width=\textwidth]{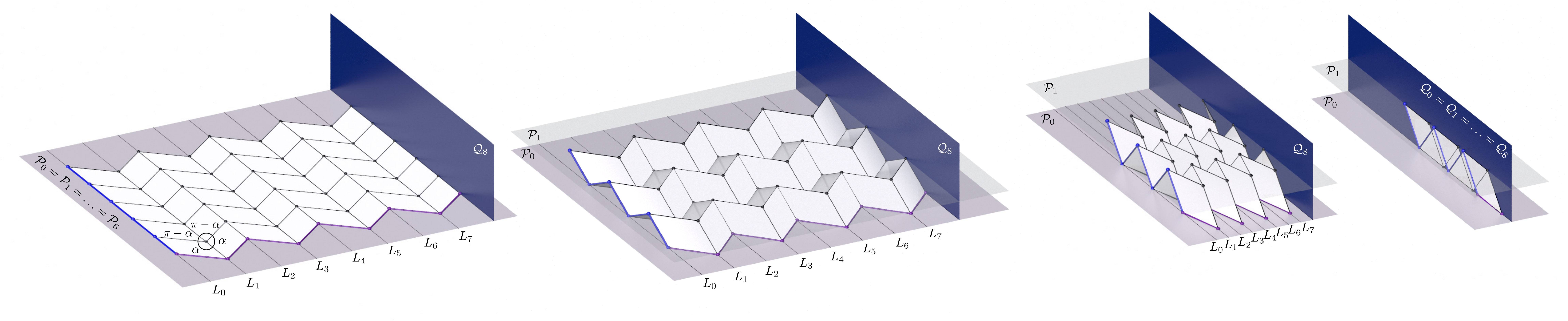}
\end{center}
\caption{Miura-ori as a translational T-hedron.}
\label{fig:MiuraOri}
\end{figure}

\begin{example}
Miura-ori is a translational T-hedron that can be constructed as follows.
Subdivide a strip between two parallel lines $L_0$ and $L_1$ into congruent parallelograms.
Reflect this strip in the line $L_1$ and denote the image of $L_0$ under this reflection by $L_2$, then reflect the strip enclosed by $L_1$ and $L_2$ in the line $L_2$, and so on.
In the end, lift every second trajectory polygon to the same level, see Figure \ref{fig:MiuraGround}.

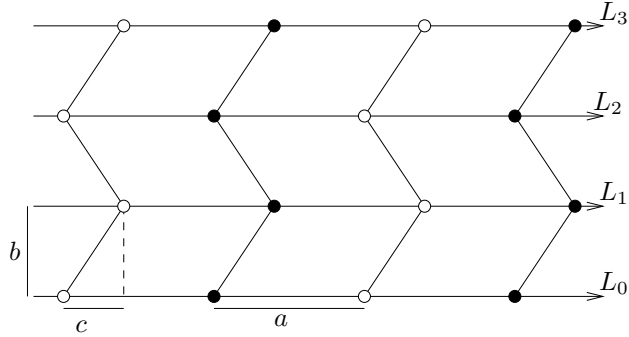
\begin{figure}[ht]
\begin{center}
\begin{picture}(0,0)%
\includegraphics{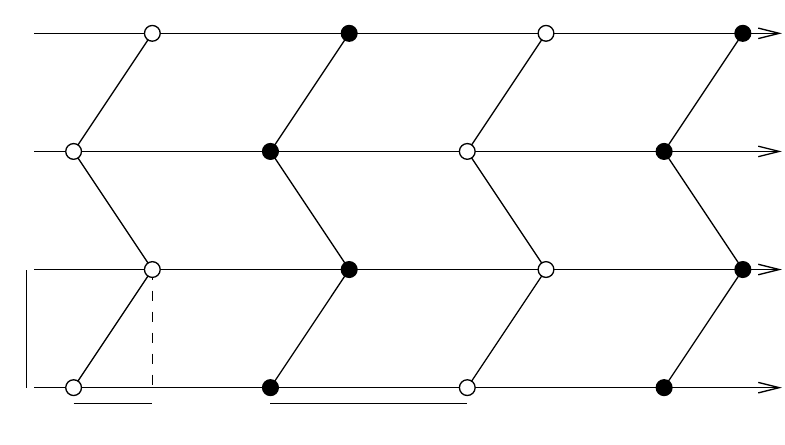}%
\end{picture}%
\setlength{\unitlength}{3315sp}%
\begingroup\makeatletter\ifx\SetFigFont\undefined%
\gdef\SetFigFont#1#2#3#4#5{%
  \reset@font\fontsize{#1}{#2pt}%
  \fontfamily{#3}\fontseries{#4}\fontshape{#5}%
  \selectfont}%
\fi\endgroup%
\begin{picture}(4482,2530)(-419,-1511)
\put(4006,-1141){\makebox(0,0)[lb]{\smash{{\SetFigFont{10}{12.0}{\rmdefault}{\mddefault}{\updefault}{\color[rgb]{0,0,0}$L_0$}%
}}}}
\put(4006,-466){\makebox(0,0)[lb]{\smash{{\SetFigFont{10}{12.0}{\rmdefault}{\mddefault}{\updefault}{\color[rgb]{0,0,0}$L_1$}%
}}}}
\put(3961,209){\makebox(0,0)[lb]{\smash{{\SetFigFont{10}{12.0}{\rmdefault}{\mddefault}{\updefault}{\color[rgb]{0,0,0}$L_2$}%
}}}}
\put( 91,-1456){\makebox(0,0)[lb]{\smash{{\SetFigFont{10}{12.0}{\rmdefault}{\mddefault}{\updefault}{\color[rgb]{0,0,0}$c$}%
}}}}
\put(1576,-1411){\makebox(0,0)[lb]{\smash{{\SetFigFont{10}{12.0}{\rmdefault}{\mddefault}{\updefault}{\color[rgb]{0,0,0}$a$}%
}}}}
\put(-404,-916){\makebox(0,0)[lb]{\smash{{\SetFigFont{10}{12.0}{\rmdefault}{\mddefault}{\updefault}{\color[rgb]{0,0,0}$b$}%
}}}}
\put(4006,884){\makebox(0,0)[lb]{\smash{{\SetFigFont{10}{12.0}{\rmdefault}{\mddefault}{\updefault}{\color[rgb]{0,0,0}$L_3$}%
}}}}
\end{picture}%
\end{center}
\caption{Miura polyhedron. Black vertices are lifted to the height $d$.}
\label{fig:MiuraGround}
\end{figure}

The vertex coordinates of the obtained polyhedral surface are
\begin{equation}
\label{eqn:MiuraCoord}
\sigma_{ij} =
\begin{pmatrix}
x_{i0} + x_{0j}\\ y_i\\ z_j
\end{pmatrix}, \
x_{0j} = ja, \ y_i = ib, \quad
x_{i0} =
\begin{cases}
0, &\text{ if }i \equiv 0 \pmod 2,\\
c, &\text{ if }i \equiv 1 \pmod 2
\end{cases}, 
z_j =
\begin{cases}
0, &\text{ if }j \equiv 0 \pmod 2,\\
d, &\text{ if }j \equiv 1 \pmod 2.
\end{cases}
\end{equation}
A Miura polyhedron is very symmetric.
If one extends the construction to $\ZZ \times \ZZ$, then the polyhedron is invariant under translations $i \mapsto i+2$, $j \mapsto j+2$, reflections in the lines $L_i$, and $180^\circ$-rotations about the vertical axes through the midpoints of trajectory edges.
It follows that all faces of this polyhedron are congruent to each other.
In particular the angles about every vertex are $\alpha, \pi - \alpha, \pi-\alpha, \alpha$.
Since
\[
\alpha + (\pi - \alpha) + (\pi-\alpha) + \alpha = 2\pi, \quad
\alpha - (\pi - \alpha) + (\pi-\alpha) - \alpha = 0,
\]
one is led to believe that during the isometric deformation a Miura polyhedron attains two flat positions: one without overlap and one with many overlaps.
This can be confirmed by analyzing the formulas of the deformation.
Substituting the data from \eqref{eqn:MiuraCoord} into the formulas of Theorem \ref{thm:DeformTranslTHedron} one obtains
\[
x_{0j}(t) = ja(t), \ y_i(t) = ib(t), \quad
x_{i0}(t) =
\begin{cases}
0, &\text{ if }i \equiv 0 \pmod 2,\\
c(t), &\text{ if }i \equiv 1 \pmod 2,
\end{cases} \
z_j(t) =
\begin{cases}
0, &\text{ if }j \equiv 0 \pmod 2,\\
d(t), &\text{ if }j \equiv 1 \pmod 2,
\end{cases}
\]
where
\[
a(t) = e^ta, \quad b(t) = \sqrt{b^2 + (1-e^{-2t})c^2}, \quad c(t) = e^{-t}c, \quad d(t) = \sqrt{d^2 + (1-e^{2t})a^2}.
\]
The deformation is illustrated in Figure \ref{fig:MiuraOri}.
At $t = t_+ := \log\frac{\sqrt{a^2+d^2}}{a} > 0$ one has $d(t) = 0$, the polyhedron flattens to a parallelogram net like in Figure \ref{fig:MiuraGround} with
\[
a(t_+) = \sqrt{a^2+d^2}, \quad b(t_+) = \sqrt{b^2 + \frac{c^2d^2}{a^2 + d^2}}, \quad c(t_+) = \frac{ac}{\sqrt{a^2+d^2}}.
\]
At $t = t_- := \log\frac{c}{\sqrt{b^2+c^2}}$ one has $b(t)= 0$, the polyhedron flattens in the $xz$-plane to a strip with multiple overlaps:
\[
a(t_-) = \frac{ac}{\sqrt{b^2+c^2}}, \quad c(t_-) = \sqrt{b^2+c^2}, \quad d(t_-) = \sqrt{d^2 + \frac{a^2b^2}{b^2+c^2}}.
\]

Setting $d=0$ in the above formulas yields a parametrization of an isometric deformation of a flat Miura-ori pattern with metric data as shown in Figure \ref{fig:MiuraGround}.
\end{example}

\begin{figure}[ht]
\begin{center}
\includegraphics[width=.8\textwidth]{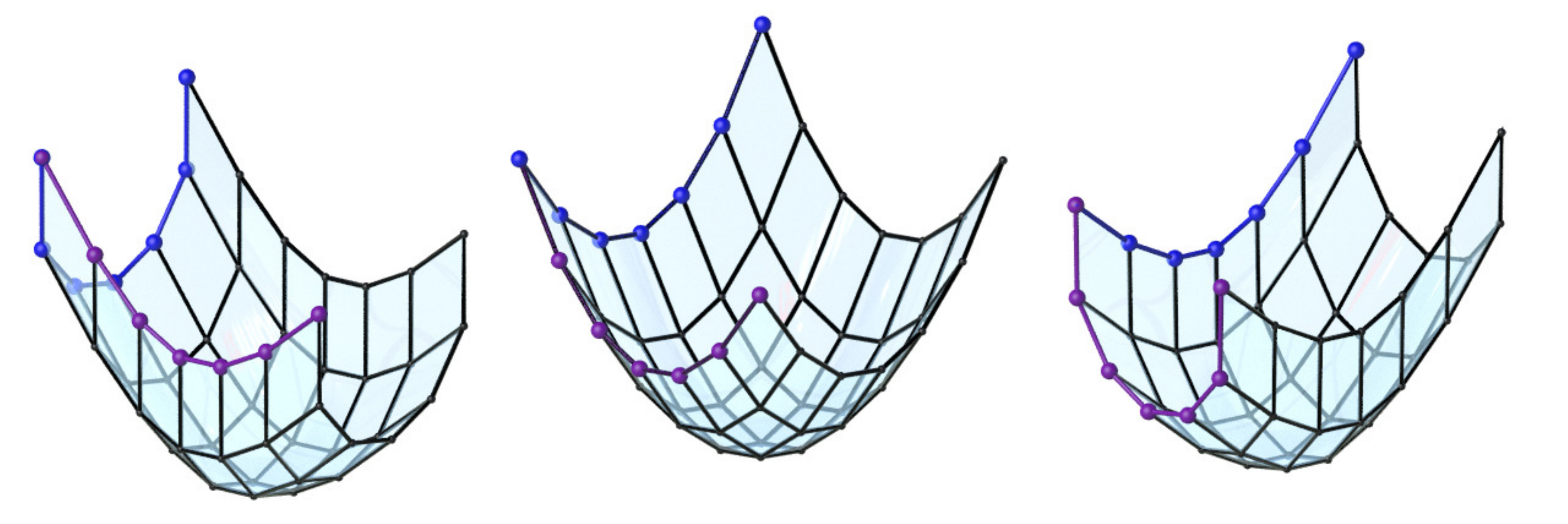}
\end{center}
\caption{Isometric deformation of a translational T-hedron inscribed in a paraboloid of revolution.}
\label{fig:DeformParabDiscrTransl}
\end{figure}

\begin{example}
\label{exl:ParTranslDiscr}
Take a square grid in the $xy$-plane and lift its vertices to the paraboloid of revolution.
One can show that the squares lift to planar quadrilaterals.
In fact, the result is a translational T-hedron.
An isometric deformation of this T-hedron is shown in Figure \ref{fig:DeformParabDiscrTransl}.
Compare this with Example \ref{exl:ParTranslSm}.
\end{example}

\subsection{Deformation of discrete molding surfaces}
For molding surfaces the formulas of Theorem \eqref{thm:IsomDefTHedron} simplify in the following way.
\begin{theorem}
A discrete molding surface with vertices
\[
\sigma_{ij} = \begin{pmatrix} \tau_{ij}\\ z_j \end{pmatrix}, \quad
\tau_{ij} = \sum_{\alpha=1}^i g_{\alpha 0}
\begin{pmatrix} -\sin\frac{\phi_{\alpha-1} + \phi_\alpha}2\\ \cos\frac{\phi_{\alpha-1} + \phi_\alpha}2 \end{pmatrix} + F_j
\begin{pmatrix} \cos\phi_i\\ \sin\phi_i \end{pmatrix}
\]
allows an isometric deformation
\begin{gather*}
g_{i0}(t) = g_{i0}, \quad F_j(t) = \sqrt{1+t}\ F_j, \quad z_j(t) = \sum_{\beta=1}^j \sign(\Delta_\beta z) \sqrt{(\Delta_\beta z)^2 - t(\Delta_\beta F)^2},\\
\phi_i(t) = 2 \sum_{\alpha=1}^i \eta_\alpha(t), \quad \sin\eta_\alpha(t) = \frac{\sin\eta_\alpha}{\sqrt{1+t}}.
\end{gather*}
\end{theorem}

\subsection{Deformation of axial T-hedra}
For an axial T-hedron not much can be simplified in the formulas of Theorem \ref{thm:IsomDefTHedron}.
More simplifications occur for discrete surfaces of revolution.

\begin{theorem}
An axial T-hedron \eqref{eqn:AnDescrAxisThedra} admits an isometric deformation
\[
\sigma_{ij}(t) = 
\begin{pmatrix}
C_i(t) F_j \cos\phi_i(t)\\ C_i(t) F_j \sin\phi_i(t)\\ z_j(t)
\end{pmatrix},
\]
where
\begin{gather*}
\phi_i(t) = \sum_{\alpha=1}^i \eta_\alpha(t) + \sum_{\alpha=1}^i \theta_\alpha(t),\\
\sin\eta_i(t) = \frac{C_{i-1}}{C_{i-1}(t)}\sin\eta_i, \quad \sin\theta_i(t) = \frac{C_i}{C_i(t)}\sin\theta_i, \quad
C_i(t) = \sqrt{C_i^2 + t},\\
z_{j}(t) = \sum_{\beta=1}^j \sign(\Delta_\beta z) \sqrt{(\Delta_\beta z)^2 - t (\Delta_\beta F)^2}.
\end{gather*}

See the top row of Figure \ref{fig:DefAxialTHedron} for an example of an isometric deformation of an axial T-hedron.

\begin{figure}[ht]
\begin{center}
\includegraphics[width=\textwidth]{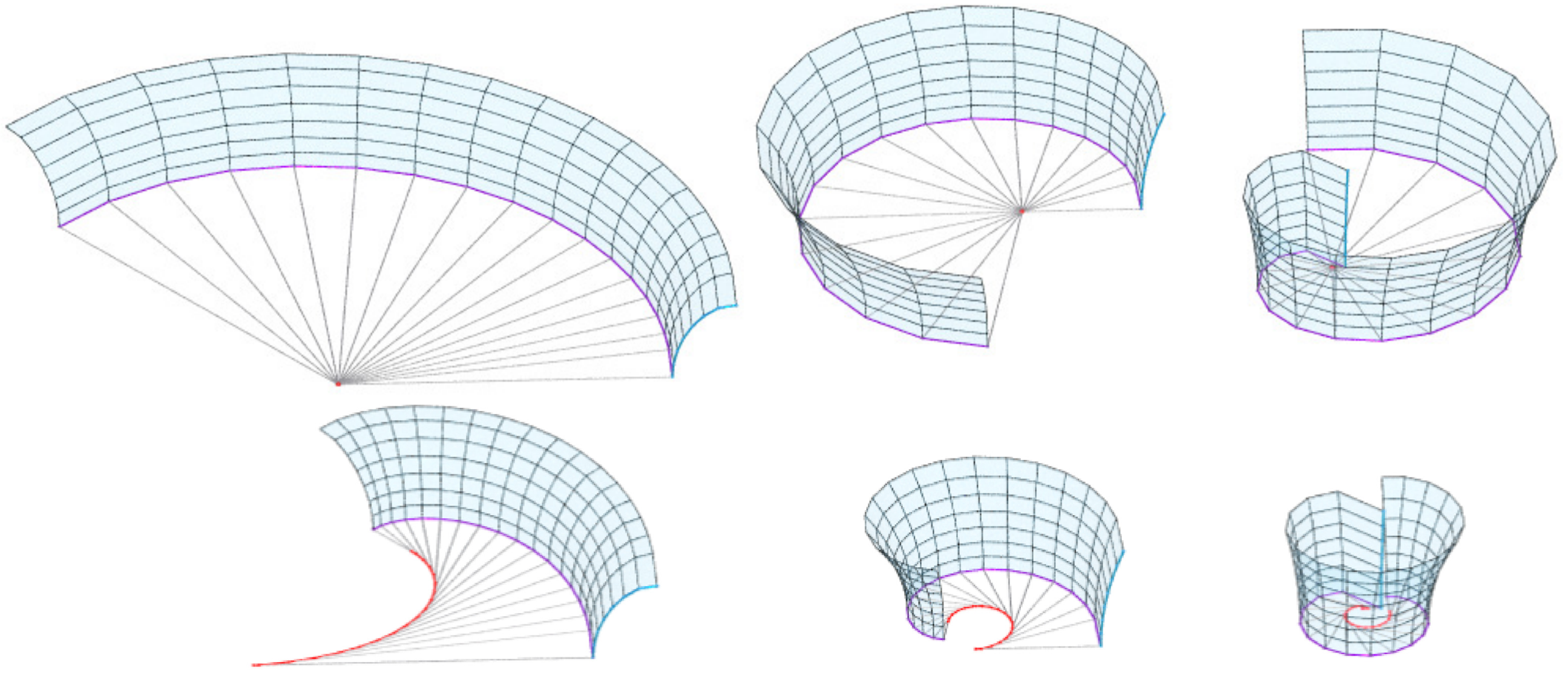}
\end{center}
\caption{Top row: an isometric deformation of an axial T-hedron. Bottom row: a general T-hedron parallel to the axial T-hedron in the top row. The parallelity is preserved during the deformation.}
\label{fig:DefAxialTHedron}
\end{figure}

A discrete surface of revolution \eqref{eqn:DiscrSurfRev} admits an isometric deformation
\[
\sigma_{ij}(t) = 
\begin{pmatrix}
F_j(t) \cos\phi_i(t)\\ F_j(t) \sin\phi_i(t)\\ z_j(t)
\end{pmatrix},
\]
where
\begin{gather*}
F_j(t) = \sqrt{1+t}\ F_j, \quad z_j(t) = \sum_{\beta=1}^j \sign(\Delta_\beta z) \sqrt{(\Delta_\beta z)^2 - t(\Delta_\beta F)^2},\\
\phi_i(t) = 2 \sum_{\alpha=1}^i \eta_\alpha(t), \quad \sin\eta_\alpha(t) = \frac{\sin\eta_\alpha}{\sqrt{1+t}}.
\end{gather*}
\end{theorem}

\begin{figure}[ht]
\begin{center}
\includegraphics[width=.8\textwidth]{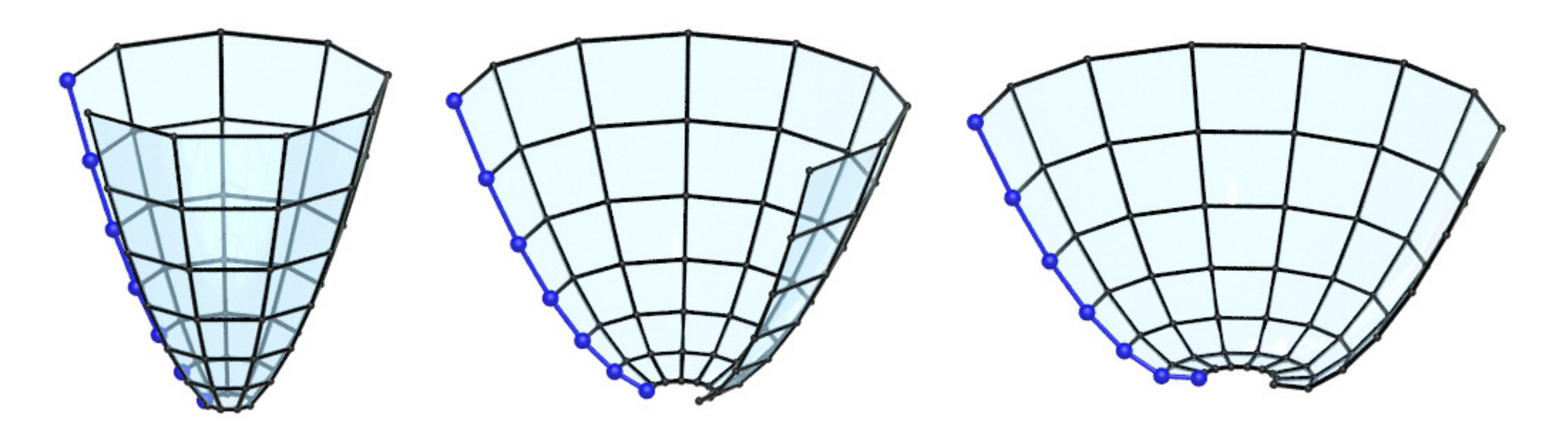}
\end{center}
\caption{Isometric deformation of a discrete paraboloid of revolution.}
\label{fig:DeformParabDiscrRev}
\end{figure}

Figure \ref{fig:DeformParabDiscrRev} shows the deformation of a discrete surface of revolution inscribed into the paraboloid of revolution.
Compare this to Figure \ref{fig:DeformParaboloidRev}.

\subsection{Parallel pairs of T-hedra}
\begin{definition}
Two T-hedra $\sigma, \sigma' \colon I \times J \to \R^3$ are called \emph{parallel} if all of their corresponding pairs of faces are parallel.
\end{definition}
Observe that all corresponding pairs of faces of $\sigma$ and $\sigma'$ are parallel if and only if all corresponding pairs of edges are parallel.

\begin{lemma}
\label{lem:ParallelCrit}
T-hedra $\sigma$ and $\sigma'$ are parallel if and only if their initial profile polygons are parallel and the lines $L_i, L'_i$ and $M_i, M'_i$ are pairwise parallel:
\[
\sigma_{0\bullet} \parallel \sigma'_{0\bullet}, \text{ that is, } \sigma_{0,j-1}\sigma_{0j} \parallel \sigma'_{0,j-1}\sigma'_{0j} \text{ for all }j, \text{ and } L_i \parallel L'_i, M_i \parallel M'_i \text{ for all }i.
\]
\end{lemma}
\begin{proof}
If $\sigma \parallel \sigma'$, then their corresponding edges are parallel, in particular one has $\sigma_{0\bullet} \parallel \sigma'_{0\bullet}$.
Also $\sigma \parallel \sigma'$ implies $P_0 \parallel P'_0$ and $Q_i \parallel Q'_i$, which implies that $L_i = P_0 \cap Q_i$ is parallel to $L'_i = P'_0 \cap Q'_i$.
Finally, the line $M_i$ is a line in $P_0$ orthogonal to all edges $\sigma_{i-1,j}\sigma_{ij}$.
Since $\sigma_{i-1,j}\sigma_{ij} \parallel \sigma'_{i-1,j}\sigma'_{ij}$, it follows that $M_i \parallel M'_i$.

In the opposite direction, the parallelity $L_i \parallel L'_i, M_i \parallel M'_i$ implies that the ground views of $\sigma$ and $\sigma'$ are parallel T-nets: $\tau \parallel \tau'$.
In view of that, the parallelity $\sigma_{0\bullet} \parallel \sigma'_{0\bullet}$ implies that the faces of the first profile strip of $\sigma$ are parallel to the faces of the first profile strip of $\sigma'$.
This, together with $\tau \parallel \tau'$, implies $\sigma_{1\bullet} \parallel \sigma'_{1\bullet}$.
Continuing in the same manner one arrives at $\sigma \parallel \sigma'$.
\end{proof}

\begin{theorem}
\label{thm:DeformParallel}
If two T-hedra $\sigma$ and $\sigma'$ are parallel, then their isometric deformations $\sigma(t)$ and $\sigma'(t)$ described in Theorem \ref{thm:IsomDefTHedron} are parallel for all $t$.
\end{theorem}
\begin{proof}
By construction, the lines $L_0(t)$ and $L'_0(t)$ remain coincident with the $x$-axis.
Formulas of Theorem \ref{thm:IsomDefTHedron} imply that $\phi_i(t) = \phi'_i(t)$ and $\psi_i(t) = \psi'_i(t)$, so that $L_i(t) \parallel L'_i(t)$ and $M_i(t) \parallel M'_i(t)$.
Due to Lemma \ref{lem:ParallelCrit} it remains to show that the profile polygons $\sigma_{0\bullet}(t)$ and $\sigma'_{0\bullet}(t)$ remain parallel during the deformation.
Let
\[
k := \frac{\Delta_j z}{\Delta_j F} = \frac{\Delta_j z'}{\Delta_j F'}
\]
be the common slope of the edges $\sigma_{0,j-1}\sigma_{0j}$ and $\sigma'_{0,j-1}\sigma'_{0j}$ over the line $L_0$.
Formulas of Theorem \ref{thm:IsomDefTHedron} imply that at time $t$ the slope of the $i$-th edge of $\sigma_{0\bullet}(t)$ equals
\[
\frac{\sign(\Delta_jz)\sqrt{(\Delta_jz)^2 - t(\Delta_jF)^2}}{\sqrt{1+t}\ \Delta_jF} =
\sign(k) \frac{\sqrt{k^2-t}}{\sqrt{1+t}}.
\]
One obtains the same value when computing the slope of the $i$-th edge of $\sigma'_{0\bullet}(t)$, and the theorem is proved.
\end{proof}

\begin{remark}
One can show that the flexibility of a simply-connected quad-surface depends only on the values of the dihedral angles between its faces.
This was shown in \cite{stachel2010kinematic} for $3\times 3$-quad-surfaces, the general case is similar.
Parallel quad-surfaces have pairwise equal dihedral angles, therefore Theorem \ref{thm:DeformParallel} is a special case of this more general statement.
\end{remark}

\begin{theorem}
Let $\sigma$ be a T-hedron such that no two consecutive profile planes $Q_{i-1}$ and $Q_i$ are parallel.
Then there is an axial T-hedron $\sigma'$ which is parallel to $\sigma$.
\end{theorem}
\begin{proof}
Geometrically $\sigma'$ can be obtained from $\sigma$ by translating the profile planes $Q_2, \ldots, Q_m$ so that they pass through the intersection line of $Q_0$ and $Q_1$, keeping the profile polygon $\sigma_{0\bullet}$ and constructing the rest of $\sigma'$ from the parallelity condition.
Analytically one can use for $\sigma'$ the data $\phi, \psi, f, g, z$ of $\sigma$ (see Theorem \ref{thm:THedronCoord}), keeping $\phi, \psi, f$, and $z$, and adjusting $g$ so that the relation \eqref{eqn:AxialCrit} holds.
Then $\sigma'$ is axial, has the same profile $\sigma_{0\bullet}$ and the same directions of lines $L_i$ and $M_i$.
By Lemma \ref{lem:ParallelCrit}, $\sigma'$ and $\sigma$ are parallel.
Observe that the faces of an axial T-hedron are never self-intersecting.
\end{proof}

Figure \ref{fig:DefAxialTHedron} shows isometric deformations of a T-hedron (bottom row) and of a parallel to it axial T-hedron (top row).

\section{T-surfaces}
\subsection{Definition and basic properties}
Let $\Omega \subset \R^2$ be an open set.
A smooth map $\sigma \colon \Omega \to \R^3$ is called a \emph{regular parametrized surface} if for every $(u,v) \in \Omega$ the vectors
\[
\sigma_u(u,v) = \frac{\partial\sigma}{\partial u}(u,v) \text{ and } \sigma_v(u,v) = \frac{\partial\sigma}{\partial v}(u,v)
\]
are linearly independent.
We do not require $\sigma$ to be injective.
However, the regularity condition implies that $\sigma$ is locally injective.

\begin{definition}
\label{dfn:TSurface}
Let $U = [0,a], V = [0,b] \subset \R$ be two intervals in the real line.
A regular parametrized surface $\sigma \colon U \times V \to \R^3$ is called a \emph{T-surface} if its family of coordinate curves $\sigma|_{u = u_0}$ and $\sigma|_{v = v_0}$ for all $u_0 \in U$ and $v_0 \in V$ satisfy the following conditions:
\begin{enumerate}
\item
They form a conjugate system, that is $\sigma_{uv} \in \Span\{\sigma_u, \sigma_v\}$ everywhere in $U \times V$.
\item
Each of them is planar but contains no straight line segments. In particular, each coordinate curve spans a unique plane.
\item
Denote the plane spanned by $\sigma|_{u=u_0}$ by $Q_{u_0}$ and the plane spanned by $\sigma|_{v=v_0}$ by $P_{v_0}$.
Then for all $u \in U$ and $v \in V$ the planes $Q_u$ and $P_v$ are orthogonal.
\item
There is no interval $U' \subset U$ for which all planes $Q_u, u \in U'$ coincide.
The same holds for the planes $P_v$.
\end{enumerate}
\end{definition}

Similarly to Lemma \ref{lem:OrthPlanesDiscr} one proves that at least one of the families
\[
\{P_v \mid v \in V\}, \quad \{Q_u \mid u \in U\}
\]
consist of parallel planes.
Without loss of generality let these be the planes $\{P_v\}$.
In this case we call $\{P_v\}$ the \emph{trajectory planes}, and $\{Q_u\}$ the \emph{profile planes} of the surface $\sigma$.
Accordingly, the coordinate curves $\sigma|_{v=v_0}$ are called the \emph{trajectory curves}, and $\sigma|_{u=u_0}$ are called the \emph{profile curves} of $\sigma$.

\begin{example}
Surfaces of revolution
\[
\sigma(u,v) =
\begin{pmatrix}
f(v) \cos u\\ f(v) \sin u\\ g(v)
\end{pmatrix}
\]
are T-surfaces.
The trajectory planes are orthogonal to the axis of rotation, the profile planes pass through the axis.
The trajectory curves are circular arcs, the profile curves are pairwise congruent.
\end{example}

In the next section we will prove that every T-surface is generated from a profile curve by revolving this curve around a variable axis and at the same time scaling it in the direction orthogonal to the axis.
If the axis lies at infinity, one obtains special surfaces of translation.

\subsection{Analytic description of T-surfaces}
\begin{theorem}
\label{thm:DescrTSurf}
In a coordinate system in which the trajectory planes are parallel to the $xy$-plane, the parametrization of a T-surface has the following form:
\begin{equation}
\label{eqn:TSurfParam}
\sigma(u,v) =
\begin{pmatrix}
\gamma(u) + f(v) \xi(u)\\
z(v)
\end{pmatrix}.
\end{equation}
Here $\gamma \colon U \to \R^2$ is a regular curve (that is, $\dot\gamma(u) \ne 0$ for all $u \in U$), $f \colon V \to \R$ a smooth function, and $\xi \colon U \to \R^2$ a vector field such that the following holds:
\begin{enumerate}
\item
$f(0) = 0$ (recall that $V = [0,b]$).
\item
Functions $f$ and $z$ are locally affinely independent, that is there is no interval $V' \subset V$ and no reals $\lambda_1, \lambda_2, \lambda_3$ such that $\lambda_1f(v) + \lambda_2z(v) + \lambda_3 = 0$ for all $v \in V'$.
\item
The curve $v \mapsto (f(v), z(v))$ is regular, that is there is no $v \in V$ such that $\dot f(v) = \dot z(v) = 0$.
\item
The curve $\gamma$ contains no straight line segments.
\item
Vectors $\xi(u)$ and $\dot\gamma(u)$ are linearly independent (in particular, $\xi(u) \ne 0$) for all $u \in U$.
\item
Vectors $\dot\xi(u)$ and $\dot\gamma(u)$ are linearly dependent for all $u \in U$, that is $\dot\xi(u) = \lambda(u) \dot\gamma(u)$ for some function $\lambda \colon U \to \R$.
Besides, $f(v)\lambda(u) \ne -1$ for all $u \in U$ and $v \in V$.
\end{enumerate}
Conversely, for every choice of $\gamma$, $\xi$, $f$, and $z$ subject to the above conditions the surface \eqref{eqn:TSurfParam} is a T-surface.
\end{theorem}

\begin{figure}[ht]
\begin{center}
\begin{picture}(0,0)%
\includegraphics{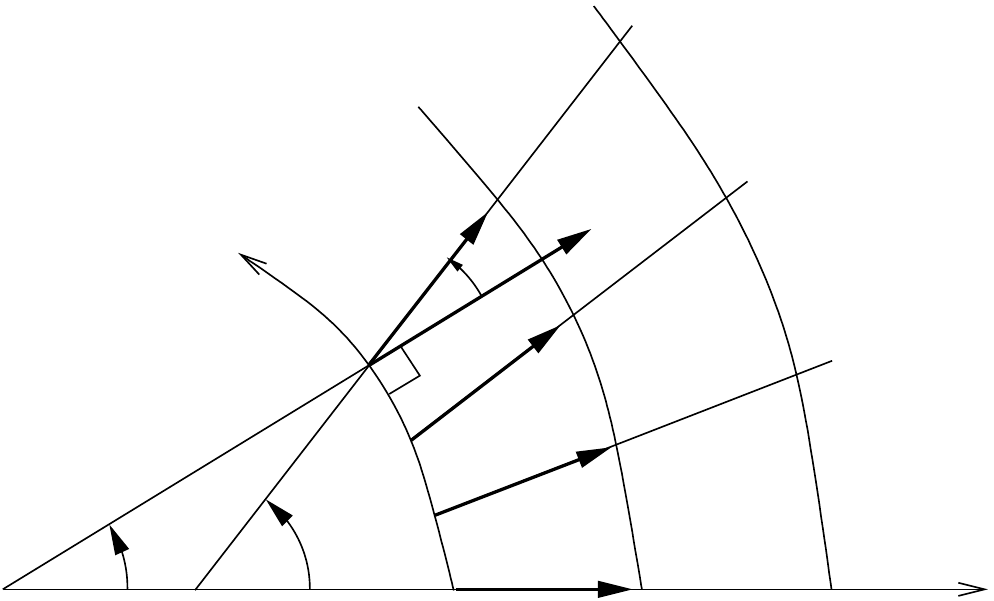}%
\end{picture}%
\setlength{\unitlength}{4144sp}%
\begingroup\makeatletter\ifx\SetFigFont\undefined%
\gdef\SetFigFont#1#2#3#4#5{%
	\reset@font\fontsize{#1}{#2pt}%
	\fontfamily{#3}\fontseries{#4}\fontshape{#5}%
	\selectfont}%
\fi\endgroup%
\begin{picture}(4530,2731)(-1316,-1463)
\put(111,-1176){\makebox(0,0)[lb]{\smash{{\SetFigFont{10}{12.0}{\rmdefault}{\mddefault}{\updefault}{\color[rgb]{0,0,0}$\phi(u)$}%
}}}}
\put(-783,174){\makebox(0,0)[lb]{\smash{{\SetFigFont{10}{12.0}{\rmdefault}{\mddefault}{\updefault}{\color[rgb]{0,0,0}$\gamma(u)=\tau(u,0)$}%
}}}}
\put(560,319){\makebox(0,0)[lb]{\smash{{\SetFigFont{10}{12.0}{\rmdefault}{\mddefault}{\updefault}{\color[rgb]{0,0,0}$\xi(u)$}%
}}}}
\put(-692,-1267){\makebox(0,0)[lb]{\smash{{\SetFigFont{10}{12.0}{\rmdefault}{\mddefault}{\updefault}{\color[rgb]{0,0,0}$\psi(u)$}%
}}}}
\put(854, 50){\makebox(0,0)[lb]{\smash{{\SetFigFont{10}{12.0}{\rmdefault}{\mddefault}{\updefault}{\color[rgb]{0,0,0}$\eta(u)$}%
}}}}
\end{picture}%
\end{center}
\caption{The ground view of a T-surface: notations for Theorem \ref{thm:DescrTSurf} and Lemma \ref{lem:XiData}.}
\label{fig:TSurfData}
\end{figure}
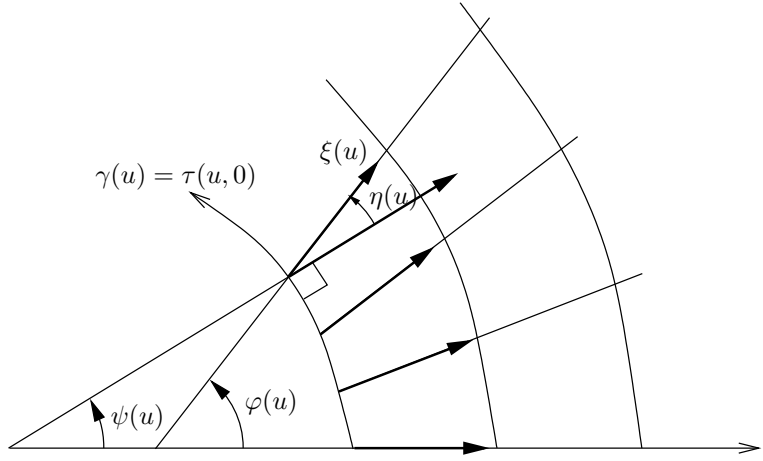

Figure \ref{fig:TSurfData} shows the entries of the formula \eqref{eqn:TSurfParam} in the ground view (that is, in the orthogonal projection to the plane $P_0$).

Observe that scaling $f$ up and $\xi$ down by the same constant does not change $\sigma$.

\begin{proof}
Let us first check that every $\sigma$ given by \eqref{eqn:TSurfParam} with $\gamma, f, \xi, z$ satisfying the conditions listed in the theorem is a T-surface.
Compute the partial derivatives:
\[
\sigma_u(u,v) =
\begin{pmatrix}
\dot\gamma(u) + f(v)\dot\xi(u)\\ 0
\end{pmatrix}, \quad
\sigma_v(u,v) =
\begin{pmatrix}
\dot f(v) \xi(u)\\ \dot z(v)
\end{pmatrix}, \quad
\sigma_{uv}(u,v) =
\begin{pmatrix}
\dot f(v) \dot\xi(u)\\ 0
\end{pmatrix}.
\]
Due to the sixth condition one has $\sigma_u(u,v) = ((1+f(v)\lambda(u))\dot\xi(u), 0) \ne 0$.
Third condition together with the fifth ensure that $\sigma_v \ne 0$ and that the vectors $\sigma_u$ and $\sigma_v$ are linearly independent.
Thus $\sigma$ is a regular parametrized surface.
Because of the sixth condition the vectors $\sigma_{uv}$ and $\sigma_u$ are linearly dependent, thus the coordinate curves form a conjugate system.
Consider the coordinate curves $\sigma|_{u=u_0}$ and $\sigma|_{v=v_0}$:
\[
\sigma(u_0,v) =
\begin{pmatrix}
\gamma(u_0) + f(v)\xi(u_0)\\ z(v)
\end{pmatrix}, \quad
\sigma(u,v_0) =
\begin{pmatrix}
\gamma(u) + f(v_0) \xi(u)\\ z(v_0)
\end{pmatrix}.
\]
The former is contained in a plane parallel to the vectors $(\xi(u_0),0)$ and $(0,1)$ and contains no straight line segments due to the second condition.
The latter is contained in the plane $z=z(v_0)$ and contains no straight line segments because its tangent $(\dot\gamma(u) + f(v_0)\dot\xi(u), 0)$ is parallel to the tangent of the curve $\gamma$ at the corresponding point, and $\gamma$ contains no straight line segments by the fourth condition.
By the above, the profile and the trajectory planes have the following form:
\[
Q_{u_0} = \left\{
\begin{pmatrix}
\gamma(u_0)\\ 0
\end{pmatrix} + s
\begin{pmatrix}
\xi(u_0)\\ 0
\end{pmatrix} + t
\begin{pmatrix}
0\\ 1
\end{pmatrix}
\right\}, \quad
P_{v_0} = \{(x,y,z) \in \R^3 \mid z = z(v_0)\}.
\]
It follows that $Q_u$ is orthogonal to $P_v$ for every $(u,v) \in U \times V$.
Planes $P_v$ coincide on some interval $V' \subset V$ if and only if the restriction of the function $z$ to this interval is constant.
But this is forbidden by the second condition in the theorem.
Planes $Q_u$ coincide on some interval $U' \subset U$ if and only if $\xi$ has a constant direction on this interval, and the curve $\gamma$ is parallel to that direction.
But then $\gamma$ contains a straight line segment which is vorbidden by the fourth condition.

Now let us prove that every T-surface has a parametrization of the form \eqref{eqn:TSurfParam} with $\gamma, \xi, f, z$ satisfying the conditions listed in the theorem.
By assumption, the coordinate system is chosen so that the trajectory planes $P_{v_0} \supset \sigma |_{v=v_0}$ are parallel to the $xy$-plane.
This implies that
\[
\sigma(u,v) =
\begin{pmatrix} \tau(u,v)\\ z(v) \end{pmatrix}
\]
for some smooth maps $\tau \colon U \times V \to \R^2$ and $z \colon V \to \R$, where $z$ by the fourth condition from Definition \ref{dfn:TSurface} is nowhere locally constant.
One has
\[
\sigma_u = \begin{pmatrix} \tau_u\\ 0 \end{pmatrix}, \quad
\sigma_v = \begin{pmatrix} \tau_v\\ \dot z \end{pmatrix}, \quad
\sigma_{uv} = \begin{pmatrix} \tau_{uv}\\ 0 \end{pmatrix}.
\]
By definition of T-surfaces one has
\[
\sigma_{uv} = \lambda \sigma_u + \mu \sigma_v
\]
for some functions $\lambda, \mu \colon U \times V \to \R$.
Since the $z$-coordinates of $\sigma_u$ and of $\sigma_{uv}$ vanish everywhere, and the $z$-coordinate of $\sigma_{v}$ is nonzero on a dense set, it follows that $\mu$ is identically zero, thus
\begin{equation}
\label{eqn:SigmaUParallel}
\sigma_{uv}(u,v) = \lambda(u,v)\sigma_u(u,v).
\end{equation}
Let us now turn to the profile planes of $\sigma$.
Each $Q_{u}$ is orthogonal to the $xy$-plane.
Thus the intersection line $L_u := Q_u \cap \{z=0\}$ is at the same time the orthogonal projection of $Q_u$ to the $xy$-plane.
It follows that for every $u_0 \in U$ the curve $\tau|_{u=u_0}$ is contained in the line $L_{u_0}$.
Let $\xi(u)$ be a vector parallel to the line $L_u$.
Then one has $L_u = \{\tau(u,0) + t\xi(u) \mid t \in \R\}$ and
\[
\tau(u,v) = \gamma(u) + f(u,v)\xi(u),\quad \text{where }\gamma(u) = \tau(u,0),\ f(0,v) = 0.
\]
Our goal now is to prove that $f(u,v) = h(u)k(v)$, because then we get
\[
\tau(u,v) = \gamma(u) + f_1(v)\xi_1(u)\quad \text{for }f_1(v) = k(v) \text{ and }\xi_1(u) = h(u)\xi(u).
\]
One has
\[
\tau_u = \dot\gamma + f_u\xi + f\dot\xi, \quad \tau_{uv} = f_{uv}\xi + f_v\dot\xi.
\]
Equation \eqref{eqn:SigmaUParallel} is equivalent to $\tau_{uv} = \lambda \tau_u$, which implies that for every $v \in V$ the vector $\tau_u(u,v)$ is a scalar multiple of the vector $\tau_u(u,0) = \dot\gamma(u)$.
It follows that for every $(u,v) \in U\times V$ the vectors
\[
f_u\xi + f \dot\xi \quad \text{and} \quad f_{uv}\xi + f_v\dot\xi
\]
are scalar multiples of $\dot\gamma$.
In particular one has
\[
\det(f_u\xi + f \dot\xi, f_{uv}\xi + f_v\dot\xi) = 0 \Rightarrow (f_uf_v - ff_{uv}) \det(\xi,\dot\xi) = 0 \quad\text{for all }u,v.
\]
If $\det(\xi,\dot\xi) = 0$ on some interval $U' \subset U$, then the vector field $\xi$ has a constant direction on this interval.
This direction must be parallel to $\dot\gamma$, which implies that $\gamma|_{U'}$ is constant, and all lines $L_u$ for $u \in U'$ coincide.
This contradicts to the fourth condition in Definition \ref{dfn:THedron}.
Thus $\det(\xi,\dot\xi) \ne 0$ on a dense subset of $U$.
It follows that
\[
f_uf_v - ff_{uv} = 0 \quad\text{everywhere in }U \times V.
\]
This equation implies that
\[
(\log f)_{uv} = 0 \quad\text{whenever }f \ne 0,
\]
which, in turn, implies that every $(u_0,v_0) \in U \times V$ such that $f(u_0,v_0) \ne 0$ has a rectangular neighborhood $U' \times V'$ such that
\[
f(u,v) = h(u)k(v) \quad \text{for all }(u,v) \in U' \times V'.
\]
For a fixed interval $V'$ take the maximal interval $U'$ such that $f$ does not vanish on $U' \times V'$.
Then the above factorization of $f$ holds on $U' \times V'$ and even on the closure of this quadrilateral.
If $U' \ne U$, then by maximality of $U'$ one has $f(u_1,v_1)=0$ for some $v_1 \in V'$ and some $u_1$ on the boundary of $U'$.
It follows that $h(u_1) = 0$.
But then $h$, and therefore $f$, vanishes on the interval $\{u_1\} \times V'$, which implies that the coordinate curve $\sigma|_{u=u_1}$ contains a straight line interval.
Thus $U' = U$.
It follows that $V$ can be represented as the union of closed intervals $V = \cup_{\alpha \in A} V_\alpha$ with disjoint interiors such that
\[
f(u,v) = h_\alpha(u) k_\alpha(v) \quad \text{for all }(u,v) \in U \times V_\alpha
\]
with smooth functions $h_\alpha$ vanishing nowhere and smooth functions $k_\alpha$ vanishing at the endpoints of the intervals $V_\alpha$.
One thus has
\[
\tau_u(u,v) = \dot\gamma(u) + k_\alpha(v)(\dot h_\alpha(u)\xi(u) + h_\alpha(u)\dot\xi(u))
\]
and, by the above, the vector $\dot h_\alpha(u)\xi(u) + h_\alpha(u) \dot\xi(u)$ must be a scalar multiple of $\dot\gamma(u)$ for all $\alpha$ and for all $u$.
This implies that functions $h_\alpha$ are proportional to each other.
That is, there is a family of non-zero reals $\{\lambda_\alpha \mid \alpha \in A\}$ such that $h_\alpha = \lambda_\alpha h$ for some function $h$.
Replacing each $k_\alpha$ by $\lambda_\alpha k_\alpha$ and setting $k(v) := k_\alpha(v)$ for $v \in V_\alpha$ one obtains $f(u,v) = h(u)k(v)$ with $h$ nowhere vanishing.
Both $h$ and $k$ are $C^\infty$ because $f$ is.
This proves that $\sigma$ has the form \eqref{eqn:TSurfParam} with $\gamma(u) = \sigma(u,0)$ and $f(0) = 0$.
The properties stated in the theorem follow from the regularity of $\sigma$ and from the conditions in Definition \ref{dfn:TSurface}.
\end{proof}

\begin{lemma}
\label{lem:XiData}
Let
\[
\xi(u) = c(u)
\begin{pmatrix} \cos\phi(u)\\ \sin\phi(u) \end{pmatrix}.
\]
Then the condition $\dot\xi \parallel \dot\gamma$ is equivalent to
\begin{equation}
\label{eqn:CEta}
\dot c \cos\eta = c\dot\phi \sin\eta,
\end{equation}
where $\eta(u) = \phi(u) - \psi(u)$ is the angle from the right-hand normal of the curve $\gamma$ at the point $\gamma(u)$ to the vector $\xi(u)$.
\end{lemma}
\begin{proof}
Write
\[
\dot\gamma(u) = g(u)
\begin{pmatrix}
-\sin\psi(u)\\ \cos\psi(u)
\end{pmatrix},
\]
so that $\eta = \phi - \psi$.
One computes
\[
\dot\xi =
\dot c \begin{pmatrix} \cos\phi\\ \sin\phi \end{pmatrix} + c \dot \phi
\begin{pmatrix} -\sin\phi\\ \cos\phi \end{pmatrix}.
\]
It follows that
\[
\det\left(\dot\xi,\frac{\dot\gamma}{\|\dot\gamma\|}\right) =
\dot c \cos(\phi-\psi) - c \dot\phi \sin(\phi-\psi) = \dot c \cos\eta - c\dot\phi \sin\eta.
\]
Therefore the condition $\dot\xi \parallel \dot\gamma$ is equivalent to \eqref{eqn:CEta}.
\end{proof}

\begin{corollary}
A T-surface is uniquely determined by one trajectory curve, one profile curve, and the tangents to the profile curves at all points of the given trajectory curve. 
\end{corollary}
\begin{proof}
The profile curve $\sigma|_{u=0}$ determines the values of $f(v)\xi(0)$ and $z(v)$ for all $v$.
The trajectory curve $\sigma|_{v=0}$ is the curve $\gamma$; it also determines the function $\psi$.
The tangents to the profile curves at all points $(u,0)$ determine the functions $\phi$ and $\eta$.
Finally, equation \eqref{eqn:CEta} implies $\frac{\dot c}{c} = \dot\phi \tan\eta$, which allows to compute $\|\xi(u)\|$ up to a constant factor:
\[
c(u) = c(0) \exp\left(\int_0^u \dot\phi(w) \tan\eta(w)\, dw\right).
\]
Although $\xi$ is known up to a constant factor only, the values of $f(v)\xi(0)$ determine the product $f(v)\xi(u)$, and this is all which is needed.
\end{proof}

\subsection{Special classes of T-surfaces}
\subsubsection{Translational T-surfaces}
A T-surface is called \emph{translational} if all of its profile planes are parallel to each other.
The parallelity of profile planes is equivalent to $\phi = \const$.
Equation \eqref{eqn:CEta} then implies that $c = \const$, that is $\xi$ is a constant vector field.
Put the $x$-axis along the intersection line of $P_0$ and $Q_0$ and scale $f$ and $\xi$ so that $\xi = (1,0)$.
If $\gamma(u) = (x(u), y(u))$, then one has
\[
\sigma(u,v) =
\begin{pmatrix}
\gamma(u) + f(v)\xi\\ z(v)
\end{pmatrix}
= \begin{pmatrix} x(u) + f(v)\\ y(u)\\ z(v) \end{pmatrix} =
\begin{pmatrix} x(u)\\ y(u)\\ 0 \end{pmatrix} +
\begin{pmatrix} f(v)\\ 0\\ z(v) \end{pmatrix}.
\]

Translational T-surfaces are a subclass of translational surfaces generated by two arbitrary nowhere parallel spatial curves $\gamma$ and $\delta$ via $\sigma(u,v) = \gamma(u) + \delta(v)$.
In turn, the translational surfaces (more exactly, their nets of coordinate curves) form a subclass of Chebyshev nets, which are intensively studied in computational geometry and widely used in computer graphics.

\subsubsection{Molding surfaces}
Molding surfaces are non-translational T-surfaces with congruent profile curves, that is $c(u) = \const$ while $\dot\phi \ne 0$.
Equation \eqref{eqn:CEta} implies $\eta = 0$, which means that $\xi$ is the right-hand unit normal to $\gamma$ (again, the length of $\xi$ can be adjusted by scaling $f$ and $\xi$).
Thus molding surfaces are parametrized as
\[
\sigma(u,v) =
\begin{pmatrix}
\gamma(u) + f(v) n(u)\\ z(v)
\end{pmatrix}.
\]

\subsubsection{Axial T-surfaces}
A T-surface is called \emph{axial}, if all of its profile planes pass through a line.
This means that $\gamma(u) = \lambda(u)\xi(u)$.
Conditions $\dot\xi \parallel \dot\gamma$ and $\xi \not\parallel \dot\gamma$ imply that $\lambda$ is constant.
It follows that $\tau(u,v) = (\lambda + f(v))\xi(u)$.
By reassigning the notation $f$ to the function $\lambda + f$ one obtains the following parametrization of an axial T-surface:
\[
\sigma(u,v) =
\begin{pmatrix} f(v)\xi(u)\\ z(v) \end{pmatrix}.
\]
Here the function $f$ must have no zeros on $V$.
The fourth and the fifth conditions in Theorem \ref{thm:DescrTSurf} must hold with $\gamma$ replaced by $\xi$; the sixth condition is automatically satisfied.
Note that the fifth condition reads as $\xi$ being a regular curve not passing through the origin and such that $\dot\phi$, the derivative of the direction of the position vector of $\xi$, nowhere vanishes.
Equation \eqref{eqn:CEta} also holds automatically, with $\eta$ denoting the angle from the right-hand normal of $\xi$ to the position vector of $\xi$.

Axial molding T-surfaces are surfaces of revolution:
\[
\sigma(u,v) =
\begin{pmatrix}
f(v) \cos\phi(u)\\ f(v) \sin\phi(u)\\ z(v)
\end{pmatrix}.
\]

\section{Isometric deformations of T-surfaces}
\begin{definition}
An \emph{isometric deformation} of a regular parametrized surface $\sigma \colon \Omega \to \R^3$ is a smooth family of regular parametrized surfaces
\begin{equation}
\label{eqn:IsomDeformSurf}
\sigma^t \colon \Omega \to \R^3, \quad t \in (-\epsilon, \epsilon), \quad \sigma^0 = \sigma, 
\end{equation}
such that for every smooth map $\gamma \colon [a,b] \to \Omega$ the lengths of the curves $\sigma \circ \gamma$ and $\sigma^t \circ \gamma$ are equal for all $t$.
We call an isometric deformation \emph{non-trivial} if no two $\sigma^{t_1}$, $\sigma^{t_2}$ are related by a rigid motion: there is no isometry $\Phi \colon \R^3 \to \R^3$ such that $\sigma^{t_2} = \Phi \circ \sigma^{t_1}$.
\end{definition}
Often a weaker version of non-triviality is used: there is a $t$ such that $\sigma^t \ne \Phi \circ \sigma$.
If in \eqref{eqn:IsomDeformSurf} $t$ ranges in an interval $[0,\epsilon)$ or $(-\epsilon,0]$ only, then the deformation is called \emph{one-sided}.
We will see examples where only a one-sided deformation in a given class of surfaces is possible.

It is well-known that the condition on equal lengths of all corresponding curves is equivalent to the coincidence of the first fundamental forms:
\[
\begin{pmatrix}
E & F\\ F & G
\end{pmatrix}
=
\begin{pmatrix}
E^t & F^t\\ F^t & G^t
\end{pmatrix},
\]
where
\[
E(u,v) = \|\sigma_u(u,v)\|^2, \quad
F(u,v) = \langle \sigma_u(u,v), \sigma_v(u,v) \rangle, \quad
G(u,v) = \|\sigma_v(u,v)\|^2,
\]
and $E^t, F^t, G^t$ are defined similarly.

\subsection{Deformation of translational T-surfaces}
Isometric deformations of general translational surfaces $\sigma(u,v) = \gamma(u) + \delta(v)$ were thoroughly studied in the 19th century, see \cite[Section 26]{voss3abbildung}.
We derive here the formulas for the deformation a translational T-surface as a warm-up before applying similar methods to other classes of T-surfaces.

\begin{theorem}
\label{thm:DeformTranslSurf}
A translational T-surface
\[
\sigma(u,v) =
\begin{pmatrix}
x(u) + f(v)\\ y(u)\\ z(v)
\end{pmatrix}
\]
can be isometrically deformed in the class of translational T-surfaces provided that the derivatives of $y$ and $z$ nowhere vanish.
If, without loss of generality, one assumes $x(0) = f(0) = y(0) = z(0) = 0$ and $\dot y > 0$, $\dot z > 0$, then the deformation is given by
\[
\sigma^t(u,v) =
\begin{pmatrix}
e^{t}x(u) + e^{-t}f(v)\\
\int_0^u \sqrt{\dot y(w)^2 + (1-e^{2t})\dot x(w)^2}\, dw\\
\int_0^v \sqrt{\dot z(w)^2 + (1-e^{-2t}) \dot f(w)^2}\, dw.
\end{pmatrix}
\]
\end{theorem}

\begin{figure}[ht]
\begin{center}
\includegraphics[width=.28\textwidth]{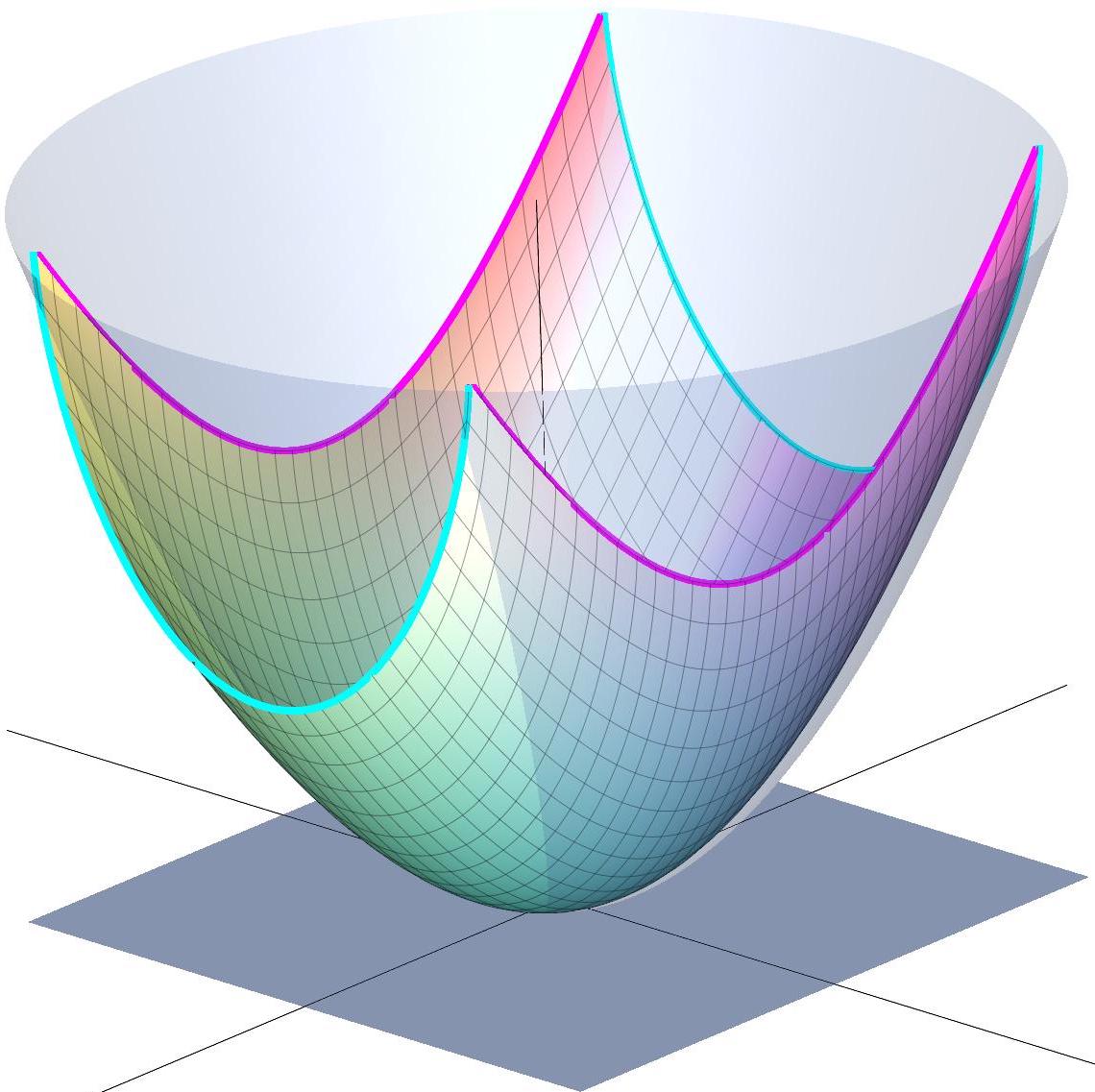}
\hspace{.04\textwidth}
\includegraphics[width=.28\textwidth]{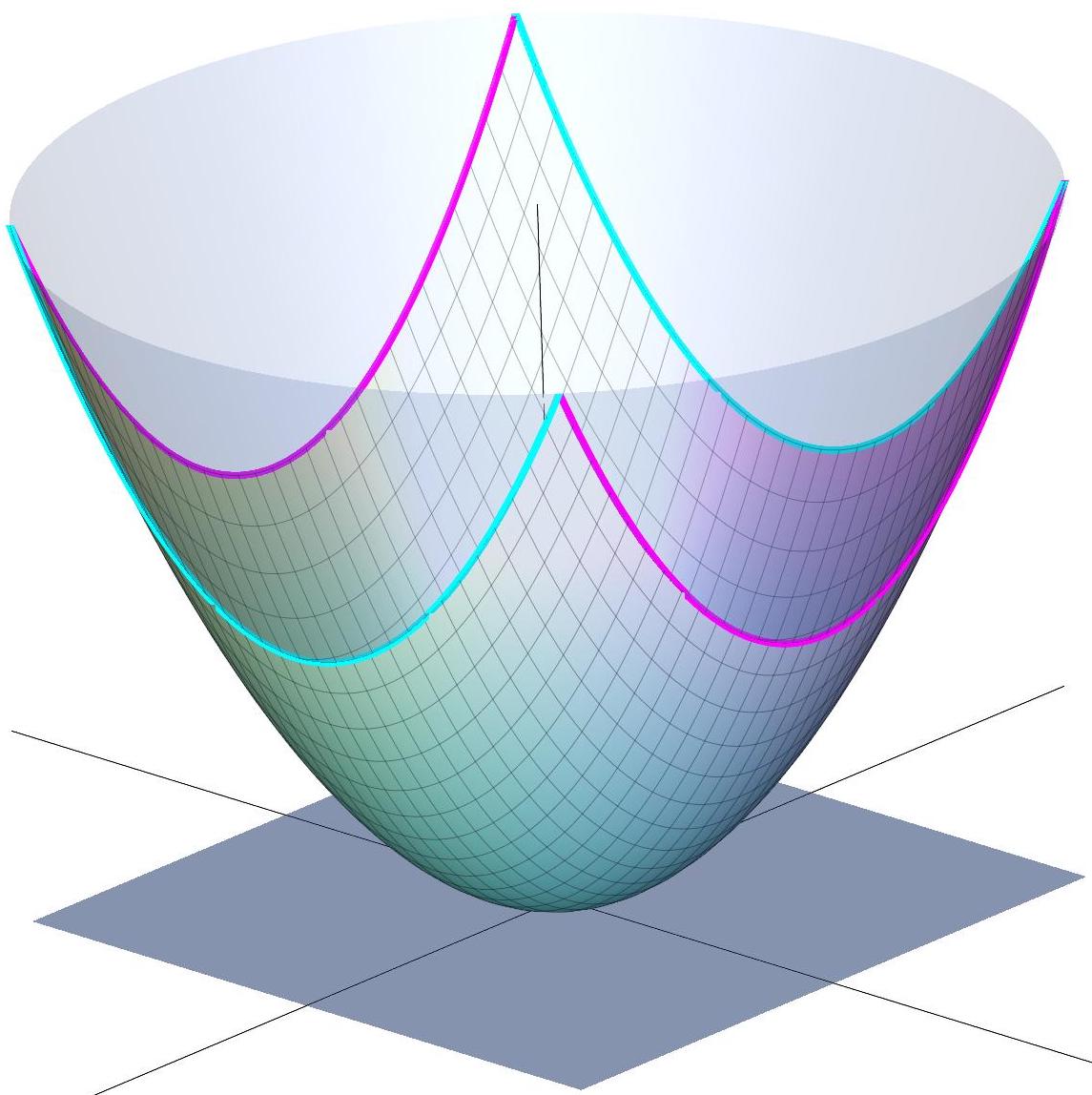}
\hspace{.04\textwidth}
\includegraphics[width=.28\textwidth]{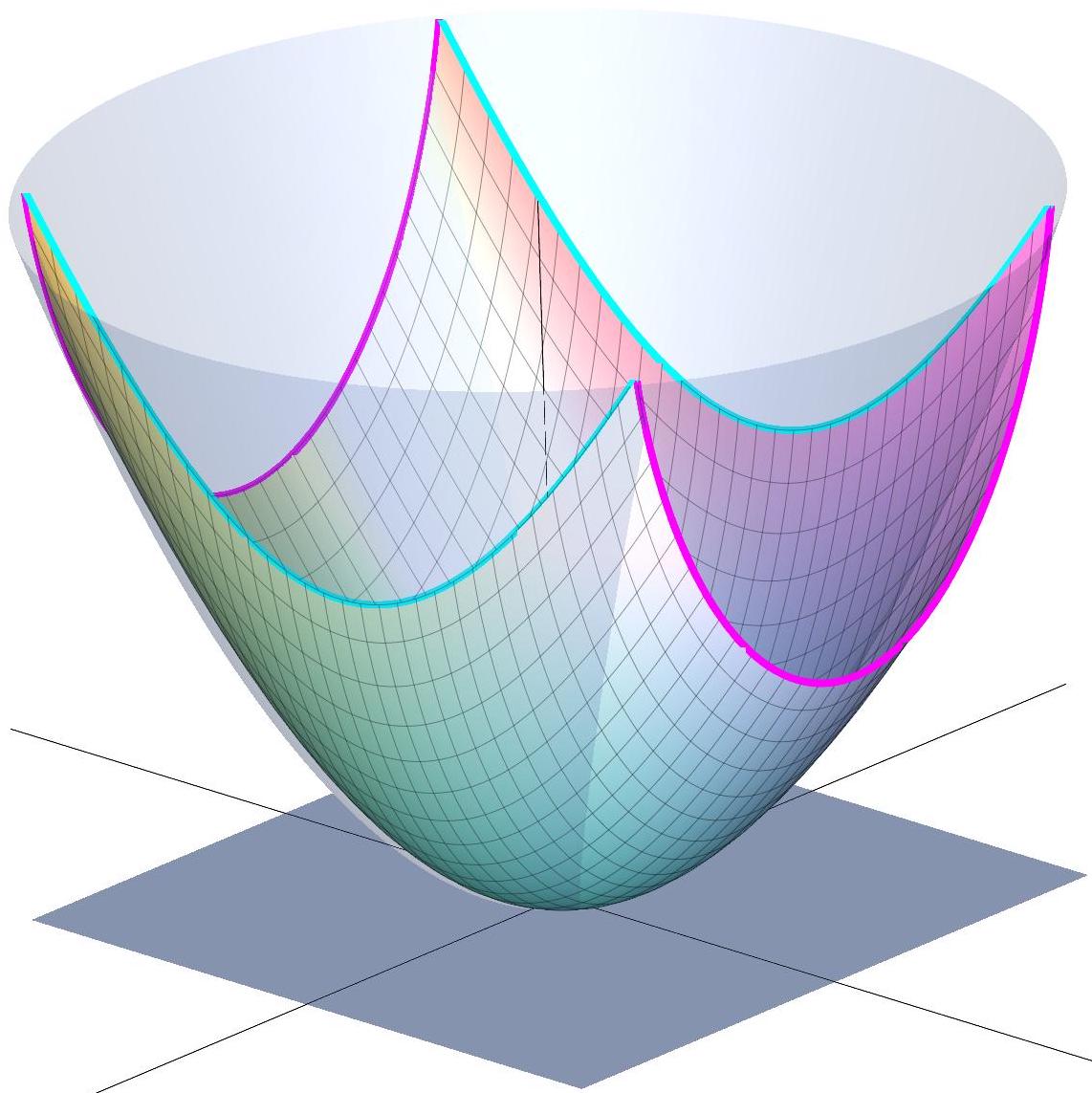}
\end{center}
\caption{An isometric deformation of a paraboloid of revolution as translational surface.}
\label{fig:DeformParaboloidTransl}
\end{figure}

\begin{example}
\label{exl:ParTranslSm}
The paraboloid of revolution $z = x^2 + y^2$ is a translation surface: the parabola $(u, 0, u^2)$ is translated along the parabola $(0, v, v^2)$.
The corresponding isometric deformation of the paraboloid is shown in Figure \ref{fig:DeformParaboloidTransl}.
Compare this with Example \ref{exl:ParTranslDiscr}.
\end{example}

\begin{proof}
From the partial derivatives
\[
\sigma_u(u,v) =
\begin{pmatrix}
\dot x(u)\\ \dot y(u)\\ 0
\end{pmatrix}, \quad
\sigma_v(u,v) =
\begin{pmatrix}
\dot f(v)\\ 0\\ \dot z(v)
\end{pmatrix}
\]
one computes the coefficients of the first fundamental form of $\sigma$:
\begin{gather*}
E = \dot x(u)^2 + \dot y(u)^2,\\
F = \dot x(u) \dot f(v),\\
G = \dot f(v)^2 + \dot z(v)^2.
\end{gather*}
Similar formulas hold for the first fundamental form of $\sigma^t$ determined by functions $x^t, f^t, y^t, z^t$ which are to be found.
One has
\[
F^t = F\ \Leftrightarrow\ \dot x^t(u) \dot f^t(v) = \dot x(u) \dot f(v)\ \Leftrightarrow\ \frac{\dot x^t(u)}{\dot x(u)} = \frac{\dot f(v)}{\dot f^t(v)}.
\]
In the last equation the left hand side is independent of $v$, the right hand side independent of $u$.
It follows that both expressions depend on $t$ only.
Thus one has
\[
F^t = F\ \Leftrightarrow\ \dot x^t(u) = k(t) \dot x(u), \ \dot f^t(v) = \frac{\dot f(v)}{k(t)} \text{ for some function }k.
\]
It follows that
\[
E^t = E\ \Leftrightarrow\ k(t)^2 \dot x(u)^2 + \dot y^t(u)^2 = \dot x(u)^2 + \dot y(u)^2\ \Leftrightarrow\ \dot y^t(u)^2 = \dot y(u)^2 + (1-k(t)^2) \dot x(u)^2.
\]
The right hand side becomes negative for $\dot y(u) = 0$ and $k(t)>1$.
Similarly,
\[
G^t = G\ \Leftrightarrow\ \frac{\dot f(v)^2}{k(t)^2} + \dot z^t(v)^2 = \dot f(v)^2 + \dot z(v)^2\ \Leftrightarrow\ \dot z^t(v)^2 = \dot z(v)^2 + (1 - k(t)^{-2}) \dot f(v)^2.
\]
Here the right hand side becomes negative for $\dot z(v) = 0$ and $k(t)<1$.
For a two-sided isometric deformation the function $k \colon (-\epsilon, \epsilon) \to \R$ must be monotone and satisfy $k(0)=1$.
Thus a deformation exists only if $\dot y$ and $\dot z$ never vanish.
One may put $k(t) = e^t$, which results in the formulas stated in the lemma.
\end{proof}

\begin{figure}[ht]
\begin{center}
\includegraphics[width=.4\textwidth]{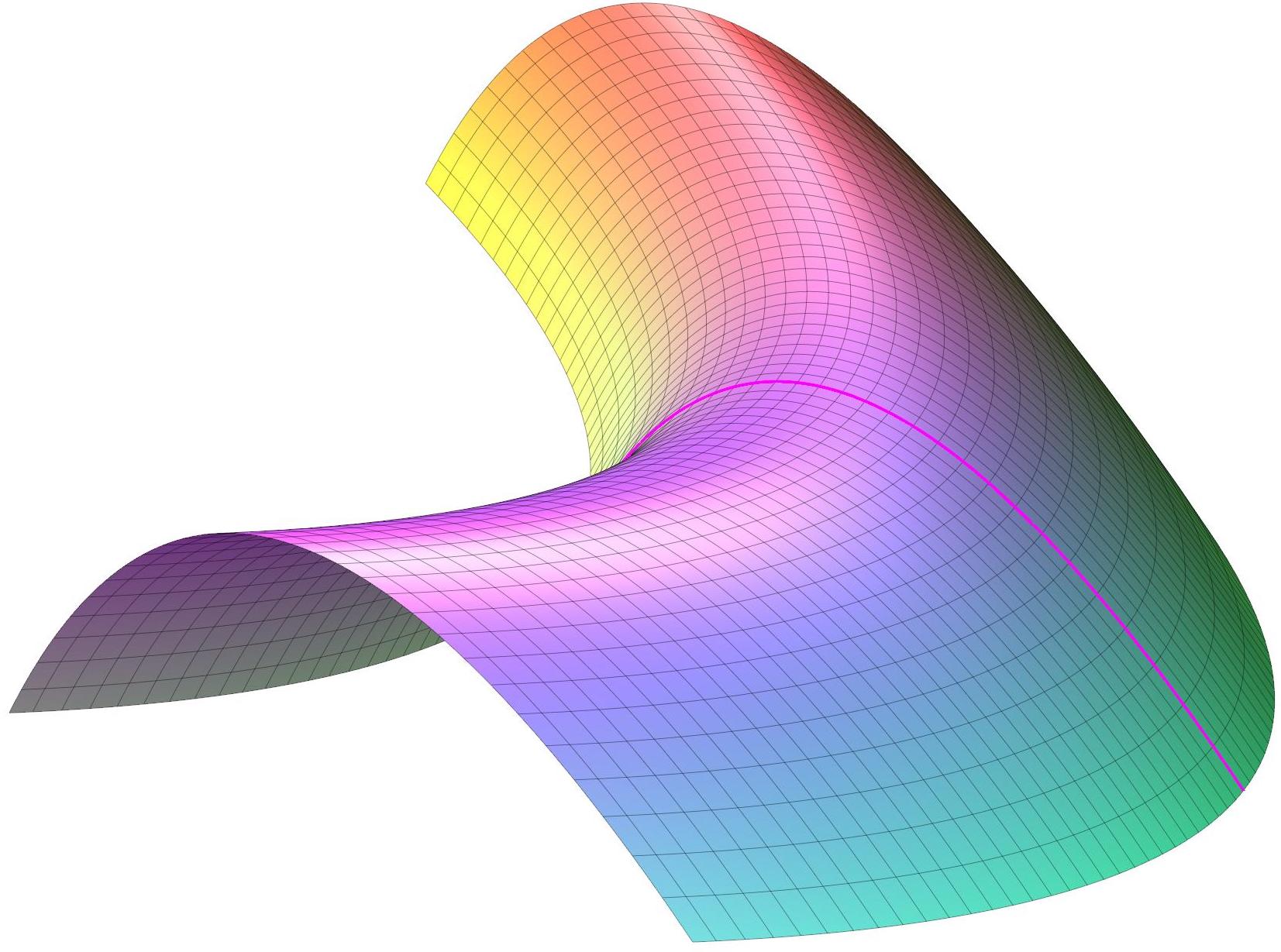}
\hspace{.1\textwidth}
\includegraphics[width=.4\textwidth]{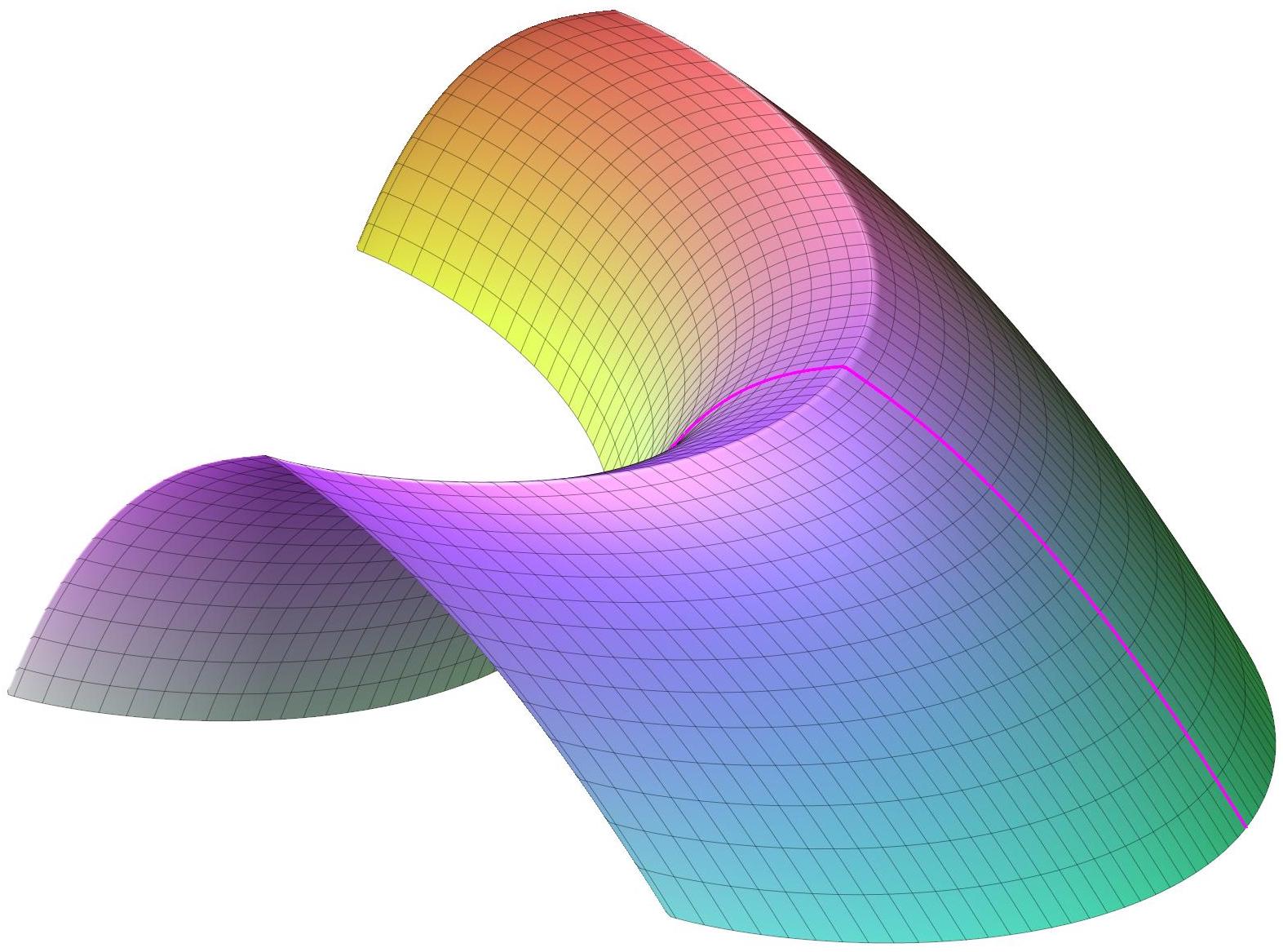}
\end{center}
\caption{A one-sided isometric deformation of a translational T-surface with a crease appearing along a coordinate curve.}
\label{fig:TranslSurfCreased}
\end{figure}

\begin{figure}[ht]
\begin{center}
\includegraphics[width=.4\textwidth]{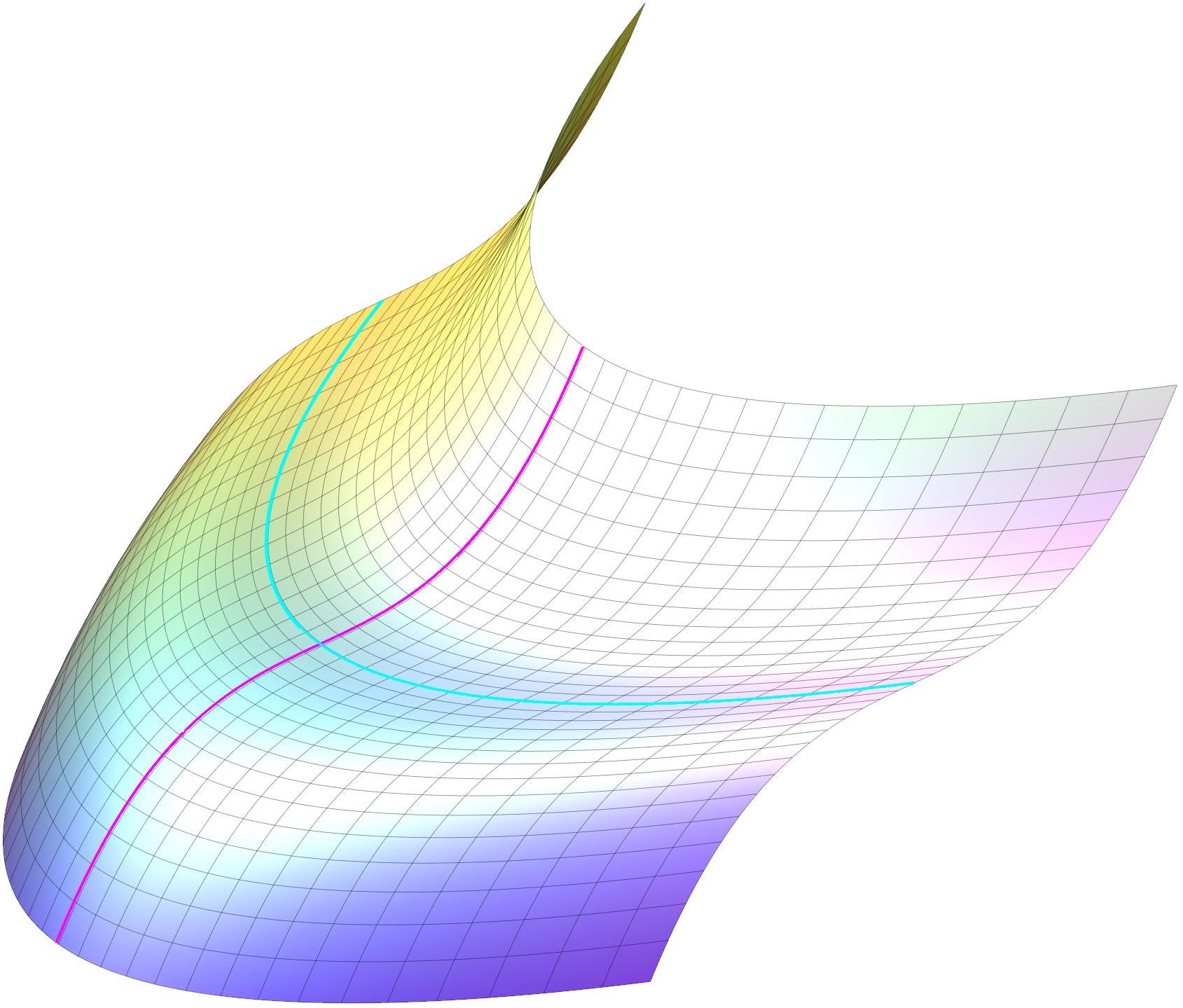}
\hspace{.1\textwidth}
\includegraphics[width=.4\textwidth]{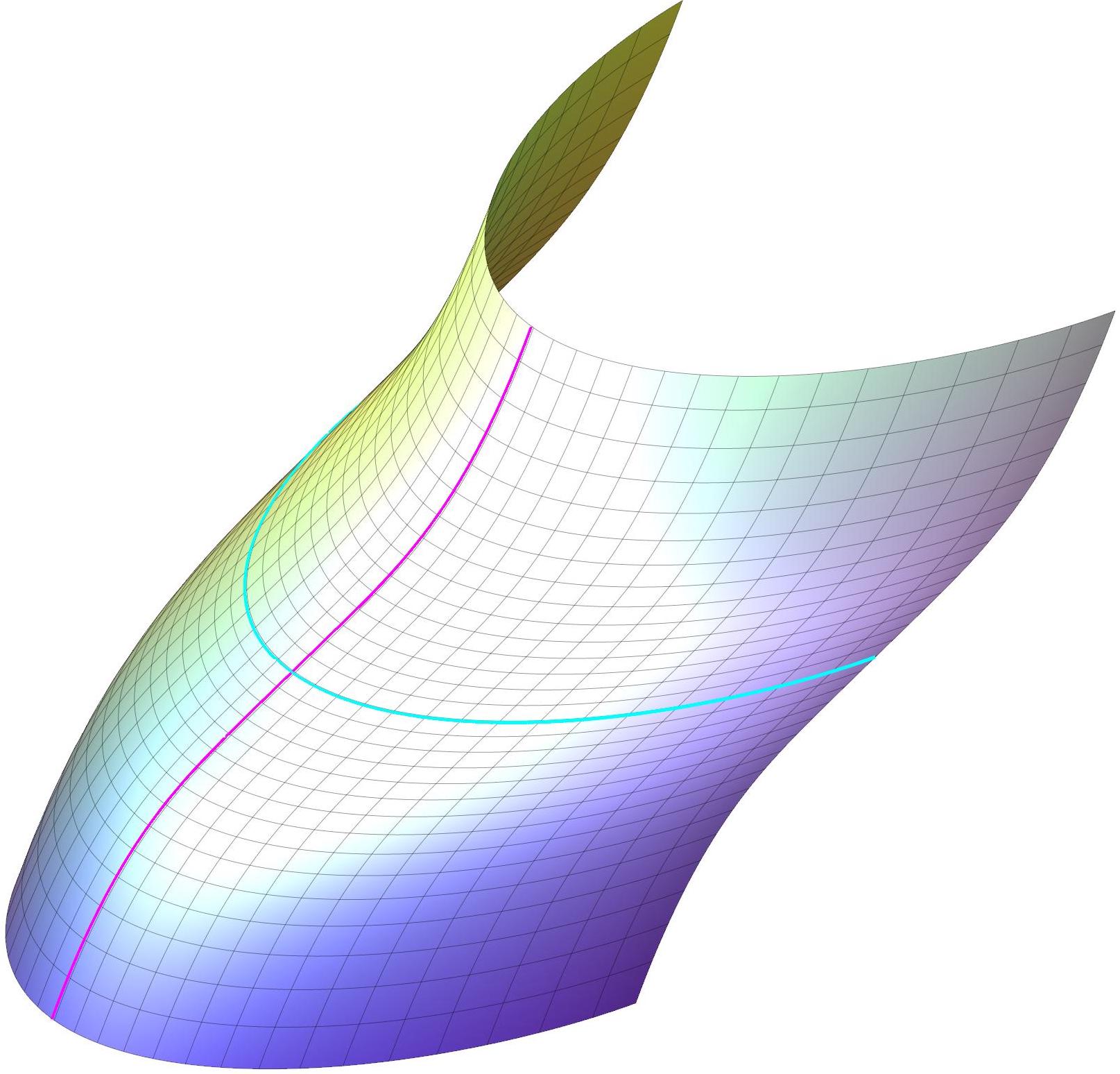}
\end{center}
\caption{A one-sided isometric deformation of a translational T-surface.}
\label{fig:TranslSurfNotCreased}
\end{figure}

If $\dot z(v_0) = 0$ for some $v_0$ but $\dot y$ nowhere vanishes, then the formulas in Theorem \ref{thm:DeformTranslSurf} yield a one-sided isometric deformation for $t \in [0,+\infty)$.
However, if $v_0$ is a local extremum, then the surface becomes creased along the coordinate curve $v=v_0$.
See Figure \ref{fig:TranslSurfCreased} for an illustration.
By contrary, $\sigma^t$ remains smooth for $t > 0$ if $v_0$ is an inflection point for the function $z$, see Figure \ref{fig:TranslSurfNotCreased}.
The same phenomenon occurs if $\dot f(u_0) = 0$ for some $u_0$.

\subsection{Deformation of molding surfaces}
\begin{theorem}
A molding surface
\[
\sigma(u,v) =
\begin{pmatrix}
\gamma(u) + f(v) n(u)\\ z(v)
\end{pmatrix}
\]
can be isometrically deformed in the class of molding surfaces if the derivative of $z$ nowhere vanishes.
If one assumes without loss of generality that $\gamma(0) = (0,0)$ and $\dot z > 0$, then the deformation has the form
\[
\sigma^t(u,v) =
\begin{pmatrix}
\gamma^t(u) + e^{-t}f(v) n^t(u)\\ z^t(v)
\end{pmatrix},
\]
where
\begin{gather*}
\gamma^t(u) = \int_0^u g(w)
\begin{pmatrix}
-\sin(e^t\psi(w))\\ \cos(e^t\psi(w))
\end{pmatrix}\, dw, \quad
n^t(u) =
\begin{pmatrix}
\cos(e^t\psi(u))\\ \sin(e^t\psi(u))
\end{pmatrix},\\
z^t(v) = \int_0^v \sqrt{\dot z(w)^2 + (1-e^{-2t})\dot f(w)^2}\, dw
\end{gather*}
with $g(u)$ the norm of $\dot\gamma(u)$ and $\psi(u)$ the direction of the right-hand normal to $\gamma$, in other words,
\[
\dot\gamma(u) = g(u)
\begin{pmatrix} -\sin\psi(u)\\ \cos\psi(u) \end{pmatrix}, \quad
n(u) =
\begin{pmatrix} \cos\psi(u)\\ \sin\psi(u) \end{pmatrix}.
\]
\end{theorem}
\begin{proof}
Compute the partial derivatives
\[
\sigma_u(u,v) =
\begin{pmatrix}
\dot\gamma(u) + f(v)\dot n(u)\\ 0
\end{pmatrix}
=
\begin{pmatrix}
-(g(u) + f(v)\dot\psi(u)) \sin\psi(u)\\
(g(u) + f(v)\dot\psi(u)) \cos\psi(u)\\
0
\end{pmatrix}, \quad
\sigma_v(u,v) =
\begin{pmatrix}
\dot f(v) n(u)\\ \dot z(v)
\end{pmatrix}
=
\begin{pmatrix}
\dot f(v) \cos\psi(u)\\
\dot f(v) \sin\psi(u)\\
\dot z(u)
\end{pmatrix}.
\]
From this one computes the coefficients of the first fundamental form:
\[
E = (g(u) + f(v)\dot\psi(u))^2, \quad F = 0, \quad G = \dot f(v)^2 + \dot z(v)^2.
\]
Similar formulas hold for $\sigma^t$.
One has
\[
E^t = E\ \Leftrightarrow\ g^t(u) + f^t(v)\dot\psi^t(u) = g(u) + f(v)\dot\psi(u).
\]
Substituting $v=0$ one gets $g^t = g$.
Thus one has
\[
\frac{f^t(v)}{f(v)} = \frac{\dot\psi(u)}{\dot\psi^t(u)}.
\]
In this equation the left hand side is independent of $u$, the right hand side independent of $v$.
It follows that both expressions depend on $t$ only.
Reparametrizing the deformation one may assume that
\[
f^t(v) = e^{-t}f(v), \quad \dot\psi^t(u) = e^t \dot\psi(u).
\]
By applying, if needed, a rigid motion, one may assume $\psi^t(0) = 0$ for all $t$, so that
\[
\dot\psi^t(u) = e^t \dot\psi(u)\ \Rightarrow\ \psi^t(u) = e^t \psi(u).
\]
Then one has
\[
\dot\gamma^t(u) = g(u)
\begin{pmatrix}
-\sin(e^t\psi(u))\\ \cos(e^t\psi(u))
\end{pmatrix}, \quad
n^t(u) =
\begin{pmatrix} \cos(e^t\psi(u))\\ \sin(e^t\psi(u)) \end{pmatrix},
\]
which implies the formulas for $\gamma^t(u)$ and $n^t(u)$ stated in the lemma (again, with a rigid motion one can achieve $\gamma^t(0) = (0,0)$ and $\psi^t(0) = 0$).

It remains to observe that
\[
G^t = G\ \Leftrightarrow\ e^{-2t}\dot f(v)^2 + \dot z^t(v)^2 = \dot f(v)^2 + \dot z(v)^2\ \Leftrightarrow\ \dot z^t(v)^2 = \dot z(v)^2 + (1-e^{-2t})\dot f(v)^2,
\]
and the theorem is proved.
\end{proof}

If the derivative of $z$ vanishes at some points, then there is a one-sided deformation where the surface gets creased along the trajectory curves corresponding to the local extrema of $z$, similarly to the situation for translational T-surfaces depicted in Figure \ref{fig:TranslSurfCreased}.

\subsection{Deformation of axial T-surfaces}
\begin{theorem}
\label{thm:DeformAxialSurf}
An axial T-surface
\[
\sigma(u,v) =
\begin{pmatrix}
f(v)\xi(u)\\ z(v)
\end{pmatrix}
=
\begin{pmatrix}
f(v)c(u) \cos\phi(u)\\ f(v)c(u) \sin\phi(u)\\ z(v)
\end{pmatrix}
\]
can be isometrically deformed in the class of axial T-surfaces if the derivative of $z$ nowhere vanishes.
If one assumes without loss of generality that $\dot \phi > 0$ and $\dot z > 0$, then the deformation has the form
\[
\sigma^t(u,v) =
\begin{pmatrix}
f(v)\sqrt{c(u)^2 + t} \cos\phi^t(u)\\
f(v)\sqrt{c(u)^2 + t} \sin\phi^t(u)\\
z^t(v)
\end{pmatrix},
\]
where
\begin{gather*}
\phi^t(u) = \int_0^u \frac{\sqrt{c(w)^4\dot\phi(w)^2 + t(c(w)^2\dot\phi(w)^2 + \dot c(w)^2)}}{c(w)^2 + t}\, dw = \int_0^u \frac{\dot\phi(w)c(w)\sqrt{c(w)^2 + \frac{t}{\cos^2\eta(w)}}}{c(w)^2 + t}\, dw,\\
z^t(v) = \int_0^v \sqrt{\dot z(w)^2 - t \dot f(w)^2}\, dw,
\end{gather*}
with $\eta$ denoting the angle from the right-hand normal of $\xi$ to the position vector of $\xi$.
\end{theorem}
\begin{proof}
The equality of the first fundamental forms of $\sigma$ and $\sigma^t$ can be checked by a direct computation.
We will demonstrate how the formulas are derived, which will also imply that there is only one degree of freedom within the space of axial T-surfaces isometric to $\sigma$.

From the partial derivatives
\[
\sigma_u(u,v) =
\begin{pmatrix}
f(v)\dot\xi(u)\\ 0
\end{pmatrix}, \quad
\sigma_v(u,v) =
\begin{pmatrix}
\dot f(v)\xi(u)\\ \dot z(v)
\end{pmatrix},
\]
one computes the coefficients of the first fundamental form of $\sigma$:
\begin{gather*}
E = f(v)^2 \|\dot\xi(u)\|^2,\\
F = \dot f(v) f(v) \langle \dot\xi(u), \xi(u) \rangle,\\
G = \dot f(v)^2 \|\xi(u)\|^2 + \dot z(v)^2.
\end{gather*}
Similar formulas hold for the first fundamental form of $\sigma^t$ determined by functions $f^t, c^t, \phi^t, z^t$ which are to be found.
One has
\[
E^t = E \Leftrightarrow f^t(v)\|\dot\xi^t(u)\| = f(v)\|\dot\xi(u)\| \Leftrightarrow \frac{f^t(v)}{f(v)} = \frac{\|\dot\xi(u)\|}{\|\dot\xi^t(u)\|}.
\]
In the last equation the left hand side is independent of $v$, the right hand side independent of $u$.
It follows that both expressions depend on $t$ only.
Since scaling $f^t$ up and $\xi^t$ down by the same factor does not change $\sigma^t$, one may assume
\[
f^t(u) = f(u), \quad \|\dot\xi^t(u)\| = \|\dot\xi(u)\|.
\]
The latter equation is equivalent to
\begin{equation}
\label{eqn:CTAxial}
(\dot c^t)^2 + (c^t)^2(\dot\phi^t)^2 = \dot c^2 + c^2 \dot\phi^2.
\end{equation}
Due to $f^t = f$ one has
\[
F^t = F \Leftrightarrow \langle \dot\xi^t(u), \xi^t(u) \rangle = \langle \dot\xi(u), \xi(u) \rangle \Leftrightarrow \frac{d(c^t(u)^2)}{du} = \frac{d(c(u)^2)}{du}.
\]
That is to say, $c^t(u)^2 - c(u)^2$ must be independent of $u$.
By changing the deformation parameter $t$ one can assume
\[
c^t(u)^2 = c(u)^2 + t.
\]
This implies
\[
c^t = \sqrt{c^2 + t}, \quad \dot c^t = \frac{\dot c c}{\sqrt{c^2 + t}}.
\]
Substituting this into \eqref{eqn:CTAxial} and solving the resulting equation for $\dot\phi^t$ one obtains
\[
(\dot\phi^t)^2 = \frac{c^4\dot\phi^2 + t(c^2\dot\phi^2 + \dot c^2)}{(c^2+t)^2}.
\]
Since the derivative of $\phi$ is positive, the right hand side remains positive for all $t$ sufficiently close to $0$.
Postcomposing $\sigma^t$, if needed, with a rigid motion of $\R^3$, one may assume that $\phi^t(0) = \phi(0)$.
It follows that $\phi^t(u) = \phi(0) + \int_0^w \dot\phi^t(w)\, dw$, so that one obtains the first formula for $\phi^t$ stated in the lemma.
This formula can be simplified using the equation \eqref{eqn:CEta} to obtain the second formula for $\phi^t$.

Finally, due to $f^t = f$ and $\|\xi^t(u)\|^2 = \|\xi(u)\|^2 + t$ one has
\[
G^t = G \Leftrightarrow f(v)^2(c(u)^2 + t) + \dot z^t(v)^2 = f(v)^2c(u)^2 + \dot z(v)^2 \Leftrightarrow \dot z^t(v)^2 = \dot z(v)^2 - tf(v)^2.
\]
Up to a rigid motion, $z^t(0) = z(0)$, which leads to the last formula of the theorem.
\end{proof}

\begin{figure}[ht]
\begin{center}
\includegraphics[width=.31\textwidth]{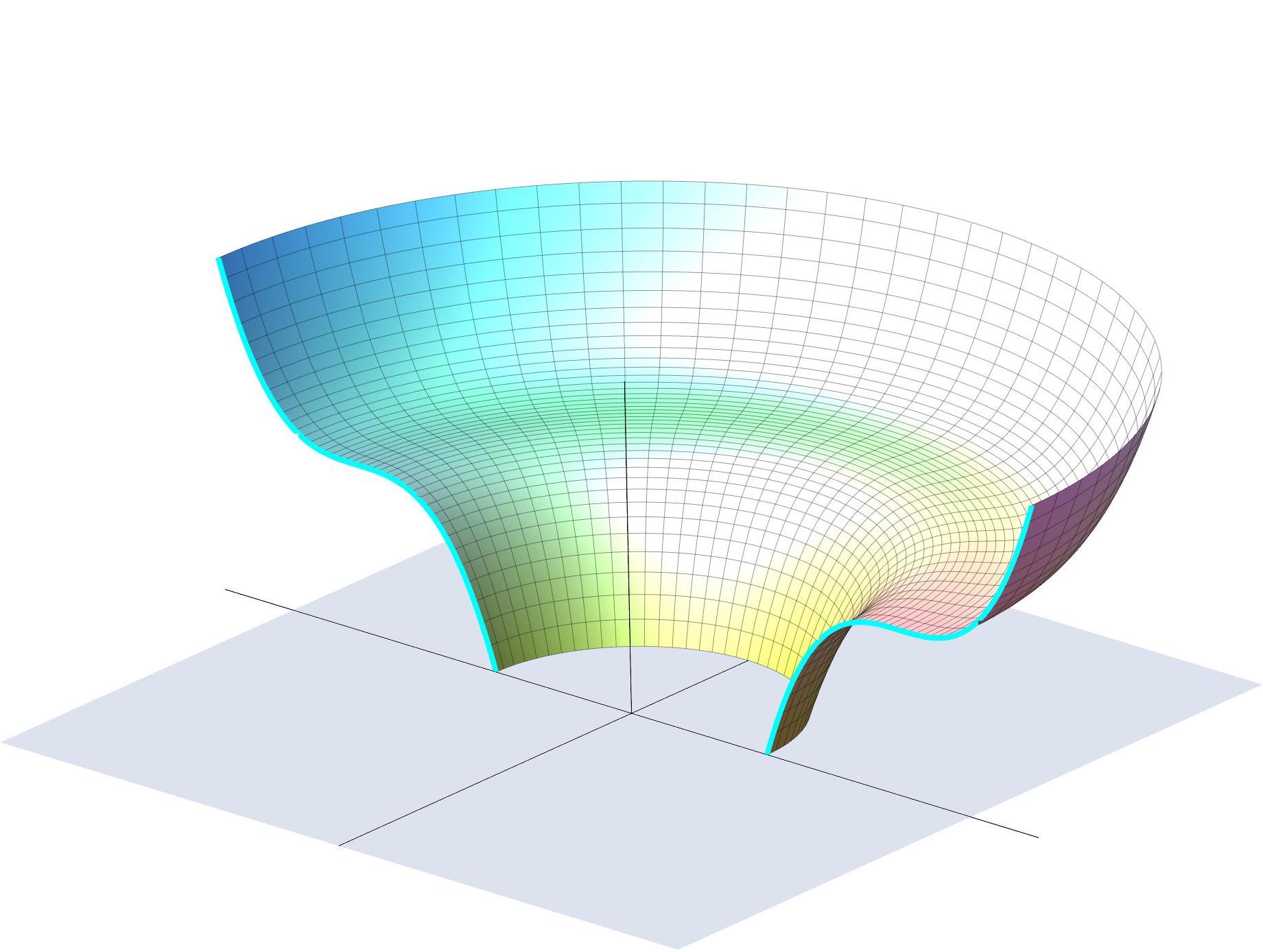}
\includegraphics[width=.31\textwidth]{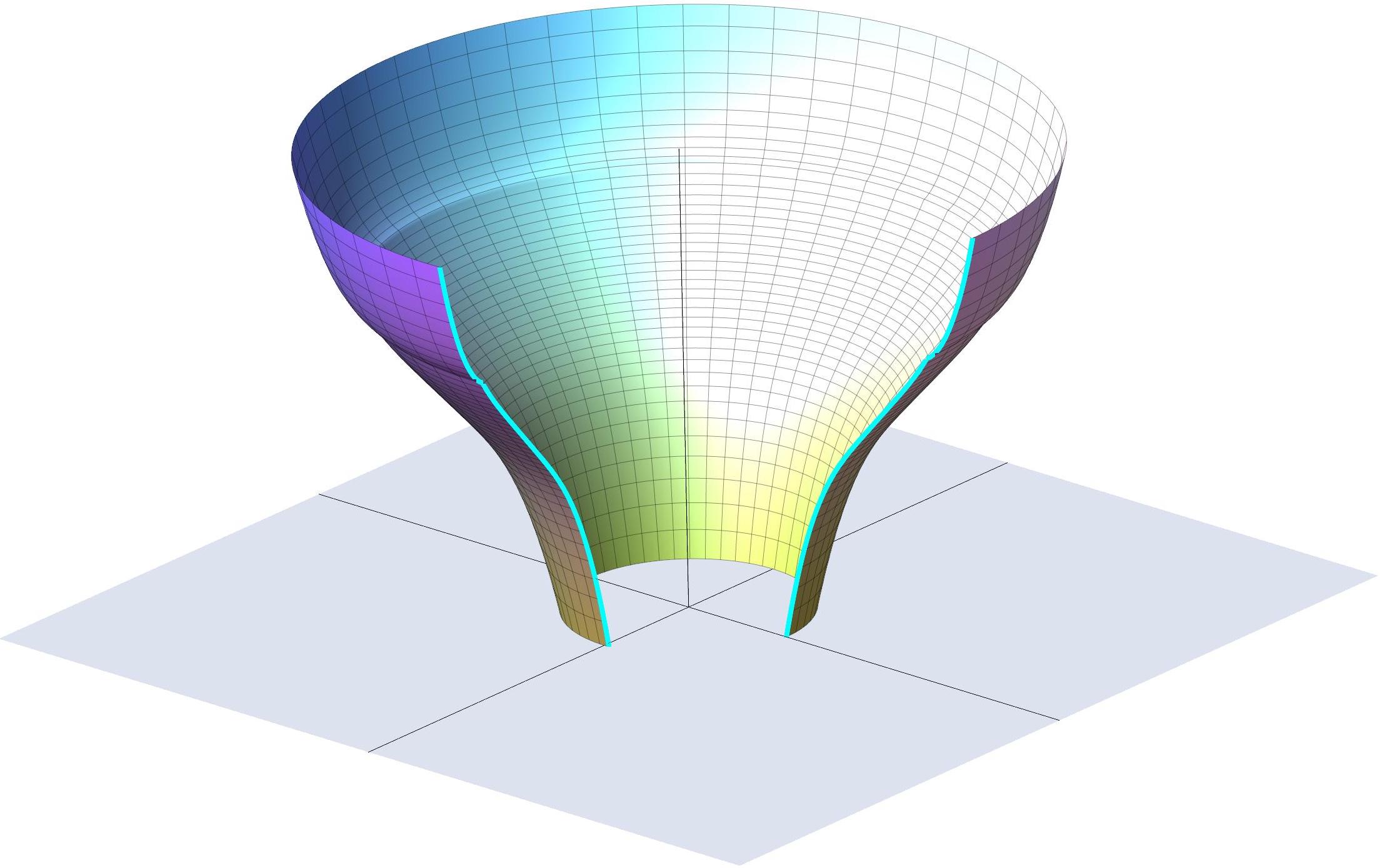}
\includegraphics[width=.31\textwidth]{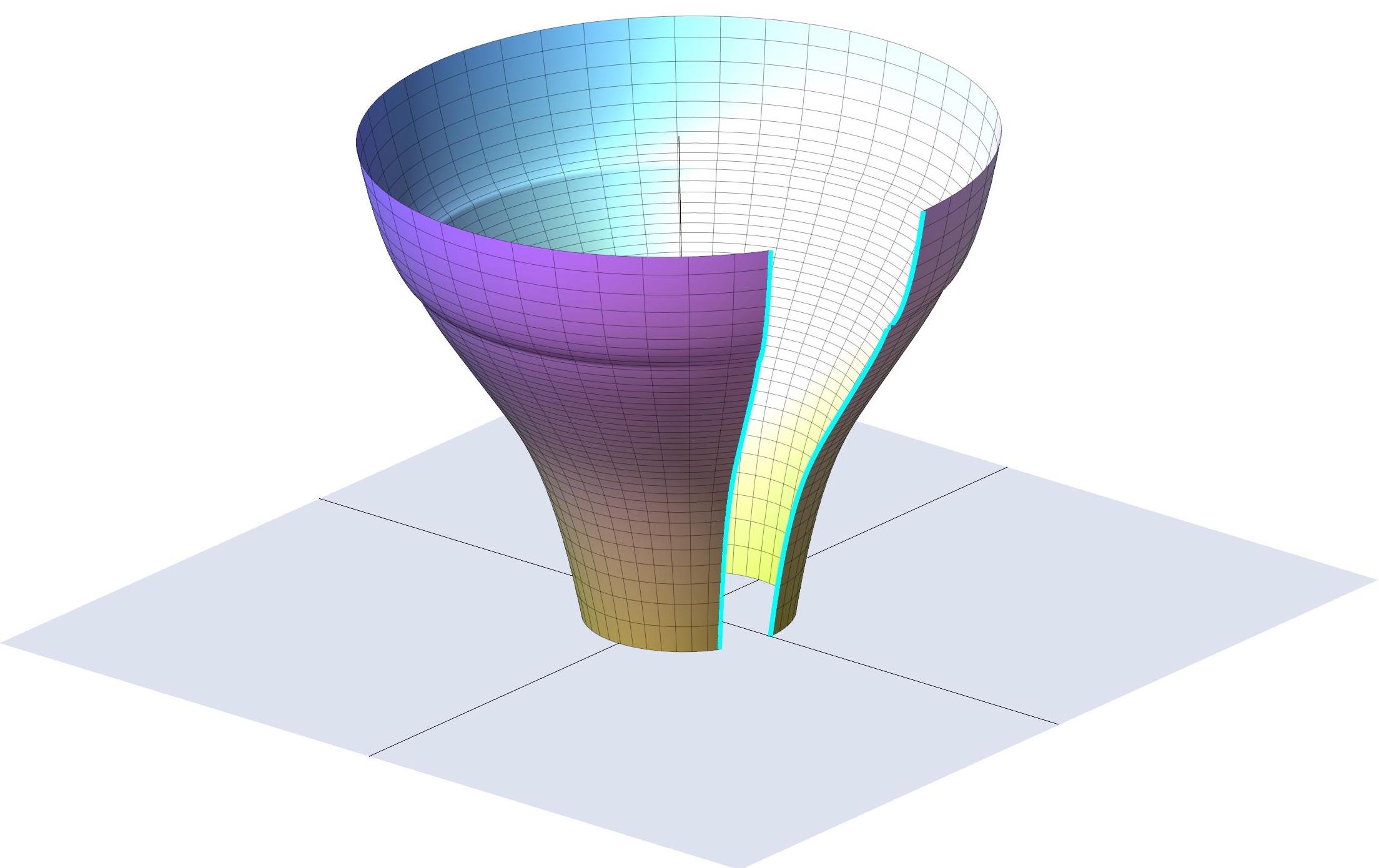}
\end{center}
\caption{A one-sided isometric deformation of a surface of revolution whose profile has an inflection point with a horizontal tangent.}
\label{fig:DeformAxial}
\end{figure}

If the derivative of $z$ vanishes at some points, then there is a one-sided deformation where the surface gets creased along the trajectory curves corresponding to the local extrema of $z$, similarly to the situation for translational T-surfaces depicted in Figure \ref{fig:TranslSurfCreased}.
If $z$ has no local extrema but only inflection points, then it remains smooth during the one-sided deformation, see Figure \ref{fig:DeformAxial}.

By specializing Theorem \ref{thm:DeformAxialSurf} to $c(u)=1$ one obtains formulas for an isometric deformation of surfaces of revolution.

\begin{theorem}
A surface of revolution
\[
\sigma(u,v) =
\begin{pmatrix}
f(v) \cos\phi(u)\\ f(v) \sin\phi(u)\\ z(v)
\end{pmatrix}
\]
can be isometrically deformed in the class of surfaces of revolution if the derivative of $z$ nowhere vanishes.
If one assumes without loss of generality that $\dot \phi > 0$ and $\dot z > 0$, then the deformation has the form
\[
\sigma^t(u,v) =
\begin{pmatrix}
f(v)\sqrt{1 + t} \cos\frac{\phi(u)}{\sqrt{1+t}}\\
f(v)\sqrt{1 + t} \sin\frac{\phi(u)}{\sqrt{1+t}}\\
\int_0^v \sqrt{\dot z(w)^2 - t \dot f(w)^2}\, dw
\end{pmatrix}.
\]
\end{theorem}

\begin{figure*}[ht]
    \begin{subfigure}[b]{.3\columnwidth}
    \centering
        \includegraphics[height= 45 mm]{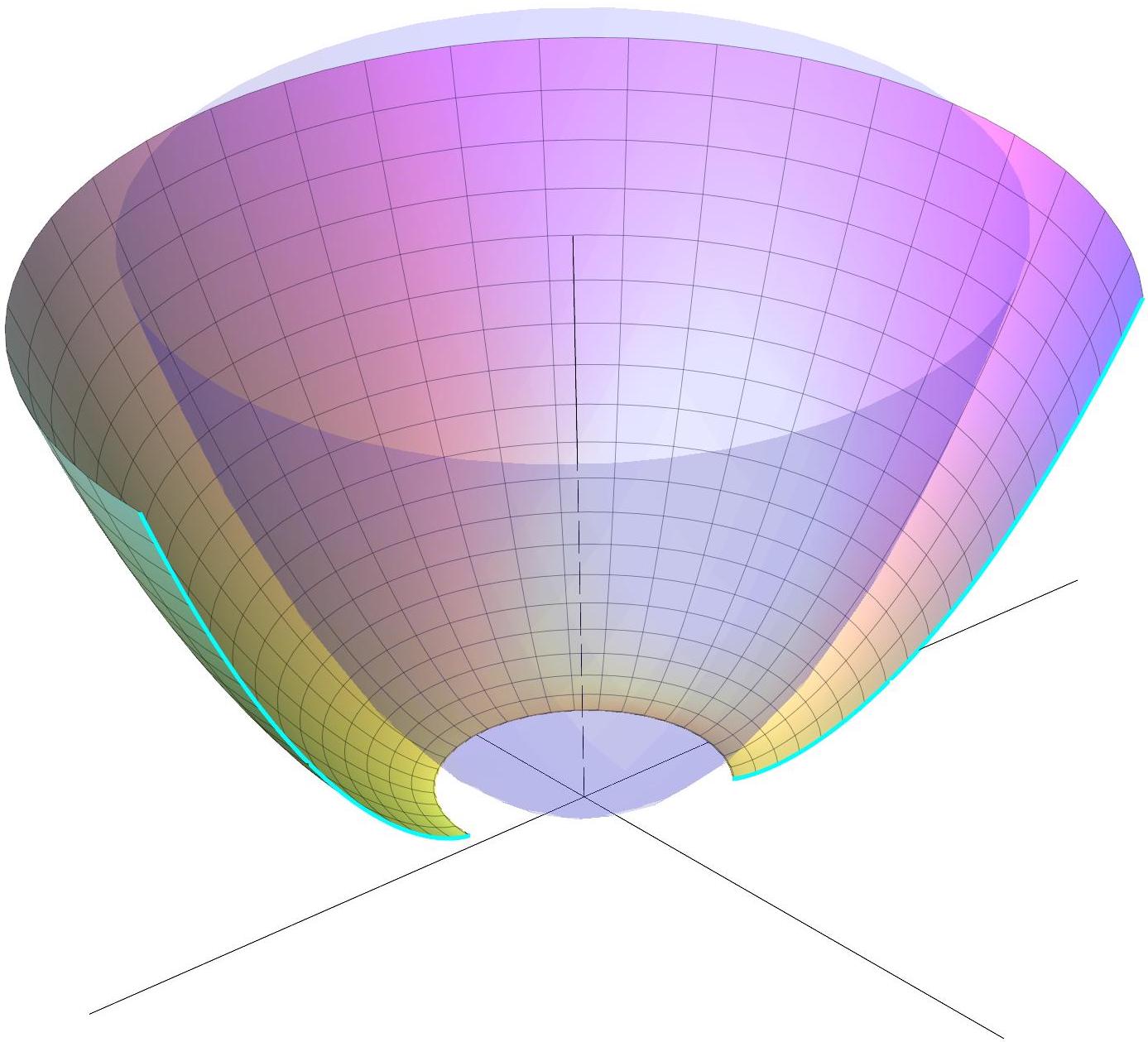}
    \end{subfigure}
    \hfill
    \begin{subfigure}[b]{.3\columnwidth}
    \centering
        \includegraphics[height= 45 mm]{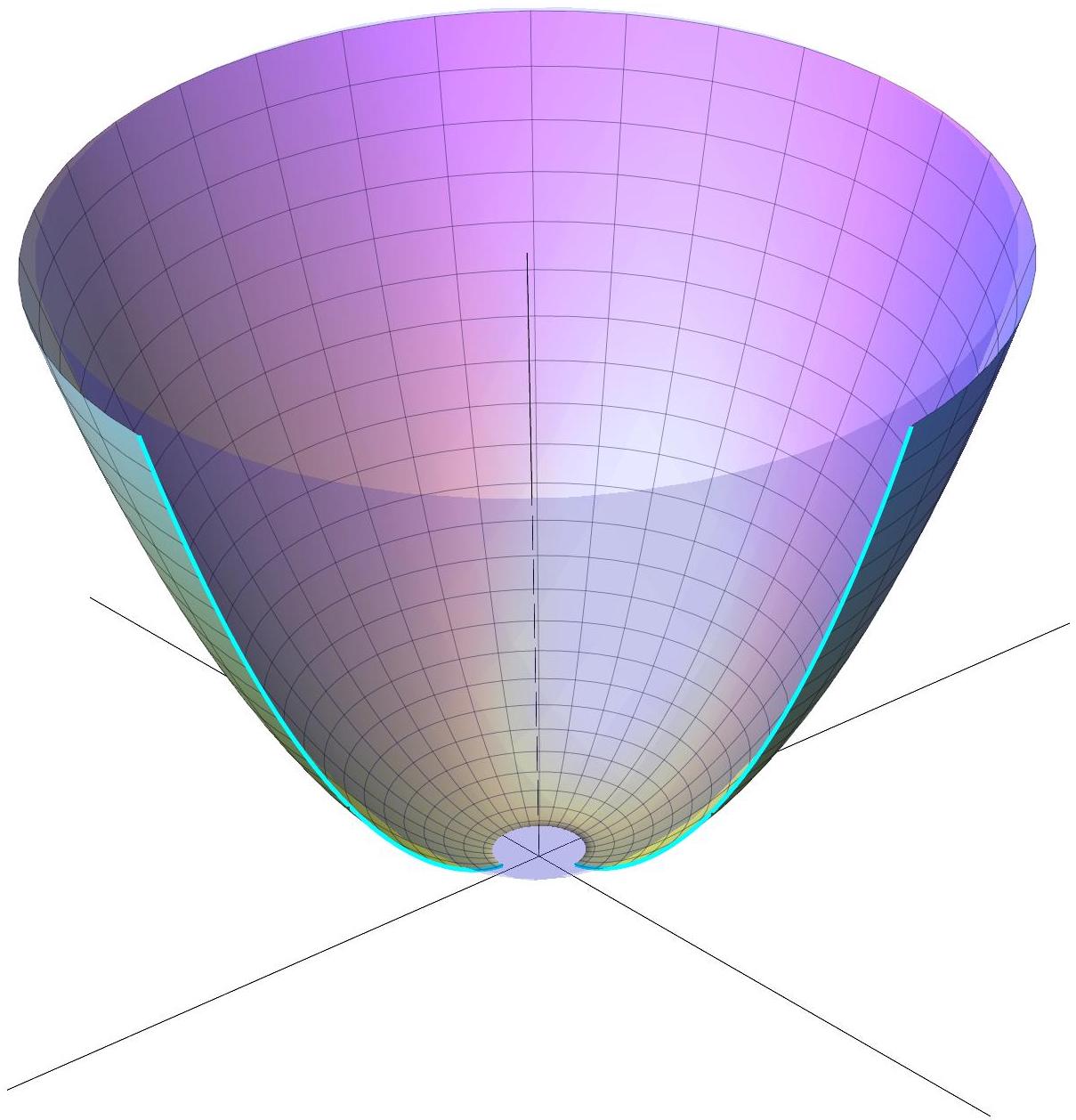}
    \end{subfigure}
    \hfill
    \begin{subfigure}[b]{.3\columnwidth}
    \centering
        \includegraphics[height= 45 mm]{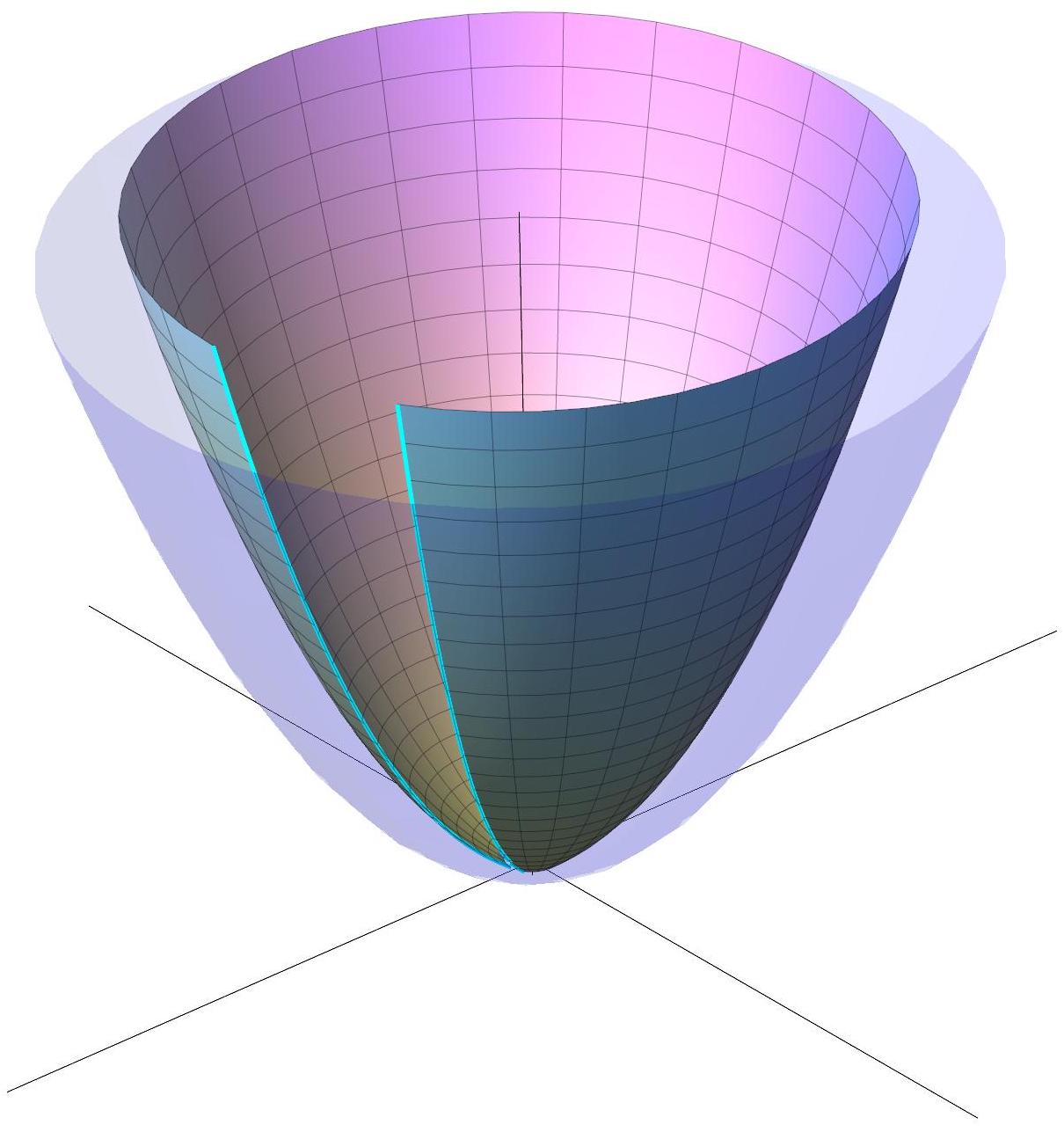}
    \end{subfigure}
    \caption{Isometric deformation of a wedge of a paraboloid of revolution in the class of surfaces of revolution.}
    \label{fig:DeformParaboloidRev}
\end{figure*}

Figure \ref{fig:DeformParaboloidRev} shows an isometric deformation of a wedge of the paraboloid $z = x^2 + y^2$.
Compare this with Figure \ref{fig:DeformParabDiscrRev}.

\subsection{Deformations of general T-surfaces}
\begin{theorem}
A T-surface
\[
\sigma(u,v) =
\begin{pmatrix}
\gamma(u) + f(v) \xi(u)\\
z(v)
\end{pmatrix}, \quad
\gamma(u) = \int_0^u g(u)
\begin{pmatrix} -\sin\psi(w)\\ \cos\psi(w) \end{pmatrix}, \quad
\xi(u) = c(u)
\begin{pmatrix} \cos\phi(u)\\ \sin\phi(u) \end{pmatrix}
\]
can be isometrically deformed in the class of T-surfaces if the derivatives of $\phi$ and $z$ nowhere vanish.
If one assumes without loss of generality that $\phi(0) = z(0) = 0$, $\dot\phi > 0$ and $\dot z > 0$, then the deformation has the form
\[
\sigma^t(u,v) =
\begin{pmatrix}
\gamma^t(u) + f(v) \xi^t(u)\\ z^t(u)
\end{pmatrix}, \quad
\gamma^t(u) = \int_0^u g(u)
\begin{pmatrix} -\sin\psi^t(w)\\ \cos\psi^t(w) \end{pmatrix}, \quad
\xi^t(u) = \sqrt{c(u)^2 + t}
\begin{pmatrix} \cos\phi^t(u)\\ \sin\phi^t(u) \end{pmatrix},
\]
where
\begin{gather*}
z^t(v) = \int_0^v \sign(\dot z(w)) \sqrt{\dot z(w)^2 - t \dot f(w)^2}\, dw,\\
\phi^t(u) = \int_0^u \frac{\dot\phi(w)c(w)\sqrt{c(w)^2 + \frac{t}{\cos^2\eta(w)}}}{c(w)^2 + t}\, dw, \quad\text{where } \eta(w) = \phi(w) - \psi(w)\\
\psi^t(u) = \phi^t(u) - \eta^t(u), \quad\text{where }
\eta^t(u) = \arctan \frac{c \sin\eta}{\sqrt{c^2\cos^2\eta + t}}.
\end{gather*}
\end{theorem}
\begin{proof}
From the partial derivatives:
\[
\sigma_u(u,v) =
\begin{pmatrix}
\dot\gamma(u) + f(v) \dot\xi(u)\\
0
\end{pmatrix}, \quad
\sigma_v =
\begin{pmatrix}
\dot f(v) \xi(u)\\ \dot z(v)
\end{pmatrix}
\]
one computes the coefficients of the first fundamental form:
\begin{gather*}
E = \|\sigma_u(u,v)\|^2 = \|\dot\gamma(u) + f(v)\dot\xi(u)\|^2\\
F = \langle \sigma_u(u,v), \sigma_v(u,v) \rangle =
\langle \dot\gamma(u) + f(v) \dot\xi(u), \dot f(v) \xi(u) \rangle\\
G = \|\sigma_v(u,v)\|^2 = (\dot f(v))^2 \|\xi(u)\|^2 + (\dot z(v))^2,
\end{gather*}
and similarly for $\sigma^t$.

Let us study the equation $E^t = E$.
Substituting $v=0$ one obtains
\begin{equation}
\label{eqn:GammaConst}
\|\dot\gamma^t(u)\| =  \|\dot\gamma(u)\|
\end{equation}
for all $t, u$.
Because of $\dot\gamma(u) \parallel \dot\xi(u)$ and $\dot\gamma^t(u) \parallel \dot\xi^t(u)$ the equation $\|\sigma^t_u(u,v)\| = \|\sigma_u(u,v)\|$ is equivalent to $f^t(v)\|\dot\xi^t(u)\| = f(v)\|\dot\xi(u)\|$.
This rewrites as
\[
\frac{f^t(v)}{f(v)} = \frac{\|\dot\xi(u)\|}{\|\dot\xi^t(u)\|}.
\]
Since the left hand side is independent of $u$, and the right hand side is independent of $v$, these quotients are equal to a function $a(t)$.
Remembering that the parametrization $\sigma^t$ does not change under scaling of $f^t$ down and $\xi^t$ up by $a(t)$ one concludes that without loss of generality it may be assumed that
\begin{equation}
\label{eqn:FXiConst}
f^t(v) = f(v), \quad \|\dot\xi^t(u)\| = \|\dot\xi(u)\| \quad\text{for all } t.
\end{equation}
Conditions \eqref{eqn:GammaConst} and \eqref{eqn:FXiConst} are necessary and sufficient for $E^t = E$, that is for an isometric deformation of all trajectory curves.

Now let us study the equation $F^t = F$:
\[
\langle \dot\gamma^t(u) + f^t(v) \dot\xi^t(u), \dot f^t(v) \xi^t(u) \rangle = \langle \dot\gamma(u) + f(v) \dot\xi(u), \dot f(v) \xi(u) \rangle.
\]
Since $\dot\xi^t(u) \parallel \dot\gamma^t(u)$ and due to \eqref{eqn:GammaConst} and \eqref{eqn:FXiConst} this is equivalent to
\[
\langle \dot\xi^t(u), \xi^t(u) \rangle = \langle \dot\xi(u), \xi(u) \rangle,
\]
which means that $\frac{d}{du}\left(\|\xi^t(u)\|^2\right) = \frac{d}{du}\left(\|\xi(u)\|^2\right)$, in other words, that the difference $\|\xi^t(u)\|^2 - \|\xi(u)\|^2$ is a function of $t$ only, independent of~$u$.
By a parameter change one may assume
\begin{equation}
\label{eqn:NormXi}
\|\xi^t(u)\|^2 = \|\xi(u)\|^2 + t \quad \text{for all }t.
\end{equation}
Reversing the above argument one sees that, under assumption of \eqref{eqn:GammaConst} and \eqref{eqn:FXiConst}, equation \eqref{eqn:NormXi} is not only necessary, but also sufficient for $F^t = F$.

Finally, the equation $G^t = G$ due to \eqref{eqn:FXiConst} and \eqref{eqn:NormXi} takes the form
\[
(\dot f(v))^2 (\|\xi(u)\|^2 + t) + (\dot z^t(v))^2 =
(\dot f(v))^2 \|\xi(u)\|^2 + (\dot z(v))^2,
\]
which solves for $\dot z^t(v)$ as
\[
(\dot z^t(v))^2 = (\dot z(v))^2 - t(\dot f(v))^2
\]
and under the assumption $\dot z > 0$ leads to the formula in the theorem.

It remains to derive the formulas for $\gamma^t$ and $\xi^t$ satisfying the conditions \eqref{eqn:GammaConst}, \eqref{eqn:FXiConst}, \eqref{eqn:NormXi}, as well as $\dot\gamma^t(u) \parallel \dot\xi^t(u)$.
Since we have
\[
\dot\gamma^t(u) = g^t(u)
\begin{pmatrix}
-\sin\psi^t(u)\\ \cos\psi^t(u)
\end{pmatrix}, \quad
\xi^t(u) = c^t(u)
\begin{pmatrix}
\cos\phi^t(u)\\ \sin\phi^t(u)
\end{pmatrix},
\]
equations \eqref{eqn:GammaConst} and \eqref{eqn:NormXi} simply say $g^t = g$ and $c^t = \sqrt{c^2 + t}$.
The second of the equations \eqref{eqn:FXiConst} was solved when we studied the deformation of axial T-surfaces, see \eqref{eqn:CTAxial} and thereafter.
This leads to the formula for $\phi^t$.
Equation \eqref{eqn:CEta} applied to $\sigma^t$ implies
\[
\tan\eta^t = \frac{\dot c^t}{c^t \dot\phi^t} = \frac{\dot c}{\dot\phi \sqrt{c^2 + \frac{t}{\cos^2\eta}}} = \frac{c \sin\eta}{\sqrt{c^2\cos^2\eta + t}},
\]
and the theorem is proved.
\end{proof}

\subsection{Parallel pairs of T-surfaces}
\begin{definition}
Two T-surfaces $\sigma, \sigma' \colon U \times V \to \R^3$ are called \emph{parallel} if for every $(u,v) \in U \times V$ one has $\sigma_u(u,v) \parallel \sigma'_u(u,v)$ and $\sigma_v(u,v) \parallel \sigma'_v(u,v)$.
\end{definition}

\begin{theorem}
Let $\sigma$ be a T-surface given by \eqref{eqn:TSurfParam}.
Then there are infinitely many T-surfaces parallel to $\sigma$.
Each of them is generated by a pair of curves, one parallel to the trajectory curve $\sigma|_{v=0}$ and one parallel to the profile curve $\sigma|_{u=0}$.
\end{theorem}
\begin{proof}
One has
\[
\sigma(u,0) =
\begin{pmatrix}
\gamma(u)\\ z(0)
\end{pmatrix}, \quad
\sigma(0,v) =
\begin{pmatrix}
\gamma(0) + f(v)\xi(0)\\ z(v).
\end{pmatrix}
\]
If $\sigma'$ is parallel to $\sigma$, then its initial trajectory curve must be given by a curve $\gamma'$ parallel to $\gamma$, and its initial profile curve by a curve $(\gamma'(0) + f'(v)\xi(0), z'(v))$ parallel to $\sigma_{u=0}$.
On the other hand, for any $\gamma'$, $f'$, $z'$ satisfying these conditions the surface
\[
\sigma'(u,v) =
\begin{pmatrix}
\gamma'(u) + f'(v)\xi(u), z'(v)
\end{pmatrix}
\]
(provided that it is regular, for which one has to take care of only the last part of the sixth condition in Theorem \ref{thm:DescrTSurf}) is parallel to $\sigma$ and is a T-surface.
Both facts follow from $\dot\xi \parallel \dot\gamma \parallel \dot\gamma'$.
\end{proof}

\begin{corollary}
For every T-surface with $\dot\phi \ne 0$ there is a parallel axial T-surface.
\end{corollary}
\begin{proof}
It suffices to put $\gamma' = \xi$.
\end{proof}

\begin{theorem}
If two T-surfaces are parallel, then their isometric deformations remain parallel.
\end{theorem}
\begin{proof}
On one hand this is a special case of a more general theorem: if smooth surfaces $\sigma$ and $\sigma'$ have parallel conjugate nets, and $\sigma$ admits an isometric deformation which preserves its net, then $\sigma'$ also admits an isometric deformation which preserves its net and keeps it parallel to the net on $\sigma$, see \cite[\S 10]{sauer1970differenzengeometrie}.

On the other hand, this follows from our formulas for isometric deformations.
By the previous theorem it suffices to keep track of the parallelity of the initial profile and trajectory curves.
For the $0$-profile curves of $\sigma$ and $\sigma'$ the formulas for $f^t$ and $z^t$ ensure that they remain parallel.
For the $0$-trajectory curve these are the formulas for $\psi$: they depend only on the vector field $\xi$, which is common for $\sigma$ and $\sigma'$.
\end{proof}
\section*{Acknowledgments}
The first and the second authors are supported by grant F77 (SFB “Advanced Computational Design”, subproject SP7) of the Austrian Science Fund FWF.
Additionally, the authors wish to thank Georg Nawratil and Kiumars Sharifmoghaddam for helpful discussions.

\end{document}